\documentclass[12pt]{article}
\usepackage[totalwidth=480pt, totalheight=680pt]{geometry}
\usepackage{amsmath, amssymb, amsthm, amsfonts, mathrsfs, eucal, enumerate, rotating, lscape,layout}
\usepackage[all,tips,2cell]{xy}
\usepackage[usenames]{color}
\usepackage{xspace}
\usepackage{ifpdf}
\usepackage{ifthen}

\newtheorem{theorem}{Theorem}[section]
\newtheorem*{theorem*}{Theorem}
\newtheorem{corollary}[theorem]{Corollary}
\newtheorem*{corollary*}{Corollary}
\newtheorem{lemma}[theorem]{Lemma}

\newtheorem*{claim*}{Claim}
\newtheorem*{lemma*}{Lemma}
\newtheorem{proposition}[theorem]{Proposition}
\newtheorem{def-prop}[theorem]{Definition-Proposition}
\newtheorem*{proposition*}{Proposition}
\theoremstyle{remark}
\newtheorem{remark}[theorem]{Remark}
\theoremstyle{definition}
\newtheorem{definition}[theorem]{Definition}

\numberwithin{equation}{section}
\newtheorem{intro_thm}[subsubsection]{Theorem}
\newtheorem*{intro_prop*}{Proposition}

\newcommand{\Ob}{\mathrm{Ob}}
\newcommand{\id}{\mathrm{id}}
\newcommand{\pr}{\mathrm{pr}}
\newcommand{\alman}[1]{\mathfrak{#1}}
\newcommand{\gotika}{\alman a}

\newcommand{\gotikl}{\alman l}
\newcommand{\gotikr}{\alman r}

\newcommand{\gotiki}{\alman i}
\newcommand{\TORS}{\textsc{Tors}}
\newcommand{\twoCh}{\textsc{2Pic}}
\newcommand{\oneCh}{\textsc{Pic}}
\newcommand{\twoch}{2\wp}
\newcommand{\onech}{\wp}
\newcommand{\Ch}{2\mathfrak{S}\textsc{tack}}
\newcommand{\ch}{2\textsc{Stack}}
\newcommand{\PCh}{2\mathfrak{P}\textsc{ic}}
\newcommand{\boyut}[1]{\mathscr{#1}}
\newcommand{\Frac}{\mathsf{Frac}}
\newcommand{\Hom}{\mathrm{Hom}}
\newcommand{\Homsf}{\mathsf{Hom}}
\newcommand{\FFrac}{\mathbb{F}\mathsf{rac}}
\newcommand{\HHom}{\mathbb{H}\mathsf{om}}
\newcommand{\rmD}{\mathrm{D}}
\newcommand{\rmC}{\mathrm{C}}
\newcommand{\rmT}{\mathrm{T}}
\newcommand{\im}{\mathrm{im}}
\newcommand{\coker}{\mathrm{coker}}
\newcommand{\oneA}{\boyut A}
\newcommand{\oneB}{\boyut B}

\newcommand{\oneE}{\boyut E}
\newcommand{\oneK}{\boyut K}

\newcommand{\oneL}{\boyut L}

\newcommand{\oneP}{\boyut P}

\newcommand{\ES}{\mathsf{S}}

\newcommand{\twoC}{\mathbb{C}}

\newcommand{\twoP}{\mathbb{P}}
\newcommand{\twoQ}{\mathbb{Q}}
\DeclareMathOperator{\twosk}{2sk}
\newcommand{\ra}{{\rightarrow}}
\newcommand{\Ra}{{\Rightarrow}}

\newcommand{\Da}{{\Downarrow}}

\newcommand{\kesir}[2]{\genfrac{}{}{0pt}{}{#1}{#2}}
\xyoption{2cell}
\xyoption{rotate}
\usepackage{hyperref}
\begin{document}
\UseAllTwocells
\title{Length 3 Complexes of Abelian Sheaves and Picard 2-Stacks}
\author{A. Emin Tatar\\
  \small Department of Mathematics, Florida State University\\
  \small Tallahassee, FL 32306-4510, USA\\
  \small \texttt{atatar@math.fsu.edu}}
\date{}
\maketitle
\begin{abstract}
We define a tricategory $\rmT^{[-2,0]}$ of length 3 complexes of abelian sheaves, whose hom-bigroupoids consist of weak morphisms of such complexes. We also define a 3-category $\twoCh(\ES)$ of Picard 2-stacks, whose hom-2-groupoids consist of additive 2-functors. We prove that these categories are triequivalent as tricategories. As a consequence we obtain a generalization of Deligne's analogous result about Picard stacks in SGA4, Expos\'{e} XVIII(Deligne (1973) \cite{Deligne}).
\end{abstract}
\tableofcontents
\section*{Introduction}\label{section:introduction}
\addcontentsline{toc}{section}{Introduction}
Let $\rmD^{[-1,0]}(\ES)$ be the subcategory of the derived category of category of complexes of abelian sheaves $A^{\bullet}$ over a site $\ES$ with $H^{-i}(A^{\bullet}) \neq 0$ only for $i=0,1$. Let $\oneCh^{\flat}(\ES)$ denote the category of Picard stacks over $\ES$ with 1-morphisms isomorphism classes of additive functors. In SGA4 Expos\'{e} XVIII, Deligne shows the following.

\begin{intro_prop*}\cite[Proposition 1.4.15]{Deligne}
The  functor
$$\onech^{\flat}:\xymatrix@1{\rmD^{[-1,0]}(\ES) \ar[r] & \oneCh^{\flat}(\ES)}$$
given by sending a length 2 complex of abelian sheaves, $A^{\bullet}:A^{-1} \ra A^0$ over $\ES$ to its associated Picard stack $[A^{-1} \ra A^0]^{\sim}$, an isomorphism class of fractions from $A^{\bullet}$ to $B^{\bullet}$ to an isomorphism class of morphisms of associated Picard stacks is an equivalence.
\end{intro_prop*}
The purpose of this paper is to generalize the above result to Picard 2-stacks over $\ES$.
Let $\twoCh^{\flat \flat}(\ES)$ denote the category of Picard 2-stacks, whose morphisms are equivalence classes of additive 2-functors. Let $\rmD^{[-2,0]}(\ES)$ be the subcategory of the derived category of category of complexes of abelian sheaves $A^{\bullet}$ over $\ES$ with $H^{-i}(A^{\bullet}) \neq 0$ for $i=0,1,2$.

\begin{intro_thm}\label{theorem:intro_corollary}
The  functor

$$\twoch^{\flat \flat}: \xymatrix@1{\rmD^{[-2,0]}(\ES) \ar[r] & \twoCh^{\flat \flat}(\ES)}$$
given by sending a length 3 complex of abelian sheaves, $A^{\bullet}:A^{-2} \ra A^{-1} \ra A^0$ over $\ES$ to its associated Picard 2-stack $[A^{-2} \ra A^{-1} \ra A^0]^{\sim}$, an equivalence class of fractions from $A^{\bullet}$ to $B^{\bullet}$ to an equivalence class of morphisms of associated Picard 2-stacks is an equivalence.
\end{intro_thm}

Basically, it gives a geometric description of the derived category of length 3 complexes of abelian sheaves. It states that any Picard 2-stack over a site $\ES$ is biequivalent to a Picard 2-stack associated to a length 3 complex of abelian sheaves and that any morphism of Picard 2-stacks comes from a fraction of such complexes.  A complex of abelian sheaves, whose only non-zero cohomology groups are placed at degrees -2,-1, and 0 can be thought as a length 3 complex of abelian sheaves, and therefore a morphism in $\rmD^{[-2,0]}(\ES)$ between any two complexes $A^{\bullet}$ and $B^{\bullet}$ is given by an equivalence class of fraction

$$(q,M^{\bullet},p):\xymatrix@1{A^{\bullet} & M^{\bullet} \ar[r]^p \ar[l]_q & B^{\bullet}}$$
with $q$ being a quasi-isomorphism.

However, we prove a much stronger statement, so that the latter theorem becomes an immediate consequence of it. Let $\twoCh(\ES)$ be the 3-category of Picard 2-stacks where 1-morphisms are additive 2-functors, 2-morphisms are natural 2-transformations, and 3-morphisms are modifications. Length 3 complexes of abelian sheaves over $\ES$ placed in degrees $[-2,0]$ form a 3-category $\rmC^{[-2,0]}(\ES)$ by adding to the regular morphisms of complexes, the degree -1 and -2 morphisms. Then we easily construct an explicit trihomomorphism

\begin{equation*}
\twoch:\xymatrix@1{\rmC^{[-2,0]}(\ES) \ar[r] & \twoCh(\ES)},
\end{equation*}
that is a 3-functor between 3-categories. Under this construction, length 3 complexes of abelian sheaves correspond to Picard 2-stacks. Although morphisms of such complexes induce morphisms between associated Picard 2-stacks, not all of them are obtained in this way. In this sense, the 1-morphisms of $\rmC^{[-2,0]}(\ES)$ are not geometric and the reason is their strictness. We resolve this problem by weakening $\rmC^{[-2,0]}(\ES)$ as follows: We introduce a tricategory $\rmT^{[-2,0]}(\ES)$ (a tricategory is a weak version of a 3-category in the sense of \cite{MR1261589}) with same objects as $\rmC^{[-2,0]}(\ES)$. For any two complexes of abelian sheaves $A^{\bullet}$ and $B^{\bullet}$, morphisms between $A^{\bullet}$ and $B^{\bullet}$ in $\rmT^{[-2,0]}(\ES)$ is the bigroupoid $\Frac(A^{\bullet},B^{\bullet})$, whose main property is that it satisfies $\pi_0(\Frac(A^{\bullet},B^{\bullet})) \simeq \Hom_{\rmD^{[-2,0]}(\ES)}(A^{\bullet},B^{\bullet})$, where $\pi_0$ denotes the isomorphism classes of objects. Roughly speaking, objects of $\Frac(A^{\bullet},B^{\bullet})$ are fractions from $A^{\bullet}$ to $B^{\bullet}$ in the ordinary sense and its 2-morphisms are certain commutative diagrams (\ref{diagram:diamond}) called ``\emph{diamonds}''. Then we prove:

\begin{intro_thm}\label{theorem:intro_main_thm}
The trihomomorphism

\begin{equation*}
\twoch:\xymatrix@1{\rmT^{[-2,0]}(\ES) \ar[r] & \twoCh(\ES)}
\end{equation*}
defined by sending $A^{\bullet}$ a length 3 complex of abelian sheaves to its associated Picard 2-stack is a triequivalence.
\end{intro_thm}
Since in particular a triequivalence is essentially surjective, every Picard 2-stack is biequivalent to a Picard 2-stack associated to a complex of abelian sheaves. Then by ignoring the 3-morphisms and passing to the equivalence class of morphisms in the triequivalence of Theorem \ref{theorem:intro_main_thm}, we deduce Theorem \ref{theorem:intro_corollary}.

\subsection*{Organization of the paper}
This paper is organized as follows:

In Section \ref{section:preliminaries}, we recall important points of butterflies in abelian context, $(A,B)$-torsors, where $A$ and $B$ are abelian sheaves, and $(\oneA,\oneB)$-torsors, where $\oneA$ and $\oneB$ are Picard stacks. We also remind the reader some important results from \cite{Aldrovandi2009687} that we will refer continuously.

In Section \ref{section:torsors}, we explain briefly the basics on 2-stacks with structures, and exact sequences of Picard 2-stacks. We also give an example of Picard 2-stack, namely $\TORS(\oneA,A^0)$, where $\oneA$ is a Picard stack and $A^0$ is an abelian sheaf. This example will be of great importance for the rest since it will be the Picard 2-stack associated to $A^{\bullet}$ a length 3 complex of abelian sheaves. We define at the end of the section the 3-category $\twoCh(\ES)$ of Picard 2-stacks, as well.

In Section \ref{section:3-category_of_complexes}, we first define the 3-category $\rmC^{[-2,0]}(\ES)$ of length 3 complexes of abelian sheaves. We construct an explicit trihomomorphism from this 3-category to the 3-category of Picard 2-stacks.

In Section \ref{section:fractions}, for any two length 3 complexes of abelian sheaves $A^{\bullet}$ and $B^{\bullet}$, we define a bigroupoid $\Frac(A^{\bullet},B^{\bullet})$. It is a weakened version of the hom-2-category $\Hom_{\rmC^{[-2,0]}}(A^{\bullet},B^{\bullet})$ in the sense that $\pi_0(\Frac(A^{\bullet},B^{\bullet})) \simeq \Hom_{\rmD^{[-2,0]}(\ES)}(A^{\bullet},B^{\bullet})$.

In Section \ref{section:bieqv_of_hom_bicats}, we show that for any two length 3 complexes of abelian sheaves $A^{\bullet}$ and $B^{\bullet}$, there exists a biequivalence as bigroupoids between $\Frac(A^{\bullet},B^{\bullet})$ and the 2-category $\Homsf(A^{\bullet},B^{\bullet})$ of morphisms of Picard 2-stacks associated to $A^{\bullet}$ and $B^{\bullet}$.

In Section \ref{section:tricategory_of_complexes}, we define the tricategory $\rmT^{[-2,0]}(\ES)$. It consists of same objects as $\rmC^{[-2,0]}(\ES)$ and for any two length 3 complexes $A^{\bullet}$ and $B^{\bullet}$ of abelian sheaves, $\Frac(A^{\bullet},B^{\bullet})$ as the hom-bigroupoid. We extend the trihomomorphism constructed in Section \ref{section:3-category_of_complexes} to a trihomomorphism on $\rmT^{[-2,0]}(\ES)$. We prove Theorem \ref{theorem:intro_main_thm} which says that the latter trihomomorphism is a triequivalence and from which Theorem \ref{theorem:intro_corollary} follows.

In Section \ref{section:stackification}, we informally discuss the stack versions of what has been done in the previous sections.

\subsection*{Acknowledgements}
I would like to express my profound gratitude to my advisor, Ettore Aldrovandi, for helping me at all stages of the paper, which is part of my Ph.D.\ thesis. I would like to thank Behrang Noohi for helpful conversations. I also thank Chris Portwood for proofreading.

\section{Preliminaries}\label{section:preliminaries}
The method that we are going to adopt to prove our results is going to use mostly the language and techniques developed in the papers of Aldrovandi and Noohi such as butterflies, torsors, etc. So it is worthwhile to mention here some of their work. We finish with a few words about bicategories and tricategories. Before, let us fix our conventions and notations.

Throughout the paper, we will work with sheaves, stacks, etc.\, defined over a site $\ES$. For simplicity, we will assume that $\ES$ has fibered products. Fibered 2-categories, 2-functors, and natural 2-transformations will be used in the sense of Hakim \cite{MR0364245}. A complex of abelian sheaves will mean a length 3 complex of abelian sheaves over the site $\ES$ unless otherwise stated. It will be denoted as

$$A^{\bullet}:\xymatrix@1{A^{-2} \ar[r]^{\delta_A} & A^{-1} \ar[r]^{\lambda_A} & A^0}.$$
For any complex of abelian sheaves $A^{\bullet}$, $A^{\bullet <0}$ denotes the complex

$$A^{\bullet <0}:\xymatrix@1{A^{-2} \ar[r]^{\delta_A} & A^{-1} \ar[r] & 0}$$
and therefore $f^{\bullet <0}:A^{\bullet <0} \ra B^{\bullet <0}$ a morphism of complexes between $A^{\bullet <0}$ and $B^{\bullet <0}$.

\subsection{Butterflies} The reader can refer to \cite{noohi-2005} and \cite{MR2280287} for details of butterflies over a point or to \cite{Aldrovandi2009687} for a treatment over a site. Here, we will remind the basic definitions following the latter point of view in an abelian context.
\begin{definition}\label{definition:butterfly}
Let $A^{\bullet}:A^{-1} \ra A^0$ and $B^{\bullet}:B^{-1} \ra B^0$ be two length 2 complexes of abelian sheaves. A butterfly from $A^{\bullet}$ to $B^{\bullet}$ is a commutative diagram of abelian sheaf morphisms of the form

\begin{equation}\label{diagram:butterfly}
\xymatrix{A^{-1} \ar[dd] \ar[dr]^{\kappa} && B^{-1} \ar[dd] \ar[dl]_{\imath}\\
            & E \ar[dl]_{\rho} \ar[dr]^{\jmath} &\\
            A^0 && B^0}
\end{equation}
where $E$ is an abelian sheaf, the NW-SE sequence is a complex, and the NE-SW sequence is an extension. $[A^{\bullet},E,B^{\bullet}]$ will denote a butterfly from $A^{\bullet}$ to $B^{\bullet}$.
\end{definition}

A morphism of butterflies $\varphi:[A^{\bullet},E,B^{\bullet}] \ra [A^{\bullet},E',B^{\bullet}]$ is an abelian sheaf isomorphism $E \ra E'$ satisfying certain commutative diagrams. Two such morphisms compose in an obvious way. Therefore butterflies from $A^{\bullet}$ to $B^{\bullet}$ form a groupoid denoted by $\textsf{B}(A^{\bullet},B^{\bullet})$. A butterfly is \emph{flippable} or \emph{reversible} if both diagonals of (\ref{diagram:butterfly}) are extensions.

For more about crossed modules and butterflies in the abelian context, we refer the reader to \cite[\S 12]{noohi-2005} and \cite[\S8]{Aldrovandi2009687}.

\subsection{$(A,B)$-torsors}
Let $A \ra B$ be a morphism of, not necessarily abelian, sheaves. An $(A,B)$-torsor is a pair $(L,x)$, where $L$ is an $A$-torsor and $x:L \ra B$ is an $A$-equivariant morphism of sheaves (see \cite{MR546620}). A morphism between two pairs $(L,x)$ and $(K,y)$ is a morphism of sheaves $F: L \ra K$ compatible with the action of $A$ such that the diagram

$$\xymatrix{L  \ar[ddr]_x \ar[rr]^{F} &\ddtwocell<\omit>{=}& K  \ar[ddl]^y\\
            &&&\\
            &B&}$$
commutes. $(A,B)$-torsors form a category denoted by $\TORS(A,B)$.

\subsection{$(\oneA,\oneB)$-torsors}
Let $\oneA$ be a gr-stack, not necessarily Picard. A stack $\oneP$ in groupoids is a (right) $\oneA$-torsor if there exists a morphism of stacks

$$\textrm{m}: \xymatrix@1{\oneP \times \oneA \ar[r] & \oneP}$$
compatible with the group laws in $\oneA$, and the morphism

$$(\pr,\textrm{m}):\xymatrix@1{\oneP \times \oneA \ar[r] & \oneP \times \oneP}$$
is an equivalence, and for all $U \in \ES$, $\oneP_U$ is not empty. \cite[\S 6.1]{MR1086889}

Let $\oneA \ra \oneB$ be a morphism of gr-stacks. An $(\oneA,\oneB)$-torsor is a pair $(\oneL,x)$, where $\oneL$ is an $\oneA$-torsor, and $x:\oneL \ra \oneB$ is an $\oneA$-equivariant morphism of stacks \cite[\S 6.1]{MR2387582}, \cite[\S 6.3.4]{Aldrovandi2009687}. A 1-morphism of $(\oneA,\oneB)$-torsors is a pair

$$(F,\mu):\xymatrix@1{(\oneL,x) \ar[r] & (\oneK,y)},$$
where $F: \oneL \ra \oneK$ is a morphism of stacks such that

$$\xymatrix{\oneL \ar[rr]^F \ar[ddr]_x \ddrrtwocell<\omit>{\mspace{20mu}\sigma_F}&& \oneK \ar[ddl]^y\\
            &&\\
            & \oneB &}$$
and $\mu$ is a natural transformation of stacks

$$\xymatrix{\oneL \times \oneA \ar[dd] \ar[rr]^{F \times 1} \ddrrtwocell<\omit>{\mu}&& \oneK \times \oneA \ar[dd]\\
            &&&\\
            \oneL \ar[rr]_F && \oneK}$$
expressing the compatibility of $F$ with the torsor structure. A 2-morphism of $(\oneA,\oneB)$-torsors $(F,\mu) \Ra (G,\nu)$ is given by a natural transformation $\phi:F \Ra G$ satisfying the conditions given in \cite[\S6.3.4]{Aldrovandi2009687}. $(\oneA,\oneB)$-torsors form a 2-stack denoted by $\TORS(\oneA, \oneB)$.

\subsection{Abelian Sheaves and Picard Stacks}
We recall Deligne's work about abelian sheaves and Picard stacks from \cite[\S 1.4]{Deligne}. They are going to be referred \textit{sans cesse} throughout the paper. These results are also revisited by Aldrovandi and Noohi in \cite{Aldrovandi2009687}. In order to be consistent with the rest of the paper, we recall them as they are announced in \cite{Aldrovandi2009687}.

\begin{theorem*}\cite[Theorem 8.3.1]{Aldrovandi2009687}
For any two length 2 complexes of abelian sheaves $A^{\bullet}$ and $B^{\bullet}$, there is an equivalence of groupoids

$$\xymatrix@1{\Homsf(A^{\bullet},B^{\bullet}) \ar[r]^{\sim} & \mathsf{B}(A^{\bullet},B^{\bullet})},$$
where $\Homsf(A^{\bullet},B^{\bullet})$ is the groupoid of additive functors between the Picard stacks associated to $A^{\bullet}$ and $B^{\bullet}$.
\end{theorem*}

\begin{proposition*}\cite[Proposition 8.3.2]{Aldrovandi2009687}
Let $\oneA$ be a Picard stack. Then there exists a length 2 complex of abelian sheaves $A^{\bullet}:A^{-1} \ra A^0$ such that $\oneA$ is equivalent to Picard stack $[A^{-1} \ra A^0]^{\sim}$.
\end{proposition*}

Let $\rmC^{[-1,0]}(\ES)$ denote the bicategory of morphisms of abelian sheaves over $\ES$ with commutative squares as 1-morphisms and homotopies as 2-morphisms. Let $\oneCh(\ES)$ denote the 2-category of Picard stacks over $\ES$ with 1-morphisms being additive functors and 2-morphisms being natural 2-transformations. Putting the above results together, Deligne proves:

\begin{theorem*}\cite[Proposition 8.4.3]{Aldrovandi2009687}
The functor

$$\xymatrix@1{\rmC^{[-1,0]}(\ES) \ar[r] & \oneCh(\ES)}$$
defined by sending a morphism of abelian sheaves $A^{\bullet}:[A_1 \ra A_2]$ to its associated Picard stack $[A_1 \ra A_2]^{\sim}$ is a biequivalence of bicategories.
\end{theorem*}

\begin{remark}
In the same paper, the authors also prove these facts in the non-abelian context by not assuming that stacks and sheaves are necessarily Picard or abelian.
\end{remark}

\subsection{Tricategories}
Even though the language of bicategories and tricategories is going to be extensively used, we are not going to remind here in full detail bicategories or tricategories. Just for motivation, a 3-category can be thought as the category of 2-categories with 2-functors or weak 2-functors in the sense of B{\'e}nabou \cite{MR0220789} and a tricategory as a weakened version of a 3-category.  However, we want to recall the triequivalence since the proof of Theorem \ref{theorem:triequivalence} will follow its definition. For more about bicategories and tricategories, we refer the reader to \cite{MR1261589,Gurski_thesis,MR2276246,MR0220789}.
\begin{definition}\label{definition:triequivalence}\cite{MR2276246}
A trihomomorphism of tricategories $T: \mathfrak{C} \ra \mathfrak{D}$ is called a triequivalence if it induces biequivalences $T_{X,Y}: \mathfrak{C}(X,Y) \ra \mathfrak{D}(TX,TY)$ of hom-bicategories for all objects $X,Y$ in $\mathfrak{C}$ ($T$ is locally a biequivalence), and every object in $\mathfrak{D}$ is biequivalent in $\mathfrak{D}$ to an object of the form $TX$ where $X$ is an object in $\mathfrak{C}$.
\end{definition}

\section{Picard 2-Stacks as Torsors}\label{section:torsors}
In this section, our goal is to give some of the fundamental facts about 2-stacks and torsors that will be needed throughout the paper. Our main references for 2-stacks with structures such as monoidal, group-like, braided, Picard are \cite{MR1191733,MR1301844} and for torsors \cite{MR2387582,MR1086889}
\subsection{2-Stacks}
\begin{definition}\cite[Definition 6.2]{breen-2006}
A 2-stack $\twoP$ is a fibered 2-category in 2-groupoids such that
\begin{itemize}
\item for all $X,Y$ objects in $\ES_U$, $\Hom_{\twoP_U}(X,Y)$ is a stack over $\ES/U$;
\item 2-descent data is effective for every object in $\twoP$.
\end{itemize}
\end{definition}

In the above definition 2-groupoids are considered in the sense of Breen \cite{MR1301844}, that is, 1-morphisms are weakly invertible.

\begin{definition}\label{definition:gr_2_stack}\cite[Definition 8.4]{MR1301844}
A gr-2-stack $\twoP$ is a 2-stack with a morphism $\otimes: \twoP \times \twoP \ra \twoP$ of 2-stacks, an associativity constraint $\gotika$ compatible with $\otimes$, a left unit $\gotikl$ and a right unit $\gotikr$ constraints compatible with $\gotika$, and an inverse constraint $\gotiki$ with respect to $\otimes$ compatible with units.
\end{definition}

A more detailed definition of gr-2-stack can be found in \cite{MR1191733}. Next, following \cite[\S8.4]{MR1301844} we add to gr-2-stacks commutativity constraints with an increasing level of strictness.

\begin{definition}
A gr-2-stack $\twoP$ is said to be:
\begin{itemize}
\item \emph{braided}, if there exists a functorial natural transformation
        \[R_{X,Y}:\xymatrix@1{X \otimes Y \ar[r] & Y \otimes X}\]
        that satisfy the 2-braiding axioms of Kapranov and Voevodsky  \cite{MR1278735} together with the additional condition that, in their terminology, the pair of 2-morphisms defining the induced $Z$-systems coincide. The corrected and full 2-braiding axioms can be found in \cite{MR1402727}.
\item \emph{strongly braided}, if it is braided and for any $X,Y$ two objects, there exists a functorial 2-morphism
        \begin{equation}\label{diagram:strongly_braided}
        \xymatrix{X \otimes Y \rrtwocell<5>^{1_{X \otimes Y}}_{R_{Y,X}R_{X,Y}}{\mspace{36mu}S_{X,Y}} && X \otimes Y}.
        \end{equation}
        such that the two compatibility conditions given in \cite[Definition 15]{MR1458415} are satisfied.
\item \emph{symmetric}, if it is strongly braided and the following whiskerings coincide:
        \begin{equation}\label{diagram:2_symmetric1}
        \xymatrix {X \otimes Y \rrtwocell<5>^{1_{X \otimes Y}}_{R_{Y,X}R_{X,Y}}{\mspace{36mu}S_{X,Y}} && X \otimes Y \ar[r]^{R_{X,Y}} & Y \otimes X},
        \end{equation}
        \begin{equation}\label{diagram:2_symmetric2}
        \xymatrix {X \otimes Y \ar[r]^{R_{X,Y}} & Y \otimes X \rrtwocell<5>^{1_{Y,X}}_{R_{X,Y}R_{Y,X}}{\mspace{36mu}S_{Y,X}} && Y \otimes X}.
        \end{equation}
\item \emph{Picard}, if it is symmetric and for any object $X$, there exists a functorial 2-morphism
        \begin{equation}{\label{diagram:2_picard}}
        \xymatrix{X \otimes X \rrtwocell^{1_{X,X}}_{R_{X,X}}{\mspace{28mu}S_X} && X \otimes X}
        \end{equation}
        additive in $X$ (i.e. there is a relation between $S_{X\otimes Y}$, $S_X$, and $S_Y$) such that $S_{X,X}=S_X * S_X$.
\end{itemize}
\end{definition}
Further down in the paper, we will be talking about the 3-category of Picard 2-stacks which requires the concept of morphism of Picard 2-stacks. Following Breen \cite{MR1301844}, we will call such a morphism additive 2-functor. It will be a cartesian 2-functor between the underlying fibered 2-categories compatible with the monoidal, braided, and Picard structures carried by the 2-categories. The compatibility with monoidal structure is already known. In Gordon, Power, Street \cite{MR1261589}, a monoidal 2-category is defined as a one-object tricategory. More in detail, one can think of a monoidal 2-category as the hom-2-category of a one-object tricategory, whose associativity and unit constraints hold up to 2-isomorphisms and whose modifications are invertible. Then the trihomomorphism \cite[Definition 3.1]{MR1261589} between such tricategories will be the right definition of morphism between monoidal 2-categories. For the compatibility with the rest of the structures, we refer the reader to the author's thesis \cite{Tatar_thesis}.

Here is a technical result that we will use several times in our proofs.

\begin{lemma}\label{lemma:pullback}
Let $\twoP$ be a Picard 2-stack and $A,B$ be two abelian sheaves with additive 2-functors
$\phi:\xymatrix{A \ar[r] & \twoP}$ and $\psi: \xymatrix{B \ar[r] & \twoP}$. Then $A \times_{\twoP} B$ is a Picard stack.
\end{lemma}

\begin{proof}
The fibered category $A \times_{\twoP} B$ with fibers $(A \times_{\twoP} B)_{|U}$ consisting of
\begin{itemize}
\item objects $(a,f,b)$, where $a \in A(U)$, $b \in B(U)$, and $f: \phi(a) \ra \psi(b)$ is a 1-morphism in $\twoP_U$;
\item morphisms $(a,f,\alpha,g,b)$, where $\xymatrix{\phi(a) \rrtwocell^{f}_{g}{\alpha} && \psi(b)}$ is a 2-morphism in $\twoP_U$;
\end{itemize}
is a prestack since for any $U \in \ES$, 1-morphisms of $\twoP$ form a stack over $\ES/U$. It is in fact a stack.

Let $((U_i \ra U),(a_i,f_i,b_i),\alpha_{i,j})_{i,j \in I}$ be a descent datum with $(U_i \ra U)_{i \in I}$ a covering of $U$, $(a_i,f_i,b_i)$ an object in $(A \times_{\twoP} B)_{U_i}$ and $\alpha_{i,j}$ a 1-morphism in $(A \times_{\twoP} B)_{U_{ij}}$ between $(a_j,f_j,b_j)_{|U_{ij}}$ and $(a_i,f_i,b_i)_{|U_{ij}}$. Since $a_{i|U_{ij}}=a_{j|U_{ij}}$, $b_{i|U_{ij}}=b_{j|U_{ij}}$ and $A$ and $B$ are sheaves, there exist $a \in A(U)$ and $b \in B(U)$ such that $a_{|U_i} = a_i$ and $b_{|U_i} = b_i$. Then the collection $((U_i \ra U),f_i,\alpha_{i,j})_{i,j \in I}$ satisfies the descent in $\Hom(\phi(a),\psi(b))$, which is effective since $\twoP$ is a Picard 2-stack. That is, there exists $f \in \Hom(\phi(a),\psi(b))$ and $\beta_i: f_{|U_i} \Ra f_i$ compatible with $\alpha_{i,j}$ such that for all $i \in I$, $(a_i,f_{|U_i},\beta_i,f_i,b_i)$ is a morphism from $(a,f,b)_{|U_i}$ to $(a_i,f_i,b_i)$. Thus, the descent $((U_i \ra U),(a_i,f_i,b_i),\alpha_{i,j})_{i,j \in I}$ is effective.

Next, we show that $A \times_{\twoP} B$ is Picard. First, let us recall the notation from Definition $\ref{definition:gr_2_stack}$. $\otimes_{\twoP}$ is the monoidal operation, $\gotika$, $\gotikl$, $\gotikr$, $\gotiki$, $R_{-,-}$, $S_{-,-}$, and $S_{-}$ are respectively associativity, left unit, right unit, inverse, braiding, symmetry, and Picard constraints. The unnamed arrows in the diagrams below are structural equivalences resulting from additive 2-functors $\phi$ and $\psi$.
\begin{itemize}
\item[] \emph{monoidal structure}: The multiplication is defined as

    $$(a_1,f_1,b_1) \otimes (a_2,f_2,b_2):=(a_1+a_2,f_1f_2,b_1+b_2),$$
    where $f_1f_2$ is the morphism that makes the diagram
    $$\xymatrix{\phi(a_1) \otimes_{\twoP} \phi(a_2) \ar[dd] \ar[rr]^{f_1 \otimes_{\twoP} f_2} && \psi(b_1) \otimes_{\twoP} \psi(b_2) \ar[dd]\\
                &&&\\
                \phi(a_1+a_2) \ar[rr]_{f_1f_2} && \psi(b_1+b_2) \uulltwocell<\omit>{\mspace{-20mu}N_m}}$$
    commute up to a 2-isomorphism $N_m$.

    For any three objects $(a_i,f_i,b_i)$ for $i=1,2,3$, the associator is given by the morphism
    $(a_1+a_2+a_3,f_1(f_2f_3),\alpha_{f_1,f_2,f_3},(f_1f_2)f_3,b_1+b_2+b_3)$,
    where $\alpha_{f_1,f_2,f_3}$ is defined as the 2-isomorphism of the bottom face that makes the following cube commutative (we ignored $\otimes_{\twoP}$ for compactness).

    \footnotesize
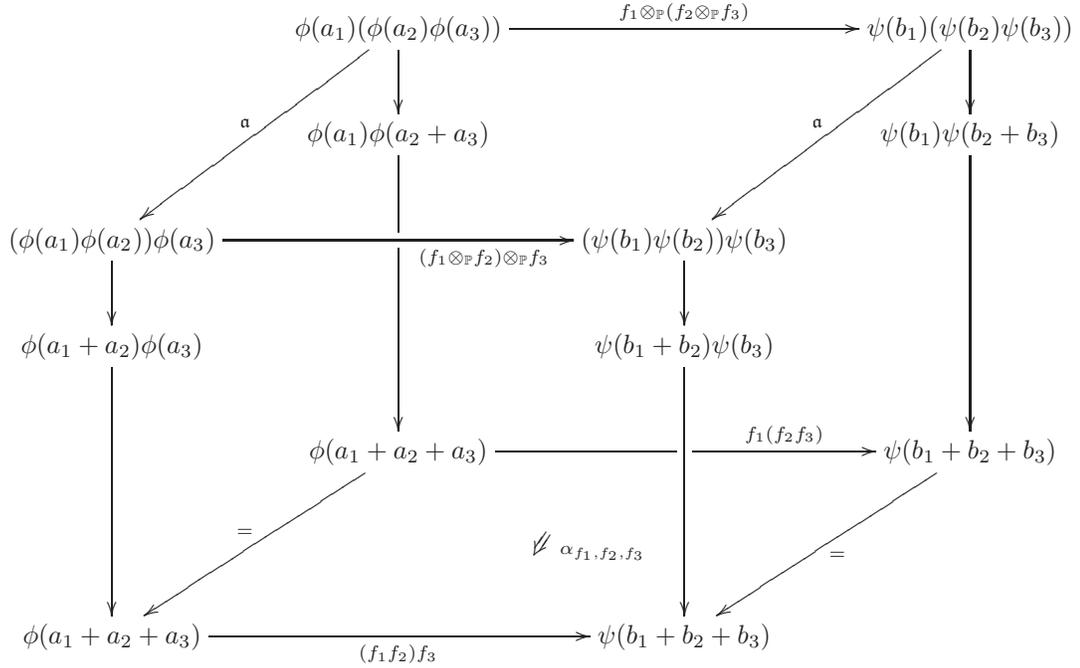
    $$\xymatrix{& \phi(a_1)(\phi(a_2)\phi(a_3)) \ar[d] \ar[ddl]_{\gotika} \ar[rr]^{f_1 \otimes_{\twoP}(f_2\otimes_{\twoP} f_3)} && \psi(b_1)(\psi(b_2)\psi(b_3)) \ar[d] \ar[ddl]_{\gotika}\\
                & \phi(a_1)\phi(a_2+a_3) \ar'[d][ddd] && \psi(b_1)\psi(b_2+b_3) \ar[ddd]\\
                (\phi(a_1)\phi(a_2))\phi(a_3) \ar[d] \ar[rr]_(0.65){(f_1 \otimes_{\twoP} f_2) \otimes_{\twoP} f_3} && (\psi(b_1)\psi(b_2))\psi(b_3) \ar[d] &\\
                \phi(a_1+a_2)\phi(a_3) \ar[ddd] && \psi(b_1+b_2)\psi(b_3) \ar[ddd]&\\
                & \phi(a_1+a_2+a_3) \ar'[r][rr]^(0.35){f_1(f_2f_3)} \ar[ddl]_{=} \ddrtwocell<\omit>{\mspace{80mu}\alpha_{f_1,f_2,f_3}} && \psi(b_1+b_2+b_3) \ar[ddl]^{=} \\
                &&&\\
                \phi(a_1+a_2+a_3) \ar[rr]_{(f_1f_2)f_3} && \psi(b_1+b_2+b_3) &}$$
    \normalsize
    The other 2-isomorphisms of the cube are, the left and right 2-isomorphisms represent the compatibility of the additive 2-functors $\psi$ and $\phi$ with the associativity constraint (see Data HTD5 in \cite{MR1261589}), the back and front ones are of the form $N_m$, the top one is given by the associativity constraint $\gotika$ of $\twoP$ on the 1-morphisms.

    The object $I:=(0_A,e,0_B)$, where $0_A$ (resp. $0_B$) is the unit object in $A$ (resp. in $B$) and $e$ is defined by the 2-commutative diagram
    $$\xymatrix{1_{\twoP} \ar[rr]^{=} \ar[dd] \ddrrtwocell<\omit>{\mspace{20mu}N_u}&& 1_{\twoP} \ar[dd]\\
                &&&\\
                \phi(0_A) \ar[rr]_{e}&& \psi(0_B)}$$

     is the unit in the fibered product $A \times_{\twoP} B$. $I$ comes with the functorial morphisms $l_{(a,f,b)}:=(a,ef,L_f,f,b)$ and $r_{(a,f,b)}:=(a,fe,R_f,f,b)$, where $L_f$ is defined as the 2-isomorphism of the front face that makes the diagram commute (similar diagram for $R_f$).
    \footnotesize
    $$\xymatrix{& \phi(0_A) \otimes_{\twoP} \phi(a) \ar[rr]^{e \otimes_{\twoP} f} \ar[dl] && \psi(0_B) \otimes_{\twoP} \psi(b) \ar[dl]\\
                 \phi(0_A+a) \ar[dd]_(0.65){=} \ar[rr]_(0.6){ef} &\rtwocell<\omit>{<4>\mspace{10mu}L_f}& \psi(0_B+b) \ar[dd]^(0.65){=} \\
                & 1_{\twoP} \otimes_{\twoP} \phi(a) \ar'[u][uu] \ar'[r][rr]^(0.3){1_{\twoP} \otimes_{\twoP} f} \ar[dl] && 1_{\twoP} \otimes_{\twoP} \psi(b)\ar[uu] \ar[dl]\\
                \phi(a) \ar[rr]_{f} && \psi(b)}$$
    \normalsize
    The other 2-isomorphisms of the cube are, the left and right 2-isomorphisms represent the compatibility of the additive 2-functors $\psi$ and $\phi$ with the unit constraint (see Data HTD6 in \cite{MR1261589}), the top and bottom ones are of the form $N_m$, and the back one is of the form $N_u$.
\item[] \emph{braiding}: The morphism between $(a_1,f_1,b_1) \otimes (a_2,f_2,b_2)$ and $(a_2,f_2,b_2) \otimes (a_1,f_1,b_1)$ is given by $(a_1+a_2,f_1f_2,\beta_{f_1,f_2},f_2f_1,b_1+b_2)$, where $\beta_{f_1,f_2}$ is the 2-isomorphism of the bottom face of the commutative cube.
    \footnotesize
    \begin{equation}\label{diagram:braiding}
    \xymatrix{& \phi(a_1) \otimes_\twoP \phi(a_2) \ar[dl]_{R} \ar[rr]^{f_1 \otimes_{\twoP} f_2} \ar'[d][dd] && \psi(b_1) \otimes_\twoP \psi(b_2) \ar[dl]_{R} \ar[dd]\\
                \phi(a_2) \otimes_\twoP \phi(a_1) \ar[rr]^(0.6){f_2 \otimes_{\twoP} f_1} \ar[dd] && \psi(b_2) \otimes_{\twoP} \psi(b_1) \ar[dd]&\\
                & \phi(a_1+a_2) \drtwocell<\omit>{\mspace{40mu}\beta_{f_1,f_2}} \ar'[r]^(0.8){f_1f_2}[rr] \ar[dl]_{=} && \psi(b_1+b_2) \ar[dl]^{=} \\
                \phi(a_1+a_2) \ar[rr]_{f_2f_1} && \psi(b_1+b_2)&}
    \end{equation}
    \normalsize
    The other 2-isomorphisms of the cube are, the left and right 2-isomorphisms represent the compatibility of the additive 2-functors $\psi$ and $\phi$ with the braiding structure \cite{Tatar_thesis}, the front and back ones are of the form $N_m$, and the top one represents the compatibility of $R_{-,-}$ with $\twoP$.
\item[] \emph{group like}: Inverse of an object $(a,f,b)$ is defined as $(-a,g,-b)$, where there exists a 2-isomorphism $\gamma: fg \Ra e$ defined as the 2-morphism of the front face that makes the cube commutative.

    \footnotesize
    $$\xymatrix{& \phi(a) \otimes_{\twoP} \phi(-a) \ar[dl] \ar[rr]^{f \otimes_{\twoP} g} \ar'[d][dd] && \psi(b) \otimes_{\twoP} \psi(-b) \ar[dl] \ar[dd]\\
                 \phi(a+(-a)) \ar[rr]_(0.6){fg} \ar[dd]_(0.65){=} & \rtwocell<\omit>{<5>\gamma}& \psi(b+(-b)) \ar[dd]^(0.65){=}\\
                & 1_\twoP \ar'[r][rr]^(0.35){=} \ar[dl] && 1_\twoP \ar[dl] \\
                 \phi(0_A) \ar[rr]_{e} && \psi(0_B) &}$$
    \normalsize
    The other 2-isomorphisms of the cube are, the left and right 2-isomorphisms represent the compatibility of the additive 2-functors $\psi$ and $\phi$ with the inverse object constraint \cite{Tatar_thesis}, the top (resp. bottom) one is of the form $N_m$ (resp. $N_u$), the back one is the inverse object constraint $\gotiki$.
\item[] \emph{symmetry}: We have to verify that the 2-morphism of the bottom face of the diagram obtained by concatenation of the appropriate two cubes of the form (\ref{diagram:braiding}) is identity. This follows from the fact that, 2-morphism of the top face of the concatenated cube pastes to identity with the help of the structural 2-morphisms of type (\ref{diagram:strongly_braided}).
\item[] \emph{Picard:} The morphism from $(a,f,b) \otimes (a,f,b)$ to itself is identity because the 2-morphism of the top face of the diagram (\ref{diagram:braiding}) becomes identity when it is pasted with (\ref{diagram:2_picard}).
\end{itemize}
The compatibility conditions for each structure are trivially satisfied.
\end{proof}

\subsection{Picard 2-Stack Associated to a Complex}\label{subsection:assoc_Picard_2_stack}
An immediate example of a Picard 2-stack is the Picard 2-stack associated to a complex of abelian sheaves which is in a sense the only example (see Lemma \ref{lemma:ess_surjective}). It is already explained in \cite{noohi-2005} and in \cite{Aldrovandi2009687} how to associate a 2-groupoid to a length 3 complex. However, this 2-groupoid is not a 2-stack. It is not even a 2-prestack (i.e. 1-morphisms only form a prestack but not a stack and 2-descent data are not effective). Therefore to obtain a 2-stack one has to apply the stackification twice. Instead, we are going to use a torsor model for associated stacks. It is more geometric, intuitive, and can be found in \cite{MR2387582} for the abelian case, and in \cite{Aldrovandi2009687} for the non-abelian case.

Consider $A^{\bullet}$ a complex of abelian sheaves. Let $\oneA$ be the associated Picard stack, that is $[A^{-2} \ra A^{-1}]^{\sim} \simeq \TORS(A^{-2},A^{-1})$ and let $\Lambda_A: \oneA \ra A^0$ be an additive functor of Picard stacks,  where $A^0$ is considered as a discrete stack (no non-trivial morphisms). It associates to an object $(L,s)$ in $\TORS(A^{-2},A^{-1})$ an element $\lambda_A(s)$ in $A^0$.

We consider $\TORS(\oneA,A^0)$ consisting of pairs $(\oneL,s)$, where $\oneL$ is an $\oneA$-torsor and $s:\oneL \ra A^0$ is an $\oneA$-equivariant map with respect to $\Lambda_A$. A morphism between any two pairs is given by another pair $(F,\gamma)$

$$(F,\gamma):\xymatrix@1{(\oneL_1,s_1) \ar[r] & (\oneL_2,s_2)},$$
where $F$ is an $\oneA$-torsor morphism compatible with the torsor structure up to $\gamma$. $F$ also fits into the commutative diagram.

$$
\xymatrix{\oneL_1 \ar[rr]^F \ar[ddr]_{s_1} &\ddtwocell<\omit>{=}& \oneL_2\ar[ddl]^{s_2}\\
          &&&\\
          & A^0&}
$$
A 2-morphism

$$\xymatrix@1{(\oneL_1,s_1) \rrtwocell^{(F,\gamma)}_{(G,\delta)}{\theta} && (\oneL_2,s_2)}$$
is given by a natural transformation $\theta: F \Ra G$ that makes the diagram

\begin{equation*}
\xymatrix{\oneL_1 \times \oneA \ar@/^0.3cm/[rr]^{F \times 1}_(0.4){}="a"_{}="d" \ar@/_0.3cm/[rr]_{G \times 1}_(0.4){}="b"^{}="e" \ar[dd]^{}="c" && \oneL_2 \times \oneA \ar[dd]\\
          &&&\\
          \oneL_1 \ar@{-->}@/^0.3cm/[rr]^F_{}="f" \ar@/_0.3cm/[rr]_{G}^{}="g" && \oneL_2
          \ar@{:>}@/_/^{\gamma}"a";"c" \ar@{=>}@/^/^{\delta}"b";"c" \ar@{=>}^{\theta \times 1}"d";"e" \ar@{=>}^{\theta}"f";"g"}
\end{equation*}
commute. It is an immediate result of the following proposition that the 2-stack  $\TORS(\oneA,A^0)$, which we have just constructed is Picard.

\begin{proposition}
For any $\oneA \ra \oneB$ morphism of Picard stacks, $\TORS(\oneA,\oneB)$ is a Picard 2-stack.
\end{proposition}

\begin{proof}
From \cite[\S 6.3.4]{Aldrovandi2009687}, it follows that $\TORS(\oneA,\oneB)$ is a 2-stack. Its group-like structure is defined in \cite[\S 4.5]{MR1086889} and \emph{Picardness} is relatively easy to verify.
\end{proof}

\begin{definition}\label{definition:assoc_Picar_2_stack}
For any complex of abelian sheaves $A^{\bullet}$, we define the Picard 2-stack associated to $A^{\bullet}$ as $\TORS(\oneA,A^0)$.
\end{definition}

\subsection{Homotopy Exact Sequence}\label{subsection:exact_sequence}
Let $\TORS(\oneA,A^0)$ be the associated Picard 2-stack to $A^{\bullet}$, then there is a sequence of Picard 2-stacks

\begin{equation}\label{exact_sequence_assoc_2_stack}
\xymatrix@1{**[l]\oneA \ar[r]^{\Lambda_A} & A^0 \ar[r]^{\pi_A} & **[r]\TORS(\oneA,A^0)},
\end{equation}
where $A^0$ is considered as discrete Picard 2-stack (no non-trivial 1-morphisms and 2-morphisms). The morphism $\pi_A$ assigns to an element $a$ of $A^0(U)$ the pair $(\oneA,a)$, where $a$ is identified with the morphism $\oneA \ra A^0$ sending $1_{\oneA}=(A^{-2},\delta_{A})$ to $a$. (\ref{exact_sequence_assoc_2_stack}) is homotopy exact in the sense that $\oneA$ satisfies the pullback diagram.

\begin{equation}\label{diagram:exact_sequance}
\xymatrix{\oneA \ar[rr] \ar[dd]_{\Lambda_A} && 0 \ar[dd]\\
          &&&\\
          A^0 \ar[rr]_{\pi_A} && **[r]\TORS(\oneA,A^0) \uulltwocell<\omit>}
\end{equation}

Since $\oneA$ is the Picard stack associated to the morphism of abelian sheaves $\delta_A:A^{-2} \ra A^{-1}$, it fits into the commutative pullback square of Picard stacks (see the proof of non-abelian version of Proposition 8.3.2 in \cite{Aldrovandi2009687}).

\begin{equation}\label{diagram:exact_sequance_0}
\xymatrix{A^{-2} \ar[rr] \ar[dd]_{\delta_A} && 0 \ar[dd]\\
          &&&\\
          A^{-1} \ar[rr]_{\pi_{\oneA}} && \oneA \uulltwocell<\omit>}
\end{equation}

Then pasting the diagrams \ref{diagram:exact_sequance} and \ref{diagram:exact_sequance_0} at $\oneA$, we obtain

\begin{equation}\label{diagram:exact_sequance_final}
\xymatrix{A^{-2} \ar[rr] \ar[dd]_{\delta_A} && 0 \ar[dd] &&\\
            &&&&& \\
            A^{-1} \ar[rr]_{\pi_{\oneA}} \ar@{-->}@/_/[ddrr]_{\lambda_A} && \oneA \ar[dd]_{\Lambda_A} \ar[rr] \uulltwocell<\omit>&& 0 \ar[dd]\\
            &&&&& \\
            && A^0 \ar[rr]_{\pi_{A}} && **[r]\TORS(\oneA,A^0) \uulltwocell<\omit>}
\end{equation}

\subsection{The 3-category of Picard 2-Stacks}\label{subsection:3_category_of_Picard_2_stacks}
Picard 2-stacks over $\ES$ form an obvious 3-category which we denote by $\twoCh(\ES)$. $\twoCh(\ES)$ has a hom-2-groupoid consisting of additive 2-functors, weakly invertible natural 2-transformations, and strict modifications. For any two Picard 2-stacks $\twoP$ and $\twoQ$, it is denoted by $\Homsf(\twoP,\twoQ)$. If $\twoP$ and $\twoQ$ are Picard 2-stacks associated to complexes of abelian sheaves $A^{\bullet}$ and $B^{\bullet}$, then the hom-2-groupoid will be denoted as $\Homsf(A^{\bullet},B^{\bullet})$.

\section{The 3-category of Complexes of Abelian Sheaves}\label{section:3-category_of_complexes}
We start with a definition of a 3-category $\rmC^{[-2,0]}(\ES)$ of complexes of abelian sheaves over $\ES$. We end with an explicit construction of a trihomomorphism $\twoch$ between $\rmC^{[-2,0]}(\ES)$ and the 3-category $\twoCh(\ES)$ of Picard 2-stacks over $\ES$.

\subsection{Definition of $\rmC^{[-2,0]}(\ES)$}\label{subsection:3_category_of_complexes}
Although the 3-category of complexes is very well known, in order to setup our notation and terminology, we will describe it explicitly. Its objects are length 3 complexes of abelian sheaves placed in degrees $[-2,0]$. For a pair of objects $A^{\bullet}$ ,$B^{\bullet}$, the hom-2-groupoid $\Hom_{\rmC^{[-2,0]}(\ES)}(A^{\bullet},B^{\bullet})$ is defined as follows:

\begin{itemize}
\item A 1-morphism $f^{\bullet}: A^{\bullet} \ra B^{\bullet}$ is a degree 0 map given by strictly commutative squares.
    \begin{equation}\label{diagram:bicategory_1-morphism}
        \xymatrix{A^{-2} \ar[rr]^{\delta_A} \ar[dd]^{f^{-2}}  && A^{-1} \ar[rr]^{\lambda_A} \ar[dd]^{f^{-1}} && A^{0} \ar[dd]^{f^0}\\
          &&&&\\
          B^{-2} \ar[rr]_{\delta_B} && B^{-1} \ar[rr]_{\lambda_B} && B^{0}}
    \end{equation}
\item A 2-morphism $s^{\bullet}:f^{\bullet} \Ra g^{\bullet}$ is a homotopy map given by the diagram \begin{equation}\label{diagram:bicategory_2-morphism}
        \xymatrix{A^{-2} \ar[rr]^{\delta_A} \ar@/_/[dd]_{g^{-2}} \ar@/^/[dd]^{f^{-2}}  && A^{-1} \ar[rr]^{\lambda_A} \ar@/_/[dd]_{g^{-1}} \ar@/^/[dd]^{f^{-1}} \ar[ddll]|{s^{-1}} && A^{0} \ar@/_/[dd]_{g^0} \ar@/^/[dd]^{f^0} \ar[ddll]|{s^0}\\
          &&&&\\
          B^{-2} \ar[rr]_{\delta_B} && B^{-1} \ar[rr]_{\lambda_B} && B^{0}}
\end{equation}
satisfying the relations
\begin{eqnarray*}
g^0-f^0= \lambda_B \circ s^0, & g^{-1}-f^{-1}=\delta_B \circ s^{-1} + s^0 \circ \lambda_A, & g^{-2}-f^{-2}=s^{-1} \circ \delta_A.
\end{eqnarray*}
\item A 3-morphism $v^{\bullet}: s^{\bullet} \Rrightarrow t^{\bullet}$ is a homotopy map between homotopies $s^{\bullet}$ and $t^{\bullet}$ given by the diagram
\begin{equation}\label{diagram:bicategory_3-morphism}
\xymatrix{A^{-2} \ar[rr]^{\delta_A} \ar@/_/[dd]_{g^{-2}} \ar@/^/[dd]^{f^{-2}}  && A^{-1} \ar[rr]^{\lambda_A} \ar@/_/[dd] \ar@/^/[dd] \ar@/_/[ddll]|{s^{-1}} \ar@/^/[ddll]|{t^{-1}} && A^{0} \ar@/_/[dd]_{g^0} \ar@/^/[dd]^{f^0} \ar@/_/[ddll]|{s^0} \ar@/^/[ddll]|{t^{0}} \ar@{-->}[ddllll]^(0.7)v\\
          &&&&\\
          B^{-2} \ar[rr]_{\delta_B} && B^{-1} \ar[rr]_{\lambda_B} && B^{0}}
\end{equation}
satisfying the relations
\begin{eqnarray*}
s^0-t^0=\delta_B \circ v, & s^{-1} - t^{-1} = -v \circ \lambda_A.
\end{eqnarray*}
\end{itemize}

\begin{remark}
In fact, the hom-2-groupoid $\Hom_{\rmC^{[-2,0]}(\ES)}(A^{\bullet},B^{\bullet})$ is the 2-groupoid associated to $\tau^{\leq 0}(\Hom^{\bullet}(A^{\bullet},B^{\bullet}))$, the smooth truncation of the hom complex $\Hom^{\bullet}(A^{\bullet},B^{\bullet})$, that is to the complex
$$\xymatrix@1{\Hom^{-2}(A^{\bullet},B^{\bullet}) \ar[r] & \Hom^{-1}(A^{\bullet},B^{\bullet}) \ar[r] & Z^0(\Hom^{0}(A^{\bullet},B^{\bullet}))}$$
of abelian groups, where for $i=1,2$ the elements of $\Hom^{-i}(A^{\bullet},B^{\bullet})$ are morphisms of complexes from $A^{\bullet}$ to $B^{\bullet}$ of degree $-i$, and where $Z^0(\Hom^{0}(A^{\bullet},B^{\bullet}))$ is the abelian group of cocycles.
\end{remark}

\subsection{Abelian Sheaves and Picard 2-Stacks}\label{subsection:trihomomorphism}
\begin{lemma}\label{lemma:trihomomorphism}
There is a trihomomorphism

\begin{equation}\label{functor:2ch}
\twoch: \xymatrix@1{\rmC^{[-2,0]}(\ES) \ar[r] & \twoCh(\ES)}
\end{equation}
between the 3-category $\rmC^{[-2,0]}(\ES)$ of complexes of abelian sheaves and the 3-category $\twoCh(\ES)$ of Picard 2-stacks.
\end{lemma}

\begin{proof}
We will give a step by step construction of the trihomomorphism and leave the verification of the axioms to the reader.
\begin{itemize}
\item Using the notations in section \ref{subsection:assoc_Picard_2_stack}, given a complex $A^{\bullet}$, we define $\twoch(A^{\bullet})$ as the associated Picard 2-stack, that is $\twoch(A^{\bullet}):=\TORS(\oneA,A^0)$.
\item For any morphism $f^{\bullet}:A^{\bullet} \ra B^{\bullet}$ of complexes (see diagram (\ref{diagram:bicategory_1-morphism})), there exists a commutative square of Picard stacks

\begin{equation}\label{diagram:1-morphism}
\xymatrix{\oneA \ar[rr]^{\Lambda_A} \ar[dd]_{F}  &\ddtwocell<\omit>{=}& A^0 \ar[dd]^{f^0}\\
          &&\\
          \oneB \ar[rr]_{\Lambda_B} && B^0}
\end{equation}
where $F$ is induced by $f^{\bullet < 0}: A^{\bullet < 0} \ra B^{\bullet <0}$. From the square (\ref{diagram:1-morphism}), we construct a 1-morphism $\twoch(f^{\bullet})$ in $\twoCh(\ES)$

$$\twoch(f^{\bullet}):\xymatrix@1{\TORS(\oneA,A^0) \ar[r] & \TORS(\oneB,B^0)}$$
that sends an $(\oneA,A^0)$-torsor $(\oneL,x)$ to $(\oneL \wedge^{\oneA}_F \oneB,f^0 \circ x+\Lambda_{B})$ where $\oneL \wedge^{\oneA}_F \oneB$ denotes the contracted product of the $\oneA$-torsors $\oneL$ and $\oneB$ such that the $\oneA$-torsor structure of $\oneB$ is induced by the morphism $F$. For details, the reader can refer to \cite[\S 6.7]{MR1086889} and \cite[\S 5.1,\S 6.1]{MR2387582}.
\item For any 2-morphism $s^{\bullet}:f^{\bullet} \Ra g^{\bullet}$ of complexes (see diagram (\ref{diagram:bicategory_2-morphism})), there exists a diagram of Picard stacks

\begin{equation}\label{diagram:2-morphism}
\xymatrix{\oneA \ddtwocell^{F}_{G}{\omit} \ar[rr]^{\Lambda_A} && A^0 \ddtwocell^{f^0}_{g^0} {\omit} \ar[ddll]|{\hat{s}^0}\\
          &&\\
          \oneB \ar[rr]_{\Lambda_B} && B^0}
\end{equation}
such that for any $(L,a)$ in $\oneA$, we have the relation

\begin{eqnarray*}
G(L,a)-F(L,a)=\hat{s}^0 \circ \Lambda_A(L,a) & \textrm{with}& \hat{s}^0(a)=(B^{-2},s^0(a)).
\end{eqnarray*}
From the relation, we construct a natural 2-transformation $\theta$:

$$\xymatrix@1{\TORS(\oneA,A^0) \ar@/^0.3cm/[rrr]^{\twoch(f^{\bullet})}_{}="a" \ar@/_0.3cm/[rrr]^{}="b"_{\twoch(g^{\bullet})} &&& \TORS(\oneB,B^0) \ar@{=>} "a";"b"^{\theta}}$$
in $\twoCh(\ES)$ that assigns to any object $(\oneL,x)$ in $\TORS(\oneA,A^0)$ a 1-morphism $\theta_{(\oneL,x)}$

\begin{equation}\label{1-morphism}
\theta_{(\oneL,x)}:\xymatrix@1{(\oneL \wedge^{\oneA}_F \oneB,x_F) \ar[r] & (\oneL \wedge^{\oneA}_G \oneB, x_G)}
\end{equation}
in $\TORS(\oneB,B^0)$, where $x_F=f^0 \circ x + \Lambda_B$ and $x_G=g^0 \circ x + \Lambda_B$. The morphism (\ref{1-morphism}) is defined by sending $(l,b)$ to $(l,b-s^0 \circ x(l))$.
\item For any 3-morphism $v^{\bullet}:s^{\bullet} \Rrightarrow t^{\bullet}$ of complexes (see diagram \ref{diagram:bicategory_3-morphism}), there exists a modification $\Gamma$:

$$\xymatrix@1{\TORS(\oneA,A^0) \ar@/^0.5cm/[rrr]^{\twoch(f^{\bullet})}_{}="a" \ar@/_0.5cm/[rrr]_{\twoch(g^{\bullet})}^{}="b"  \ar@{}|{\theta \Da \mspace{10mu}{\kesir{\Rrightarrow}{\Gamma}} \mspace{10mu} \Da \phi} "a";"b" &&& \TORS(\oneB,B^0)}$$
in $\twoCh(\ES)$ that assigns to any $(\oneL,x)$ object of $\TORS(\oneA,A^0)$ a natural 2-transformation $\Gamma_{(\oneL,x)}$,

$$\xymatrix@1{(\oneL \wedge^{\oneA}_F \oneB,x_F) \ar@/^0.5cm/[rrr]^{\theta_{(\oneL,x)}}_{}="a" \ar@/_0.5cm/[rrr]_{\phi_{(\oneL,x)}}^{}="b" &&& (\oneL \wedge^{\oneA}_G \oneB,x_G) \ar@{}|{\Da \Gamma_{(\oneL,x)}} "a";"b"}$$
in $\TORS(\oneB,B^0)$, where $\theta_{(\oneL,x)}$, $\phi_{(\oneL,x)}$ are of the form (\ref{1-morphism}). The natural 2-transformation $\Gamma_{(\oneL,x)}$ is defined by assigning to any object $(l,b)$ in $(\oneL \wedge^{\oneA}_F \oneB,x_F)$ a morphism

$$\Gamma_{(\oneL,x)}(l,b):\xymatrix{(l,b-s^0 \circ x(l)) \ar[r] & (l,b-t^0 \circ x(l))}$$
in $(\oneL \wedge^{\oneA}_G \oneB,x_G)$ given by the triple $(\textit{id}_l,1_{\oneA},\beta)$ with $\beta$ being the isomorphism

$$\xymatrix{b-s^0 \circ x(l) \ar[r]& b-s^0 \circ x(l)+\delta_B \circ v \circ x(l)},$$
and $\textit{id}_l$ the identity of $l$ in $\oneL$, and $1_{\oneA}$ the unit element in $\oneA$.
\end{itemize}
\end{proof}

\section{Weak Morphisms of Complexes of Abelian Sheaves}\label{section:fractions}
We fix two complexes of abelian sheaves $A^{\bullet}$ and $B^{\bullet}$. We define $\Frac(A^{\bullet},B^{\bullet})$ a weakened analog of the hom-2-groupoid $\Hom_{\rmC^{[-2,0]}(\ES)}(A^{\bullet},B^{\bullet})$. We also prove that $\Frac(A^{\bullet},B^{\bullet})$ is a bigroupoid.

\subsection{Definition of $\Frac(A^{\bullet},B^{\bullet})$}
$\Frac(A^{\bullet},B^{\bullet})$ is a consists of objects, 1-morphisms, and 2-morphisms such that:
\begin{itemize}
\item An object is an ordered triple $(q,M^{\bullet},p)$, called fraction

$$(q,M^{\bullet},p):\xymatrix@1{A^{\bullet} & M^{\bullet} \ar[l]_q \ar[r]^p & B^{\bullet}}$$
with $M^{\bullet}$ a complex of abelian sheaves, $p$ a morphism of complexes, and $q$ a quasi-isomorphism.
\item A 1-morphism from the fraction $(q_1,M_1^{\bullet},p_1)$ to the fraction $(q_2,M_2^{\bullet},p_2)$ is an ordered triple $(r,K^{\bullet},s)$ with $K^{\bullet}$ a complex of abelian sheaves, $r$ and $s$ quasi-isomorphisms making the diagram

\begin{equation}\label{diagram:1-morphism_of_fractions}
\xymatrix{&& M_1^{\bullet} \ar[drr]^{p_1} \ar[dll]_{q_1} && \\
            A^{\bullet} && K^{\bullet} \ar[rr]|{p} \ar[ll]|q \ar@{-->}[u]_s \ar@{-->}[d]^r && B^{\bullet}\\
            && M_2^{\bullet} \ar[urr]_{p_2} \ar[ull]^{q_2} &&}
\end{equation}
commutative.
\item A 2-morphism from the 1-morphism $(r_1,K_1^{\bullet},s_1)$ to the 1-morphism $(r_2,K_2^{\bullet},s_2)$ is an isomorphism $t^{\bullet}: K_1^{\bullet} \ra K_2^{\bullet}$ of complexes of abelian sheaves such that the diagram that we will call ``\emph{diamond}''

\begin{equation}\label{diagram:diamond}
\xymatrix@u{&&A^{\bullet}&&\\
            &&&&\\
            &&&K_1^{\bullet}\ar[uul] \ar[ddddl] \ar[dlll]|{r_1} \ar[dr]|{s_1} \ar@{-->}[ddll]|{t^{\bullet}}&\\
            M_2^{\bullet} \ar[uuurr]|{q_2} \ar[dddrr]|{p_2}&&&& M_1^{\bullet} \ar[uuull]|{q_1} \ar[dddll]|{p_1}\\
            &K_2^{\bullet}\ar[uuuur] \ar[ddr] \ar[ul]|{r_2} \ar[urrr]|(0.4){s_2} &&&\\
            &&&&\\
            &&B^{\bullet}&&}
\end{equation}
commutes.
\end{itemize}

\begin{remark}
For reasons of clarity, we will represent the above 2-morphism by the following planar commutative diagram

\begin{equation*}
\xymatrix{&&&M_1^{\bullet} \ar[ddrrr]^{p_1} \ar[ddlll]_{q_1}&&&\\
          &&&&&&\\
          A^{\bullet}&& K_1^{\bullet} \ar@{-->}[rr]^{t^{\bullet}} \ar[uur]|{s_1} \ar[ddr]|{r_1} && K_2^{\bullet} \ar[uul]|{s_2} \ar[ddl]|{r_2} && B^{\bullet}\\
          &&&&&&&&\\
          &&&M_2^{\bullet}\ar[uurrr]_{p_2} \ar[uulll]^{q_2}&&&}
\end{equation*}
where we have ignored the maps from $K^{\bullet}$'s to $A^{\bullet}$ and $B^{\bullet}$.
\end{remark}

\begin{remark}\label{remark:bigroupoid}
From the definition of 2-morphisms, it is immediate that all 2-morphisms are isomorphisms.
\end{remark}

\subsection{$\mathsf{Frac}(A^{\bullet},B^{\bullet})$ is a bigroupoid}
\begin{proposition}\label{proposition:bicategory_Frac}
Let $A^{\bullet}$ and $B^{\bullet}$ be two complexes of abelian sheaves. Then $\Frac(A^{\bullet},B^{\bullet})$ is a bigroupoid.
\end{proposition}

\begin{proof}
We will describe the necessary data to define the bigroupoid without verifying that they satisfy the required axioms.
\begin{itemize}
\item For any two composable morphisms $(r_1,K^{\bullet}_1,s_1):(q_1,M^{\bullet}_1,p_1) \ra (q_2,M^{\bullet}_2,p_2)$ and $(r_2,K^{\bullet}_2,s_2):(q_2,M^{\bullet}_2,p_2) \ra (q_3,M^{\bullet}_3,p_3)$ shown by the diagram

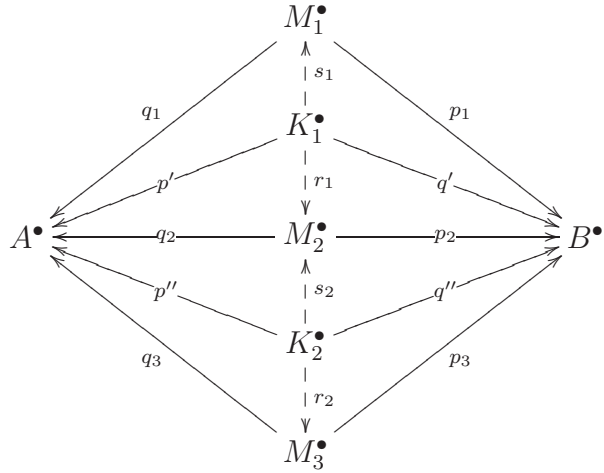
\begin{equation*}
\xymatrix{&&& M_1^{\bullet} \ar[ddrrr]^{p_1} \ar[ddlll]_{q_1} &&& \\
          &&& K_1^{\bullet} \ar[drrr]|{q'} \ar[dlll]|{p'} \ar@{-->}[u]_{s_1} \ar@{-->}[d]^{r_1} &&&\\
          A^{\bullet} &&& M_2^{\bullet} \ar[rrr]|{p_2} \ar[lll]|{q_2} &&& B^{\bullet}\\
          &&& K_2^{\bullet} \ar[urrr]|{q''} \ar[ulll]|{p''} \ar@{-->}[u]_{s_2} \ar@{-->}[d]^{r_2} &&&\\
          &&& M_3^{\bullet} \ar[uurrr]_{p_3} \ar[uulll]^{q_3} &&&}
\end{equation*}
the composition is defined by the pullback diagram.

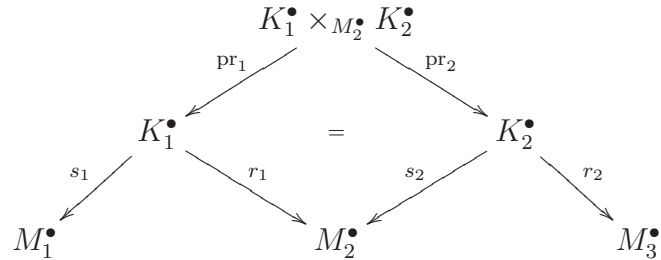
\begin{equation*}
\xymatrix{&& K_1^{\bullet} \times_{M_2^{\bullet}} K_2^{\bullet}\ar[dr]^{\pr_2} \ar[dl]_{\pr_1} &&\\
          &K_1^{\bullet} \ar[dr]^{r_1} \ar[dl]_{s_1} \ar@{}[rr]|{=}&& K_2^{\bullet} \ar[dr]^{r_2} \ar[dl]_{s_2}&\\
          M_1^{\bullet}&&M_2^{\bullet}&& M_3^{\bullet}}
\end{equation*}
That is the composition is the triple $(r_2 \circ \pr_2, K_1^{\bullet} \times_{M_2^{\bullet}} K_2^{\bullet}, s_1 \circ \pr_1)$.

\item For two 2-morphisms $t_1^{\bullet}:(r_1,K_1^{\bullet},s_1) \Ra (r_2,K_2^{\bullet},s_2)$ and $t_2^{\bullet}:(r_2,K_2^{\bullet},s_2) \Ra (r_3,K_3^{\bullet},s_3)$ shown by the diagram

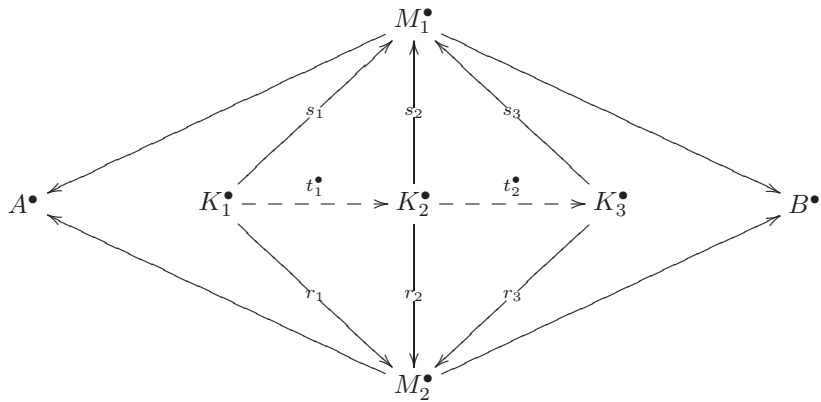
\footnotesize
\begin{equation*}
\xymatrix{&&&&M_1^{\bullet}\ar[ddrrrr] \ar[ddllll]&&&&\\
          &&&&&&&&\\
          A^{\bullet}&&K_1^{\bullet} \ar@{-->}[rr]^{t_1^{\bullet}} \ar[uurr]|{s_1} \ar[ddrr]|{r_1} && K_2^{\bullet} \ar[uu]|{s_2} \ar[dd]|{r_2} \ar@{-->}[rr]^{t_2^{\bullet}} && K_3^{\bullet} \ar[uull]|{s_3} \ar[ddll]|{r_3} && B^{\bullet}\\
          &&&&&&&&\\
          &&&&M_2^{\bullet}\ar[uurrrr] \ar[uullll]&&&&}
\end{equation*}
\normalsize the vertical composition is defined by $t_2^{\bullet} \circ t_1^{\bullet}$.

\item For two 2-morphisms $t^{\bullet}:(r_1,K_1^{\bullet},s_1) \Ra (r_2,K_2^{\bullet},s_2)$ and $u^{\bullet}:(r_1',L_1^{\bullet},s_1') \Ra (r_2',L_2^{\bullet},s_2')$ shown by the diagram

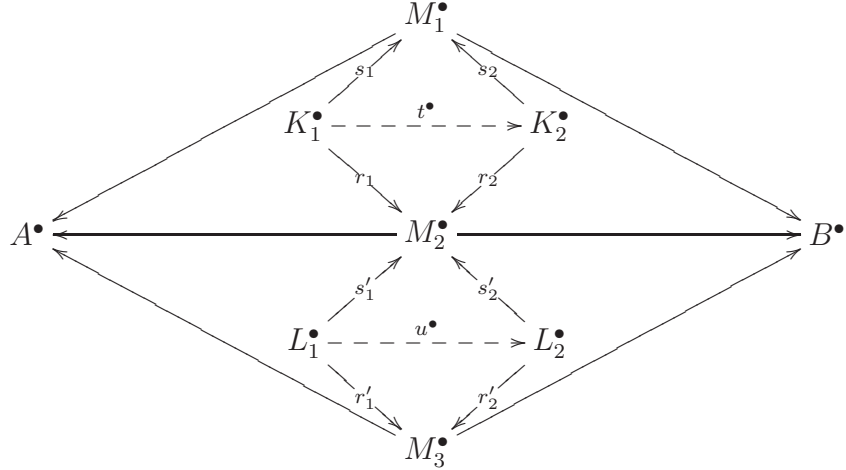
\begin{equation*}
\xymatrix{&&&&M_1^{\bullet}\ar[ddrrrr] \ar[ddllll]&&&&\\
          &&& K_1^{\bullet} \ar[dr]|{r_1} \ar[ur]|{s_1} \ar@{-->}[rr]^{t^{\bullet}} &&K_2^{\bullet}\ar[ul]|{s_2} \ar[dl]|{r_2} &&&\\
          A^{\bullet}&&&& M_2^{\bullet} \ar[rrrr] \ar[llll] &&&& B^{\bullet}\\
          &&& L_1^{\bullet} \ar[ur]|{s_1'} \ar[dr]|{r_1'} \ar@{-->}[rr]^{u^{\bullet}} && L_2^{\bullet}\ar[dl]|{r_2'} \ar[ul]|{s_2'} &&&\\
          &&&&M_3^{\bullet}\ar[uurrrr] \ar[uullll]&&&&}
\end{equation*}
the horizontal composition is given by the natural morphism  $K_1^{\bullet}\times_{M_2^{\bullet}} L_1^{\bullet} \ra K_2^{\bullet}\times_{M_2^{\bullet}} L_2^{\bullet}$ between the pullbacks of pairs $(r_1,s_1')$ and $(r_2,s_2')$ over $M^{\bullet}_2$.
\end{itemize}

Any three composable 1-morphisms $(r_1,K_1^{\bullet},s_1)$, $(r_2,K_2^{\bullet},s_2)$, and $(r_3,K_3^{\bullet},s_3)$ can be pictured as a sequence of three fractions

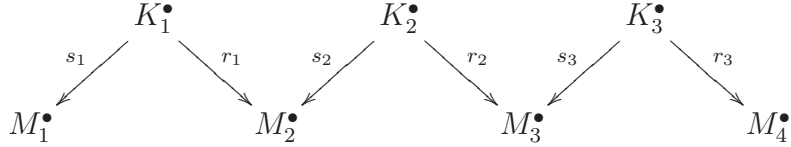
\begin{equation*}
\xymatrix{& K_1^{\bullet} \ar[dr]^{r_1} \ar[dl]_{s_1} && K_2^{\bullet} \ar[dr]^{r_2} \ar[dl]_{s_2} && K_3^{\bullet} \ar[dr]^{r_3} \ar[dl]_{s_3}&\\
          M_1^{\bullet} && M_2^{\bullet} && M_3^{\bullet} && M_4^{\bullet}}
\end{equation*}
simply by ignoring the maps to $A^{\bullet}$ and $B^{\bullet}$. They can be composed in two different ways, either first by pulling back over $M_2^{\bullet}$ then over $M_3^{\bullet}$ or vice versa. The resulting fractions will be $(r, (K_1^{\bullet} \times_{M_2^{\bullet}} K_2^{\bullet}) \times_{M_3^{\bullet}} K_3^{\bullet},s)$ and $(r',K_1^{\bullet} \times_{M_2^{\bullet}} (K_2^{\bullet} \times_{M_3^{\bullet}} K_3^{\bullet}),s')$, respectively, where $r$ and $r'$ (resp.$s$ and $s'$) are equal to $r_3$ (resp.$s_1$) composed with appropriate projection maps. The 2-isomorphism between these fractions is given by the natural isomorphism between the pullbacks. Thus, the associativity of composition of 1-morphisms is weak.

We also observe that 1-morphisms are weakly invertible. Let $(r,K^{\bullet},s)$ be a 1-morphism from $(q_1,M_1^{\bullet},p_1)$ to $(q_2,M_2^{\bullet},p_2)$, then $(s,K^{\bullet},r)$ is a weak inverse of $(r,K^{\bullet},s)$ in the sense that the composition $(r \circ \pr, K^{\bullet} \times_{M_2^{\bullet}} K^{\bullet}, r \circ \pr)$ is equivalent to the identity, that is there is a natural 2-transformation $\theta: r \circ \pr \Ra \id \circ (r \circ \pr)$ as shown in the below diagram.

\begin{equation*}
\xymatrix{&&&M_1^{\bullet} \ar[ddrrr]^{p_1} \ar[ddlll]_{q_1}&&&\\
          &&&&&&\\
          A^{\bullet}&& K^{\bullet} \ar@{-->}[rr]^{r \circ \pr} \ar[uur]|{r \circ \pr}  \ar[ddr]|{r \circ \pr} &\uutwocell<\omit>{\theta}& M_1^{\bullet} \ar[uul]|{\id} \ar[ddl]|{\id} && B^{\bullet}\\
          &&&&&&&&\\
          &&&M_1^{\bullet}\ar[uurrr]_{p_1} \ar[uulll]^{q_1} \uutwocell<\omit>{\theta}&&&}
\end{equation*}

Thus, $\Frac(A^{\bullet},B^{\bullet})$ is a bigroupoid.
\end{proof}

\begin{remark}
In the terminology of \cite{Aldrovandi2009687-III}, what we have called fractions are called in the non-abelian context weak morphisms of 2-crossed modules or butterflies of gr-stacks or bats of sheaves.
\end{remark}

\section{Biequivalence of $\Frac(A^{\bullet},B^{\bullet})$ and $\Homsf(A^{\bullet},B^{\bullet})$}\label{section:bieqv_of_hom_bicats}

Fix again two complexes of abelian sheaves $A^{\bullet}$ and $B^{\bullet}$. In this section, we prove that the bigroupoid $\Frac(A^{\bullet}, B^{\bullet})$ of fractions defined in Section \ref{section:fractions} is biequivalent to the 2-groupoid $\Homsf(A^{\bullet}, B^{\bullet})$ of additive 2-functors from $\twoch(A^{\bullet})$ to $\twoch(A^{\bullet})$ defined in Section \ref{subsection:3_category_of_Picard_2_stacks}.

\subsection{Morphisms of Picard 2-Stacks as Fractions}

\begin{lemma}\label{lemma:quasi_isomorphism_vs_equivalences}
A morphism $f: A^{\bullet} \ra B^{\bullet}$ is a quasi-isomorphism if and only if

$$\twoch(f): \xymatrix{\twoch(A^{\bullet}) \ar[r] & \twoch(B^{\bullet})}$$
is a biequivalence.
\end{lemma}

\begin{proof}
Given $f: A^{\bullet} \ra B^{\bullet}$ a morphism of complexes, we know how to induce a morphism of Picard 2-stacks (see construction of trihomomorphism $\twoch(f)$). It is also known that a 2-stack (not necessarily Picard) can be seen as a 2-gerbe over its own $\pi_0$ bounded by the stack $\oneA\mspace{-3mu}\mathit{ut}(\mathrm{I})$ of automorphisms of identity \cite[\S8.1]{MR1301844}. In particular, the Picard 2-stacks $\TORS(\oneA,A^0)$ and $\TORS(\oneB,B^0)$ are 2-gerbes over their own $\pi_0$ bounded by $\oneA\mspace{-3mu}\mathit{ut}(\mathrm{I}_{\twoch(A^{\bullet})}) \simeq {[A^{-2} \ra \ker(\delta_A)]}^{\sim}$ and  $\oneA\mspace{-3mu}\mathit{ut}(\mathrm{I}_{\twoch(B^{\bullet})}) \simeq {[B^{-2} \ra \ker(\delta_B)]}^{\sim}$, respectively. Furthermore, if $f$ is a quasi-isomorphism, then $H^{-i}(A^{\bullet}) \simeq H^{-i}(B^{\bullet})$ for $i=0,1,2$ and thus, $\pi_i(\twoch(A^{\bullet})) \simeq \pi_i(\twoch(B^{\bullet}))$ for $i=0,1,2$. So $\TORS(\oneA,A^0)$ and $\TORS(\oneB,B^0)$ are 2-gerbes with equivalent bands. Therefore they are equivalent.
\end{proof}

Given an additive 2-functor $F$ in $\Homsf(A^{\bullet},B^{\bullet})$, we will show in the next lemma that there is a corresponding object in $\Frac(A^{\bullet},B^{\bullet})$.

\begin{lemma}\label{lemma:fully_faithful}
For any additive 2-functor $F:\twoch(A^{\bullet}) \ra \twoch(B^{\bullet})$, there exists a fraction $(q,M^{\bullet},p)$ such that $F \circ \twoch(q) \simeq \twoch(p)$.
\end{lemma}

\begin{proof}
From the sequences
$$\xymatrix{\oneA \ar[r]^{\Lambda_A} & A^0 \ar[r]^(.4){\pi_A} & \twoch(A^{\bullet})} \textrm{ and } \xymatrix{\oneB \ar[r]^{\Lambda_B} & B^0 \ar[r]^(0.4){\pi_B} & \twoch(B^{\bullet})},$$
we can construct the commutative diagram
\begin{equation}\label{diagram:lemma_fully_faithful}
\xymatrix{& \oneA \times \oneB \ar[dd]^{\mu_F} \ar[dr] \ar[dl]& \\
            \oneA \ar[dr]^{\nu_F} \ar[dd]_{\Lambda_A} & & \oneB \ar[dd]^{\Lambda_B} \ar[dl]_{\xi_F}\\
            & \oneE_F \ar[dr]^{\pr_2} \ar[dl]_{\pr_1}&\\
            A^0 \ar[d]_{\pi_A} \ar[drr]^{F \circ \pi_A}&& B^0 \ar[d]^{\pi_B}\\
            \twoch(A^{\bullet}) \ar[rr]_F & & \twoch(B^{\bullet})}
\end{equation}
where $\oneE_F:=A^0 \times_{F,B} B^0$. It follows from the commutativity of the above diagram that $\mu_F=(\Lambda_A,\Lambda_B)$. The sequence

\begin{equation}\label{pullback_sequence}
\xymatrix@1{\oneB \ar[r]^{\xi_F} & \oneE_F \ar[r]^{\pr_1} & A^0}
\end{equation}
is homotopy exact since it is the pullback of the exact sequence $\oneB \ra B^0 \ra \twoch(B^{\bullet})$. From Lemma \ref{lemma:pullback}, it follows that $\oneE_F$ is a Picard stack. Therefore by \cite[Proposition 8.3.2]{Aldrovandi2009687}, there exists a length 2 complex $E^{\bullet}=[\delta_E:E^{-1}_F \ra E^0_F]$ of abelian sheaves such that the associated Picard stack $\TORS(E^{-1}_F,E^0_F)$ is equivalent to $\oneE_F$. Then by \cite[Theorem 8.3.1]{Aldrovandi2009687}, there exists a butterfly representing $\mu_F$:

\begin{equation}\label{diagram:butterfly(AxB)}
\xymatrix{ A^{-2} \times B^{-2} \ar[dd]_{\delta_A \times \delta_B} \ar[dr]^{\kappa} && E^{-1}_F \ar[dd]^{\delta_E} \ar[dl]_{\imath}\\%
           & P_F \ar[dr]^{\jmath} \ar[dl]_{\rho} & \\%
           A^{-1} \times B^{-1} \ar[d]_{\pi_{\oneA} \times \pi_{\oneB}} &&  E^0_F \ar[d]^{\pi_{\oneE_F}}\\
           \oneA \times \oneB \ar[rr]_{\mu_F} && \oneE_F}
\end{equation}
with $P_F \simeq (A^{-1} \times B^{-1})\times_{\oneE_F} E^0_F$. From a different perspective, this butterfly can be seen as

\begin{equation}\label{diagram:different_perspective}
\xymatrix{0 \ar[r] \ar[d] & E^{-1}_F \ar[r]^{\delta_E} \ar[d]_{\imath} & E^0_F \ar[d]^{\id}\\
          A^{-2} \times B^{-2} \ar[r]^{\kappa} \ar[d]_{\id} & P_F \ar[r]^{\jmath} \ar[d]_{\rho} & E^0_F \ar[d]\\
          A^{-2} \times B^{-2} \ar[r]_{\delta_A \times \delta_B} & A^{-1} \times B^{-1} \ar[r] & 0}
\end{equation}
where each column is an exact sequence of abelian sheaves. The only non-trivial sequence is the second column and its  exactness follows from the definition of a butterfly (\ref{definition:butterfly}). So we have a short exact sequence of complexes of abelian sheaves
\begin{equation}\label{exact_sequence}
\xymatrix{0 \ar[r] & E^{\bullet}_F \ar[r] & M^{\bullet}_F \ar[r] & A^{\bullet < 0} \times B^{\bullet < 0} \ar[r] & 0},
\end{equation}
where

\begin{eqnarray}
\label{the_complex} M^{\bullet}_F & := & \xymatrix@1{A^{-2} \times B^{-2} \ar[r] & P_F \ar[r] & E^0_F},\\
\nonumber E^{\bullet}_F & := & \xymatrix@1{0 \ar[r] & E^{-1}_F \ar[r] & E^0_F},\\
\nonumber A^{\bullet < 0} \times B^{\bullet < 0}& := &\xymatrix@1{A^{-2} \times B^{-2} \ar[r] & A^{-1} \times B^{-1} \ar[r] & 0}.
\end{eqnarray}

From the lower part of the diagram (\ref{diagram:different_perspective}) and the definition of $P_F$, we deduce that there are morphisms of complexes

\begin{equation}\label{diagram:two_morphisms_of_complexes}
\xymatrix@C=3pc{&A^{-2} \times B^{-2} \ar[dd]_{\kappa} \ar[dl]_{\pr_1} \ar[dr]^{\pr_2}&\\
          A^{-2} \ar[dd]_{\delta_A}&&B^{-2}\ar[dd]^{\delta_B}\\
          &P_F \ar[dd]_{\jmath} \ar[dl]_{\pr_2 \circ \rho} \ar[dr]^{\pr_1 \circ \rho}&\\
          A^{-1} \ar[dd]_{\lambda_A} && B^{-1}\ar[dd]^{\lambda_B}\\
          &E^0_F \ar[dl] \ar[dr]&\\
          A^{0} && B^{0}\\
          &M^{\bullet}_F \ar[dl]_{q} \ar[dr]^{p}&\\
          A^{\bullet}&&B^{\bullet}}
\end{equation}

We claim that $q$ is a quasi-isomorphism, that is
$$
\begin{tabular}{ccc}
$H^{-2}(M^{\bullet}_F) \simeq \ker(\delta_A),$ & $H^{-1}(M^{\bullet}_F) \simeq \ker(\lambda_A) / \im(\delta_A),$ & $H^0(M^{\bullet}_F) \simeq \coker(\lambda_A).$
\end{tabular}
$$
Indeed, from the exact sequence (\ref{exact_sequence}), we obtain the long exact sequence of homology sheaves

\begin{equation}\label{long_exact_sequnce_cohomology}
\xymatrix{0 \ar[r] & H^{-2}(M^{\bullet}_F) \ar[r] & H^{-2}(A^{\bullet < 0}) \times H^{-2}(B^{\bullet < 0}) \ar[r] & H^{-1}(E^{\bullet}_F) \ar`[d]`[dlll]`[ddlll]_{\partial}[ddlll]& \\ &&&&\\ H^{-1}(M^{\bullet}_F) \ar[r] & H^{-1}(A^{\bullet <0}) \times H^{-1}(B^{\bullet <0}) \ar[r] & H^0(E^{\bullet}_F) \ar[r] & H^0(M^{\bullet}_F) \ar[r] & 0}.
\end{equation}

On the other hand, by \cite[Proposition 6.2.6]{Aldrovandi2009687} applied to the exact sequence (\ref{pullback_sequence}), we get a long exact sequence of homotopy groups

\begin{equation}\label{rouseau_sequence}
\xymatrix{0 \ar[r] & \pi_1(\oneB) \ar[r] & \pi_1(\oneE_F) \ar[r] & \pi_1(A^0) \ar[r] & \pi_0(\oneB) \ar[r] & \pi_0(\oneE_F) \ar[r] & \pi_0(A^0) \ar[r] & 0.}
\end{equation}

Since $\pi_1(A^0) = H^{-1}(A^0) = 0$ and $\pi_0(A^0) = H^0(A^0) = A^0$, it follows from (\ref{rouseau_sequence}) that we have an isomorphism

\begin{equation}\label{rouseau_sequence_part1}
\xymatrix{H^{-2}(B^{\bullet <0}) \ar[r]^{\simeq} & H^{-1}(E^{\bullet}_F)}
\end{equation}
and an exact sequence

\begin{equation}\label{rouseau_sequence_part2}
\xymatrix{0 \ar[r] & H^{-1}(B^{\bullet <0}) \ar[r] & H^0(E^{\bullet}_F) \ar[r] & A^0 \ar[r] & 0.}
\end{equation}

(\ref{rouseau_sequence_part1}) implies that $\partial=0$ in (\ref{long_exact_sequnce_cohomology}). Therefore from (\ref{long_exact_sequnce_cohomology}) again, we obtain a short exact sequence
$$\xymatrix{0 \ar[r] & H^{-2}(M^{\bullet}_F) \ar[r] & H^{-2}(A^{\bullet <0}) \times H^{-2}(B^{\bullet <0}) \ar[r] & H^{-1}(E^{\bullet}_F) \ar[r] & 0}$$
from which we deduce that $H^{-2}(M^{\bullet}_F) \simeq H^{-2}(A^{\bullet <0})=\ker(\delta_A)$.

Now, apply the snake lemma to the short exact sequence (\ref{rouseau_sequence_part2}) and to
$$\xymatrix{0 \ar[r] & H^{-1}(B^{\bullet <0}) \ar[r] & H^{-1}(A^{\bullet <0}) \times H^{-1}(B^{\bullet <0}) \ar[r] & H^{-1}(A^{\bullet <0}) \ar[r] & 0}$$
in order to get the dashed exact sequence

$$\xymatrix{&& 0 \ar[dd] \ar[r] & H^{-1}(M^{\bullet}_F) \ar[dd] \ar[r] & \ker(\lambda_A)/\im(\delta_A) \ar[dd] \ar@{-->}`[ddd]`[dddlll]`[ddddddlll]`[ddddddlllr][ddddddlllr]& \\&&&&&\\
            0 \ar[rr]  && H^{-1}(B^{\bullet <0}) \ar[r] \ar[dd] & H^{-1}(A^{\bullet <0}) \times H^{-1}(B^{\bullet <0}) \ar[r] \ar[dd] & H^{-1}(A^{\bullet <0}) \ar[r] \ar[dd] & 0\\&&&&&\\
            0 \ar[rr] && H^{-1}(B^{\bullet <0}) \ar[r] \ar[dd] & H^0(E^{\bullet}_F) \ar[r] \ar[dd] & A_0 \ar[r] \ar[dd] & 0\\ &&&&&\\
            && 0 \ar[r] & H^0(M^{\bullet}_F) \ar[r] & \coker(\lambda_A)}$$
from which it follows $H^{-1}(M^{\bullet}_F) \simeq \ker(\lambda_A)/\im(\delta_A)$, and $H^0(M^{\bullet}_F) \simeq \coker(A^0)$ as wanted.

We end this proof by showing that $F \circ \twoch(q) \simeq \twoch(p)$. (\ref{diagram:two_morphisms_of_complexes}) induces a diagram of Picard 2-stacks

\begin{equation}\label{diagram:assoc_2_stacks}
\xymatrix{& \twoch(M^{\bullet}_F) \ar[dr]^{\twoch(p)} \ar[dl]_{\twoch(q)}&\\
            \twoch(A^{\bullet}) \ar[rr]_F&& \twoch(B^{\bullet})}.
\end{equation}
We claim that (\ref{diagram:assoc_2_stacks}) commutes up to a natural 2-transformation. To show that, it is enough to look at $\twoch(M^{\bullet}_F)$ locally. Given $U \in \ES$,  $\twoch(M^{\bullet}_F)_{U}$ is the 2-groupoid associated to the complex of abelian groups (for the definition of the 2-groupoid associated to a complex see \cite{Aldrovandi2009687} or \cite{noohi-2005})
$$\xymatrix{A^{-2}(U)\times B^{-2}(U) \ar[r]^(0.65){\delta} & P_F(U) \ar[r]^{\lambda} & E^0_F(U)}$$
Then, an object of $\twoch(M^{\bullet}_F)_{U}$ is an element $e$ of $E_F^0(U)$. Since $\oneE_F:=A^0 \times_{F,B} B^0 \simeq \TORS(E_F^{-1}, E_F^0)$, $e$ can be taken as $(a,f, b)$, where $a \in A^0(U)$, $b \in B^0(U)$, and $f:F(a)\ra b$ is a 1-morphism in $\twoch(B^{\bullet})_{U}$.

A 1-morphism of $\twoch(M^{\bullet}_F)_{U}$ from $e_1$ to $e_2$ is given by an element $p$ of $P_F(U)$ such that $\lambda(p)+e_1=e_2$ in $E_F^0(U)$. We can again take $\lambda(p)$, $e_1$, and $e_2$ as $(a,f,b)$, $(a_1,f_1,b_1)$, and $(a_2,f_2,b_2)$, respectively. Therefore, the addition in $E_F^0(U)$ should be replaced by the monoidal operation on $\oneE_F$ between the triples, that is  $(a,f,b) \otimes_{\oneE_F} (a_1,f_1,b_1) = (a_2,f_2,b_2)$. This monoidal operation is described in the proof of the technical Lemma \ref{lemma:pullback}. It creates a diagram commutative up to a 2-isomorphism in $\twoCh(B^{\bullet})_{U}$ that defines $f_2$.

\begin{equation*}
\xymatrix{F(a_2) \ar[dd]_{\simeq} \ar[rr]^{f_2} && b_2 \ar[dd]^{\simeq}\\
          &&&\\
          F(a) \otimes_{B} F(a_1) \ar[rr]_(0.6){f \otimes_{B} f_1} && b \otimes_{B} b_1 \uulltwocell<\omit> {\mspace{-10mu}\theta}}
\end{equation*}

The collection $(f, \theta)$ gives the natural 2-transformation between $\twoch(q) \circ F$ and $\twoch(p)$.

\begin{remark}
Since $q$ is a quasi-isomorphism in $\rmC^{[-2,0]}(\ES)$, the technical lemma \ref{lemma:quasi_isomorphism_vs_equivalences} implies that $\twoch(q)$ is a biequivalence in $\twoCh(\ES)$. Therefore, by choosing an inverse of $\twoch(q)$ up to a natural 2-transformation we can write $F$ as $F \simeq \twoch(p) \circ \twoch(q)^{-1}$.
\end{remark}
\end{proof}

\subsection{Hom-categories of $\Frac(A^{\bullet},B^{\bullet})$ and $\Homsf(A^{\bullet},B^{\bullet})$}
In the next two lemmas, we are going to explore the relation between 1-morphisms (resp. 2-morphisms) of $\Frac(A^{\bullet},B^{\bullet})$ and natural 2-transformations (resp. modifications) of Picard 2-stacks.

Suppose we have a natural 2-transformation $\theta$:

\begin{equation}\label{diagram:natural_2_transformation}
\xymatrix@1{\twoch(A^{\bullet}) \rrtwocell^{F}_{G}{\theta}&& \twoch(B^{\bullet})}
\end{equation}
between the two additive 2-functors $F,G:\twoch(A^{\bullet}) \ra \twoch(B^{\bullet})$. By Lemma \ref{lemma:fully_faithful}, we know that there are fractions $(q_F,M^{\bullet}_F,p_F)$ and $(q_G,M^{\bullet}_G,p_G)$ associated to $F$ and $G$.

\begin{lemma}\label{lemma:isomorphism_classes_of_morphisms}
For any natural 2-transformation $\theta$ as in (\ref{diagram:natural_2_transformation}), there is a  1-morphism in $\Frac(A^{\bullet},B^{\bullet})$ between the fractions $(q_F,M^{\bullet}_F,p_F)$ and $(q_G,M^{\bullet}_G,p_G)$.
\end{lemma}

\begin{proof}
For $F$ and $G$, we have the following diagrams similar to (\ref{diagram:lemma_fully_faithful})

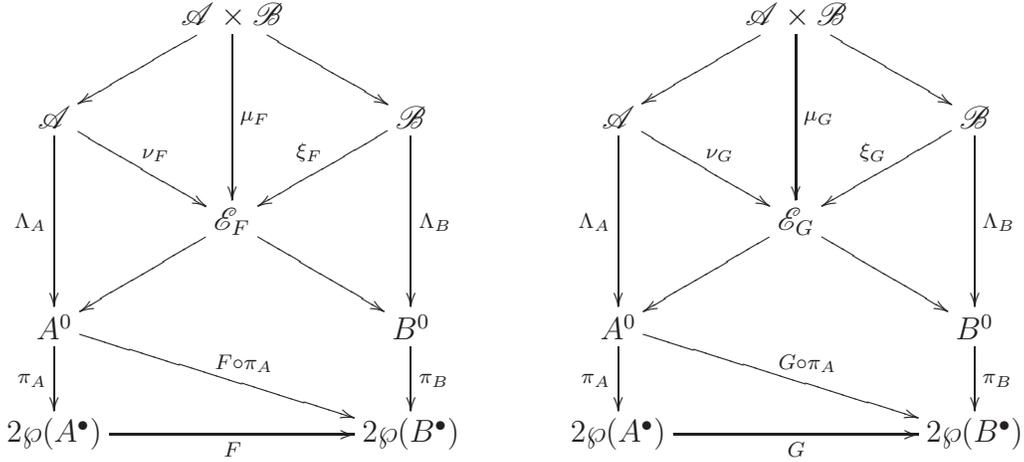
$$
\xymatrix{& \oneA \times \oneB \ar[dd]^{\mu_F} \ar[dr] \ar[dl]& \\
            \oneA \ar[dr]^{\nu_F} \ar[dd]_{\Lambda_A} & & \oneB \ar[dd]^{\Lambda_B} \ar[dl]_{\xi_F}\\
            & \oneE_F \ar[dr] \ar[dl]&\\
            A^0 \ar[d]_{\pi_{A}} \ar[drr]^{F \circ \pi_{A}}&& B^0 \ar[d]^{\pi_{B}}\\
            \twoch(A^{\bullet}) \ar[rr]_F & & \twoch(B^{\bullet})}
\xymatrix{&}
\xymatrix{& \oneA \times \oneB \ar[dd]^{\mu_G} \ar[dr] \ar[dl]& \\
            \oneA \ar[dr]^{\nu_G} \ar[dd]_{\Lambda_A} & & \oneB \ar[dd]^{\Lambda_B} \ar[dl]_{\xi_G}\\
            & \oneE_G \ar[dr] \ar[dl]&\\
            A^0 \ar[d]_{\pi_{A}} \ar[drr]^{G \circ \pi_{A}}&& B^0 \ar[d]^{\pi_{B}}\\
            \twoch(A^{\bullet}) \ar[rr]_G & & \twoch(B^{\bullet})}
$$
where $\oneE_F:=A^0 \times_{F,B} B^0$ and $\oneE_G:=A^0 \times_{G,B} B^0$ are Picard stacks by Lemma \ref{lemma:pullback}. Therefore by \cite[Proposition 8.3.2]{Aldrovandi2009687}, there exist $E_F^{-1} \ra E_F^0$ and $E_G^{-1} \ra E_G^0$ morphisms of abelian sheaves such that the Picard stack associated to them are respectively $\oneE_F$ and $\oneE_G$. The natural 2-transformation $\theta: F \Ra G$ induces an equivalence $H:\oneE_G \ra \oneE_F$ of Picard stacks defined as follows:
\begin{itemize}
\item For any $(a,g,b)$ object of $(\oneE_G)_U$, $H((a,g,b)):=(a,f,b)$, where $f$ fits into the commutative diagram
$$\xymatrix{F(a) \ar[dd]_{\theta_a} \ar[rr]^f \ar@{}[ddrr]|{=} && b \ar[dd]^{\sim}\\
            &&&\\
            G(a) \ar[rr]_{g} && b}$$
\item For any $(a,g,\sigma,g',b)$ morphism of $(\oneE_G)_U$, $H((a,g,\sigma,g',b)):=(a,f,\tau,f',b)$, where $\tau$ is defined by the following whiskering.
$$\xymatrix{F(a) \ar[r]^{\theta_a} & G(a) \rrtwocell^{g}_{g'}{\sigma} && b}$$
\end{itemize}
By \cite[Theorem 8.3.1]{Aldrovandi2009687}, $H$ corresponds to a butterfly $[E^{\bullet}_G, N ,E^{\bullet}_F]$. Since $H$ is an equivalence, this butterfly is flippable.

We compose $H$ and $\mu_G$ by composing their corresponding butterflies

\begin{equation*}
\xymatrix{A^{-2} \times B^{-2} \ar[dd]_{\delta_A \times \delta_B} \ar[dr]^{\kappa'} && E_F^{-1} \ar[dd]^{\delta_E} \ar[dl]_{\imath'}\\%
           & P_G \times_{E_G^0}^{E_G^{-1}} N \ar[dr]^{\jmath'} \ar[dl]_{\rho'} & \\%
           A^{-1} \times B^{-1} \ar[d]_{\pi_{\oneA} \times \pi_{\oneB}} &&  E_F^0 \ar[d]^{\pi_{\oneE_F}}\\
           \oneA \times \oneB \ar[rr]_{H \circ \mu_G} && \oneE_F}
\end{equation*}

where $P_G \times_{E_G^0}^{E_G^{-1}} N$ is pull-out/pull-back construction as defined in \cite[\S 5.1]{Aldrovandi2009687}.

There is also a direct morphism $\mu_F$ from $\oneA \times \oneB$ to $\oneE_F$. $\mu_F$ is equivalent to $H \circ \mu_G$ since they both map an object of $\oneA \times \oneB$ to an object in $\oneE_F$ which is isomorphic to the unit object in $\twoch(B^{\bullet})$. Then by \cite[Theorem 8.3.1]{Aldrovandi2009687}, there exists an isomorphism $k$ between the corresponding butterflies of $\mu_F$ and $H \circ \mu_G$, that is the dotted arrow in the diagram below such that all regions commute.

\begin{equation}\label{diagram:butterlfy_isomorphism}
\xymatrix{A^{-2} \times B^{-2} \ar[dd]_{\delta_A \times \delta_B} \ar[drr]^{\kappa} \ar[rr]^{\kappa'} && P_G \times_{E_G^0}^{E_G^{-1}} N \ar[ddrr]|(0.35){\hole}^(0.65){\jmath'} \ar[ddll]|(0.35){\hole}_(0.65){\rho'} \ar@{-->}[d]^k && E_F^{-1} \ar[dd]^{\delta_E} \ar[dll]_{\imath} \ar[ll]_{\imath'}\\%
           && P_F \ar[drr]_{\jmath} \ar[dll]^{\rho} && \\%
           A^{-1} \times B^{-1} \ar[d]_{\pi_{\oneA} \times \pi_{\oneB}} &&&&  E_F^0 \ar[d]^{\pi_{\oneE_F}}\\
           \oneA \times \oneB \ar@/^0.3cm/[rrrr]^{H \circ \mu_G} \ar@/_0.3cm/[rrrr]_{\mu_F} &&&& \oneE_F}
\end{equation}

Let $M_F^{\bullet}:A^{-2} \times B^{-2} \ra P_F \ra E_F^0$ and $M_G^{\bullet}:A^{-2} \times B^{-2} \ra P_G \ra E_G^0$. We claim that, there exists a complex $K^{\bullet}$ with quasi-isomorphisms $r_F$ and $r_G$ such that all regions in the diagram

\begin{equation}\label{diagram:rhombus}
\xymatrix{&& M_F^{\bullet} \ar^{p_F}[drr] \ar_{q_F}[dll] && \\
            A^{\bullet} && K^{\bullet} \ar[rr]^p \ar[ll]_q \ar[u]_{r_F} \ar[d]^{r_G} && B^{\bullet}\\
            && M_G^{\bullet} \ar_{p_G}[urr] \ar^{q_G}[ull] &&}
\end{equation}
commute.

\emph{Proof of the claim}: Let $K^{\bullet}:A^{-2} \times B^{-2} \ra P_G \times_{E_G^0} N \ra N$ and define $r_F$ by the composition

\begin{equation}\label{diagram:r_F}
\xymatrix{K^{\bullet} \ar[dd]_{r_F} & A^{-2} \times B^{-2} \ar[rr] \ar@{=}[d] && P_G \times _{E^0_G} N \ar[rr] \ar[d]^{\textrm{quotient}} && N \ar[d]^{\textrm{quotient}}\\
          &A^{-2} \times B^{-2} \ar[rr] \ar@{=}[d] &&P_G \times_{E_G^0}^{E_G^{-1}} N \ar[rr] \ar[d] && N/E_G^{-1} \ar[d]\\
          M^{\bullet}_F&A^{-2} \times B^{-2} \ar[rr] &&P_F \ar[rr] &&E^0_F}
\end{equation}
and $r_G$ by the diagram

\begin{equation}\label{diagram:r_G}
\xymatrix{K^{\bullet} \ar[d]_{r_G}&A^{-2} \times B^{-2} \ar[rr] \ar@{=}[d] && P_G \times_{E_G^0}N \ar[d] \ar[rr] && N \ar[d]\\
          M^{\bullet}_G&A^{-2} \times B^{-2} \ar[rr]&& P_G \ar[rr]  && E^0_G}
\end{equation}
The commutativity of the diagram (\ref{diagram:r_F}) follows from composition of butterflies. Since $P_G \times_{E_G^0}^{E_G^{-1}} N \simeq P_F$ and the butterfly $[E^{\bullet}_G,N,E^{\bullet}_F]$ is flippable, $r_F$ is a quasi-isomorphism. The diagram (\ref{diagram:r_G}) commutes because its left square is a pullback. This implies that $r_G$ is a quasi-isomorphism.

It remains to show that $q_F \circ r_F = q_G \circ r_G$, that is in the diagram below each column closes to a commutative square.

\begin{equation*}
\xymatrix{A^{\bullet}&A^{-2}\ar[rr] && A^{-1} \ar[rr] && A^0\\
          M^{\bullet}_F\ar[u]^{q_F}&A^{-2} \times B^{-2} \ar[u] \ar[rr] && P_F \ar[rr] \ar[u] && E^0_F \ar[u]\\
          K^{\bullet} \ar[u]^{r_F}\ar[d]_{r_G}&A^{-2} \times B^{-2} \ar@{=}[u] \ar[rr] \ar@{=}[d] && P_G \times_{E_G^0}N \ar[u] \ar[d] \ar[rr] && N \ar[u] \ar[d]\\
          M^{\bullet}_G\ar[d]_{q_G}&A^{-2} \times B^{-2} \ar[rr] \ar[d] && P_G \ar[d] \ar[rr] && E^0_G \ar[d]\\
          A^{\bullet}&A^{-2}\ar[rr] && A^{-1} \ar[rr] && A^0}
\end{equation*}
It is obvious for the first column. The commutativity of the triangles

\begin{eqnarray*}
\xymatrix{P_G \times^{E^{-1}_G}_{E^0_G} N \ar[rr]^k \ar[ddrr]_{\rho'} && P_F \ar[dd]^{\rho}\\
            &&\\
            &&A^{-1} \times B^{-1}}
&
\xymatrix{\oneE_G \ar[rr]^H \ar[ddr]_{\pr_1} && \oneE_F \ar[ddl]^{\pr_2}\\
            &&\\
            &A^0&}
\end{eqnarray*}
imply that the middle and last columns close to a commutative square, respectively (the first triangle is extracted from diagram (\ref{diagram:butterlfy_isomorphism})).

In the same way, we also show that $p_F \circ r_F = p_G \circ r_G$.
\end{proof}

Now, suppose we have a modification $\Gamma$:

\begin{equation}\label{diagram:modification}
\xymatrix@1{\twoch(A^{\bullet}) \ar@/^0.5cm/[rr]^F_{}="a" \ar@/_0.5cm/[rr]_G^{}="b"  \ar@{}|{\theta \Da \mspace{10mu}{\kesir{\Rrightarrow} {\Gamma}} \mspace{10mu} \Da \phi} "a";"b"&& \twoch(B^{\bullet})}
\end{equation}
between two natural 2-transformations $\theta,\phi: F \Ra G$. We have proved in Lemmas \ref{lemma:fully_faithful} and \ref{lemma:isomorphism_classes_of_morphisms} that both $\theta$ and $\phi$ correspond to a 1-morphism in $\Frac(A^{\bullet},B^{\bullet})$.

\begin{lemma}\label{lemma:modifications}
Given a modification $\Gamma$ as in (\ref{diagram:modification}), there exists a 2-morphism between the two 1-morphisms corresponding to $\theta$ and $\phi$.
\end{lemma}

\begin{proof}
Using the same notations as in Lemma \ref{lemma:isomorphism_classes_of_morphisms}, we construct a diagram of Picard stacks

\begin{equation*}
\xymatrix@1{\oneE_G \rrtwocell^{H_{\theta}}_{H_{\phi}} {T} && \oneE_F},
\end{equation*}
where $T$ is a natural transformation. For any object $(a,g,b)$ in $\oneE_G$, $T_{(a,g,b)}$ is a morphism in $\oneE_F$ defined by

$$\xymatrix@1{F(a) \rrtwocell<5>^{f_{\theta}}_{f_{\phi}} {\mspace{36mu}1_g * \Gamma_{a}} && b},$$
where
$$\xymatrix@1{F(a) \rrtwocell<5>^{\theta_a}_{\phi_a} {\mspace{9mu}\Gamma_a} && G(a)},$$
and $H_{\theta}(a,g,b)=(a,f_{\theta},b)$, $H_{\phi}(a,g,b)=(a,f_{\phi},b)$ . By \cite[Theorem 5.3.6]{Aldrovandi2009687}, the natural transformation $T$ corresponds to an isomorphism $t$ between the centers of the butterflies associated to $H_{\theta}$ and $H_{\phi}$.

\begin{equation}\label{diagram:modification_2}
\xymatrix{E_G^0 \ar[dd]_{\delta_{E_G}} \ar[drr]^{\kappa_{\phi}} \ar[rr]^{\kappa_{\theta}} && N_{\theta} \ar[ddrr]|(0.34){\hole}^(0.65){\jmath_{\theta}} \ar[ddll]|(0.34){\hole}_(0.65){\rho_{\theta}} && E_F^{-1} \ar[dd]^{\delta_{E_F}} \ar[dll]_{\imath_{\phi}} \ar[ll]_{\imath_{\theta}}\\%
           && N_{\phi} \ar[drr]_{\jmath_{\phi}} \ar[dll]^{\rho_{\phi}} \ar@{-->}[u]_t&& \\%
           E_G^0 \ar[d]_{\pi_{\oneE_G}} &&&&  E_F^0 \ar[d]^{\pi_{\oneE_F}}\\
           \oneE_G \ar@/^0.2cm/[rrrr]^{H_{\theta}}_{}="a" \ar@/_0.2cm/[rrrr]_{H_{\phi}}^{}="b" &&&& \oneE_F \ar@{} |{\Da \mspace{10 mu} T}"a";"b"}
\end{equation}
$t$ induces an isomorphism of complexes $t^{\bullet}$.

$$\xymatrix{K^{\bullet}_{\phi} \ar[d]_{t^{\bullet}}&A^{-2} \times B^{-2} \ar[rr] \ar@{=}[d] && P_G \times_{E_G^0}N_{\phi} \ar[d]^{\textrm{id} \times t} \ar[rr] && N_{\phi} \ar[d]^t\\
          K^{\bullet}_{\theta}&A^{-2} \times B^{-2} \ar[rr]&& P_G \times_{E_G^0}N_{\theta}\ar[rr] && N_{\theta}}$$
The proof finishes by showing that all the regions in the diagram (\ref{diagram:diamond}) commute. The only regions, whose commutativity are non-trivial, are the triangles in the middle sharing an edge marked by the isomorphism $t^{\bullet}$. They commute as well since in the diagram below

$$\xymatrix{M^{\bullet}_G &A^{-2} \times B^{-2} \ar[rr] && P_G \ar[rr]&& E^0_G\\
            K^{\bullet}_{\phi} \ar[u]^{r_{G,\phi}}\ar[d]_{t^{\bullet}}& A^{-2} \times B^{-2} \ar@{=}[u] \ar[rr] \ar@{=}[d] && P_G \times_{E_G^0}N_{\phi} \ar[u]_{\pr_1} \ar[d]^{\id \times t} \ar[rr] && N_{\phi} \ar[u]_{\rho_{\phi}} \ar[d]^t\\
            K^{\bullet}_{\theta} \ar[d]_{r_{G,\theta}} & A^{-2} \times B^{-2} \ar[rr] \ar@{=}[d] && P_G \times_{E_G^0}N_{\theta} \ar[d]^{\pr_1} \ar[rr] && N_{\theta} \ar[d]^{\rho_{\theta}}\\
            M^{\bullet}_G & A^{-2} \times B^{-2} \ar[rr]&& P_G \ar[rr] && E^0_G}$$
each column closes to a commutative triangle. This is immediate for the first two columns. The triangle formed by the last column commutes as well, since it is a piece of the commutative diagram (\ref{diagram:modification_2}).
\end{proof}

For any two complexes of abelian sheaves $A^{\bullet}$ and $B^{\bullet}$, the proofs of Lemmas \ref{lemma:fully_faithful} and \ref{lemma:isomorphism_classes_of_morphisms} define us a 2-functor

\begin{equation}\label{functor:hom_bicategories}
\twoch_{(A^{\bullet},B^{\bullet})}:\xymatrix@1{\Frac(A^{\bullet},B^{\bullet}) \ar[r] & \Homsf(A^{\bullet},B^{\bullet})}
\end{equation}
between the bigroupoid $\Frac(A^{\bullet},B^{\bullet})$ and the 2-groupoid $\Homsf(A^{\bullet}, B^{\bullet})$ of additive 2-functors between $\twoch(A^{\bullet})$ and $\twoch(B^{\bullet})$ considered as a bigroupoid. In fact, we have proved:

\begin{theorem}\label{theorem:biequivalence_of_Hom_categories}
For any two complexes of abelian sheaves $A^{\bullet}$ and $B^{\bullet}$, $\twoch_{(A^{\bullet},B^{\bullet})}$ is a biequivalence of bigroupoids.
\end{theorem}

\section{The Tricategory of Complexes of Abelian Sheaves}\label{section:tricategory_of_complexes}
After proving in Section \ref{section:bieqv_of_hom_bicats} that for any two complexes of abelian sheaves $A^{\bullet}$ and $B^{\bullet}$, $\Frac(A^{\bullet},B^{\bullet})$ is biequivalent as a bigroupoid to $\Homsf(A^{\bullet},B^{\bullet})$, it is clear that the trihomomorphism $\twoch$ (\ref{functor:2ch}) defined in Section \ref{subsection:trihomomorphism} cannot be a triequivalence. To attain the triequivalence, we need to consider at least a tricategory with same objects as $\rmC^{[-2,0]}(\ES)$ and with hom-bicategories of the form $\Frac(A^{\bullet},B^{\bullet})$. Furthermore, there is the question of essential surjectivity which we deal with in this section.

\subsection{Definition of $\rmT^{[-2,0]}(\ES)$}\label{subsection:tricategory_T}
We define the tricategory $\rmT^{[-2,0]}(\ES)$ promised at the beginning of the section.

\begin{def-prop}\label{proposition:tricategory}
$\rmT^{[-2,0]}(\ES)$ with objects complexes of abelian sheaves, and hom-bigroupoids $\Frac(A^{\bullet},B^{\bullet})$, for any two complexes of abelian sheaves $A^{\bullet}$ and $B^{\bullet}$, is a tricategory.
\end{def-prop}

\begin{proof}
We have to verify that $\rmT^{[-2,0]}(\ES)$ has the data given in \cite[Definition 3.3.1]{Gurski_thesis}.
\begin{itemize}
\item Objects are complexes of abelian sheaves.
\item For any two complexes of abelian sheaves $A^{\bullet}$ and $B^{\bullet}$, $\Frac(A^{\bullet},B^{\bullet})$ is the hom-bicategory.
\item For any three complexes of abelian sheaves $A^{\bullet}$, $B^{\bullet}$, and $C^{\bullet}$, the composition is given by the weak functor
$$\otimes_T:\xymatrix@1{\Frac(A^{\bullet}, B^{\bullet}) \times \Frac{(B^{\bullet},C^{\bullet})} \ar[r] & \Frac{(A^{\bullet},C^{\bullet}})},$$
which is defined on
\footnotesize
\begin{enumerate}
\item \normalsize{objects, by}
$$\xymatrix@!=1.5pc{&M_1^{\bullet} \ar[dl]_{q_1} \ar[dr]^{p_1}="c" &&& M_2^{\bullet} \ar[dl]_{q_2}="d" \ar[dr]^{p_2}="a" &&&& M_1^{\bullet} \times_{B^{\bullet}} M_2^{\bullet} \ar[dll]_{q_1 \circ \pr_1}="b" \ar[drr]^{p_2 \circ \pr_2} &&\\
                        A^{\bullet} && B^{\bullet} &B^{\bullet}&& C^{\bullet} & A^{\bullet}&&&&C^{\bullet} \ar@{}|{=}"a";"b" \ar@{}|{\otimes_T}"c";"d"}$$
\item \normalsize{1-morphisms, by}
$$\xymatrix@!=1.5pc{&M_1^{\bullet} \ar[dl]_{q_1} \ar[dr]^{p_1} &&& M_2^{\bullet} \ar[dl]_{q_2} \ar[dr]^{p_2} &&&& M_1^{\bullet} \times_{B^{\bullet}} M_2^{\bullet} \ar[dll]_{q_1 \circ \pr_1} \ar[drr]^{p_2 \circ \pr_2}&&\\
                        A^{\bullet} &K^{\bullet} \ar[u]|{s_1} \ar[d]|{r_1} \ar[r]|{y_1} \ar[l]|{x_1} & B^{\bullet} \ar@{}[r]|{\otimes_T} & B^{\bullet} & L^{\bullet} \ar[u]|{s_2} \ar[d]|{r_2} \ar[r]|{y_2} \ar[l]|{x_2} & C^{\bullet} \ar@{}[r]|{=}& A^{\bullet}&& K^{\bullet} \times_{B^{\bullet}} L^{\bullet} \ar[u]|{s_1 \times s_2} \ar[d]|{r_1 \times r_2} \ar[rr]|(0.6){y_1 \circ \pr_2} \ar[ll]|(0.6){x_1 \circ \pr_1} && C^{\bullet}\\
                        &N_1^{\bullet} \ar[ul]^{q_1'} \ar[ur]_{p_1'} &&& N_2^{\bullet} \ar[ul]^{q_2'} \ar[ur]_{p_2'} &&&& N_1^{\bullet} \times_{B^{\bullet}} N_2^{\bullet} \ar[ull]^{q_1' \circ \pr_1} \ar[urr]_{p_2' \circ \pr_2}&&}$$
\item{2-morphisms, by}
$$\xy \xymatrix@!=1.5pc{&M_1^{\bullet} \ar[dl]_{q_1} \ar[dr]^{p_1} &&& M_2^{\bullet} \ar[dl]_{q_2} \ar[dr]^{p_2} &&&& M_1^{\bullet} \times_{B^{\bullet}} M_2^{\bullet} \ar[dll]_{q_1 \circ \pr_1} \ar[drr]^{p_2 \circ \pr_2} &&\\
                        A^{\bullet} &K_1^{\bullet} \ra K_2^{\bullet} & B^{\bullet} \ar@{}[r]|{\otimes_T} &B^{\bullet}&L_1^{\bullet} \ra L_2^{\bullet}& C^{\bullet}\ar@{}[r]|{=}& A^{\bullet}&&K_1^{\bullet} \times_{B^{\bullet}} L_1^{\bullet} \ra K_2^{\bullet} \times_{B^{\bullet}} L_2^{\bullet} && C^{\bullet}\\
                        &N_1^{\bullet} \ar[ul]^{q_1'} \ar[ur]_{p_1'} &&& N_2^{\bullet} \ar[ul]^{q_2'} \ar[ur]_{p_2'} &&&& N_1^{\bullet} \times_{B^{\bullet}} N_2^{\bullet} \ar[ull]^{q_1' \circ \pr_1} \ar[urr]_{p_2' \circ \pr_2}&&}
                        \ar (10,-12)*{};(14,-3)*{}
                        \ar (19,-12)*{};(15,-3)*{}
                        \ar (10,-17)*{};(14,-26)*{}
                        \ar (19,-17)*{};(15,-26)*{}
                        \ar (54,-12)*{};(58,-3)*{}
                        \ar (63,-12)*{};(59,-3)*{}
                        \ar (54,-17)*{};(58,-26)*{}
                        \ar (63,-17)*{};(59,-26)*{}
                        \ar (109,-12)*{};(115,-3)*{}
                        \ar (127,-12)*{};(121,-3)*{}
                        \ar (109,-17)*{};(115,-26)*{}
                        \ar (127,-17)*{};(121,-26)*{}
                        \endxy$$
\end{enumerate}
\end{itemize}
\normalsize
We leave defining the rest of the data as well as verifying that they satisfy the axioms to the reader.
\end{proof}

The trihomomorphism (\ref{functor:2ch}) extends to a trihomomorphism

\begin{equation}\label{functor:triequivalence}
\twoch: \xymatrix{\rmT^{[-2,0]}(\ES) \ar[r] & \twoCh(\ES)}
\end{equation}
on the tricategory $\rmT^{[-2,0]}(\ES)$ as follows\footnote[1]{We commit an abuse of notation by calling both functors (\ref{functor:2ch}) and (\ref{functor:triequivalence}) by $\twoch$.}: On objects, it is defined as explained in Section \ref{subsection:trihomomorphism}. On 1-, 2-, 3-morphisms, by the biequivalence $\twoch_{(A^{\bullet},B^{\bullet})}$, where $A^{\bullet}$ and $B^{\bullet}$ are any two complexes of abelian sheaves.

Theorem \ref{theorem:biequivalence_of_Hom_categories} implies that (\ref{functor:triequivalence}) is already fully faithful in the appropriate sense. In order to prove the triequivalence, one needs to show that it is  essentially surjective, as well.

The essential surjectivity depends on the following technical lemma, which is similar to Lemme 1.4.3 in \cite{Deligne}. We give its proof in the Appendix (\ref{Annex}).

\begin{proposition}\label{proposition:ab_group_embedded_in_2_category}
For any set $E$, denote by $\mathbb{Z}(E)$ the free abelian group generated by $E$. Let $\twoC$ be a Picard 2-category and $F_0: E \ra \twoC$ be a set map. Then $F_0$ extends to an additive 2-functor $F: \mathbb{Z}(E) \ra \twoC$ where $\mathbb{Z}(E)$ is considered as a 2-category (trivially Picard).
\end{proposition}

\begin{lemma}\label{lemma:ess_surjective}
Let $\twoP$ be a Picard 2-stack, then there exists a complex of abelian sheaves $A^{\bullet}$ such that $\twoch(A^{\bullet})$ is biequivalent to $\twoP$.
\end{lemma}

\begin{proof}
There is a construction analogous to the skeleton of categories. For any 2-category $\twoP$, we construct $\twosk(\twoP)$ a 2-category that has one object per equivalence class in $\twoP$. We observe that $\twosk(\twoP)$ is a full sub 2-category of $\twoP$, that is the inclusion $\twosk(\twoP) \ra \twoP$ is a biequivalence. Let $\twoP$ be a Picard 2-stack. We note that $\Ob\twosk(\twoP):U \ra \Ob(\twosk(\twoP_U))$ is a presheaf of sets. We consider $A^0$ the abelian sheaf over $\ES$ associated to the presheaf $\{U \ra \mathbb{Z}(\Ob(\twosk(\twoP_U)))\}$ where $\mathbb{Z}(\Ob(\twosk(\twoP_U)))$ is the free abelian group associated to $\Ob(\twosk(\twoP_U))$. By Proposition \ref{proposition:ab_group_embedded_in_2_category}, the inclusion $i:\Ob\twosk(\twoP) \ra \twoP$ extends to

$$\pi_{\twoP}:\xymatrix@1{A^0 \ar[r] & \twoP}$$
an essentially surjective additive 2-functor on $A^0$.

Define $\oneA$ by the pullback diagram

\begin{equation}\label{diagram:pullback1}
\xymatrix{\oneA \ar[rr] \ar[dd]_{\Lambda_A} && 0 \ar[dd]\\
            &&&\\
            A^0 \ar[rr]_{\pi_{\twoP}} && \twoP \uulltwocell<\omit>}
\end{equation}
of morphisms of Picard 2-stacks, which is similar to (\ref{diagram:exact_sequance}). Then, the sequence of Picard 2-stacks
$$\xymatrix{\oneA \ar[r] & A^0 \ar[r] & \twoP}$$
is exact sequence in the sense of Section \ref{subsection:exact_sequence}.

On the other hand, from Lemma \ref{lemma:pullback}, it follows that $\oneA$ is a Picard stack. Therefore by \cite[Proposition 8.3.2]{Aldrovandi2009687}, there exists a morphism of
abelian sheaves $\delta_A:A^{-2} \ra A^{-1}$, where $A^{-2}$ is defined by the pullback diagram

\begin{equation}\label{diagram:pullback2}
\xymatrix{A^{-2} \ar[rr] \ar[dd]_{\delta_A} && 0 \ar[dd]\\
            &&&\\
            A^{-1} \ar[rr]_{\pi_{\oneA}} && \oneA \uulltwocell<\omit>}
\end{equation}
and $\oneA := \TORS(A^{-2},A^{-1})$.

Now putting the diagrams (\ref{diagram:pullback1}) and (\ref{diagram:pullback2}) together,

\begin{equation}
\xymatrix{A^{-2} \ar[rr] \ar[dd]_{\delta_A} && 0 \ar[dd] &&\\
            &&&&& \\
            A^{-1} \ar[rr]_{\pi_{\oneA}} \ar@{-->}@/_/[ddrr]_{\lambda_A} && \oneA \ar[dd]_{\Lambda_A} \ar[rr] \uulltwocell<\omit>&& 0 \ar[dd]\\
            &&&&& \\
            && A^0 \ar[rr]_{\pi_{\twoP}} && \twoP \uulltwocell<\omit>}
\end{equation}
we have a diagram of Picard 2-stacks. It implies that $A^{\bullet}:\xymatrix@1{A^{-2} \ar[r]^{\delta_A} & A^{-1} \ar[r]^{\lambda_A} & A^0}$ is a complex.

The Picard 2-stack associated to $A^{\bullet}$, that is $\twoch(A^{\bullet}):=\TORS(\oneA,A^0)$, verifies by definition the above diagram (see \ref{diagram:exact_sequance_final}).

The biequivalence  $\twoch(A^{\bullet}) \simeq \twoP$ is almost immediate. Essential surjectivity follows from the definition of $\pi_{\twoP}$ and equivalence of hom-categories from the fact that $A^0$ and $0$ pull back to $\oneA$ over $\twoch(A^{\bullet})$ and over $\twoP$.
\end{proof}

\subsection{Main Theorem}
Considering $\twoCh(\ES)$ as a tricategory, our main result follows from Theorem \ref{theorem:biequivalence_of_Hom_categories} and Lemma \ref{lemma:ess_surjective}.

\begin{theorem}\label{theorem:triequivalence}
The trihomomorphism (\ref{functor:triequivalence}) is a triequivalence.
\end{theorem}

An immediate consequence of Theorem \ref{theorem:triequivalence}, which was also the motivation for this paper, is the following.

Let $\twoCh^{\flat \flat}(\ES)$ denote the category of Picard 2-stacks obtained from $\twoCh(\ES)$ by ignoring the modifications and taking as morphisms the equivalence classes of additive 2-functors. Let $\rmD^{[-2,0]}(\ES)$ be the subcategory of the derived category of category of complexes of abelian sheaves $A^{\bullet}$ over $\ES$ with $H^{-i}(A^{\bullet}) \neq 0$ for $i=0,1,2$. We deduce from Theorem \ref{theorem:triequivalence} the following, which generalizes Deligne's result \cite[Proposition 1.4.15]{Deligne} from Picard stacks to Picard 2-stacks.

\begin{corollary}\label{corollary}
The functor (\ref{functor:triequivalence}) induces an equivalence

\begin{equation}
\twoch^{\flat \flat}: \xymatrix@1{\rmD^{[-2,0]}(\ES) \ar[r] & \twoCh^{\flat \flat}(\ES)}
\end{equation}
of categories.
\end{corollary}

\begin{proof}
It is enough to observe from the calculations in Section \ref{section:fractions} that $\pi_0(\Frac(A^{\bullet},B^{\bullet})) \simeq \Hom_{\rmD^{[-2,0]}(\ES)}(A^{\bullet},B^{\bullet})$. Since the objects of $\rmD^{[-2,0]}(\ES)$ are same as the objects of $\rmT^{[-2,0]}(\ES)$, the essential surjectivity follows from Lemma \ref{lemma:ess_surjective}.
\end{proof}

\section{Stackification}\label{section:stackification}
We want to conclude with an informal discussion of stack versions of some of our results. We will assume that all structures are strict unless otherwise stated. Throughout the paper, we dealt with 2- and 3-categories and their weakened versions bi- and tricategories. They can be stackified.

2-stacks over a site are well known \cite{MR1301844}. The collection of 2-stacks over $\ES$, denoted by $\ch(\ES)$, comprise a 3-category structure. We can consider the fibered 3-category $\Ch(\ES)$, whose fiber over $U$ is the 3-category $\ch(\ES /U)$ of 2-stacks over $\ES/U$. In \cite[Remark 1.12]{MR1301844}, Breen claims that $\Ch(\ES)$ is a 3-stack. Hirschowitz and Simpson in \cite{hirschowitz-1998}, generalize this result to weak $n$-stacks.

\begin{theorem*}\cite[Th\'{e}or\`{e}me 20.5]{hirschowitz-1998}
The weak $(n+1)$-prestack of weak $n$-stacks $\mathit{nW}\mathfrak{S}\textsc{tack}(\ES)$ is a weak $(n+1)$-stack over $\ES$.
\end{theorem*}

We can use the above facts to deduce that the 3-prestack of Picard 2-stacks $\PCh(\ES)$ with fibers $\twoCh(\ES /U)$ over $U$ is a 3-stack.

\begin{claim*}
$\HHom(A^{\bullet},B^{\bullet})$ fibered over $\ES$ in 2-groupoids is a 2-stack where for any $U \in \ES$, the 2-groupoid $\Homsf(A^{\bullet}_{|U},B^{\bullet}_{|U})$ of additive 2-functors from $\twoch(A^{\bullet})_{|U}$ to $\twoch(B^{\bullet})_{|U}$ defines the fiber over $U$.
\end{claim*}

We have also fibered analogs for each hom-bicategory $\Frac(A^{\bullet},B^{\bullet})$ and for $\rmT^{[-2,0]}(\ES)$. It follows from the above claim and Theorem \ref{theorem:biequivalence_of_Hom_categories} that the bi-prestack $\FFrac(A^{\bullet},B^{\bullet})$ of fractions from $A^{\bullet}$ to $B^{\bullet}$ with fibers defined by $\Frac(A^{\bullet}_{|U},B^{\bullet}_{|U})$ is a bistack. Then, once an appropriate notion of 3-descent has been specified and all descent data are shown to be effective, we conclude by the characterization proposition \cite[Proposition 10.2]{hirschowitz-1998} for $n$-stacks that the tri-prestack of complexes $\mathfrak{T}^{[-2,0]}(\ES)$ with fibers $\rmT^{[-2,0]}(\ES/U)$ is a tristack. The characterization proposition cited above briefly says that $\mathfrak{P}$ is an $n$-stack over $\ES$ if and only if all descent data are effective and for any $X,Y$ objects of $\mathfrak{P}_U$, $\Hom_{\mathfrak{P}_U}(X,Y)$ is an $n-1$ stack over $\ES/U$.

\begin{remark}
The characterization proposition in \cite[Proposition 10.2]{hirschowitz-1998} is originally enounced for Segal $n$-categories, $n$-prestacks, and $n$-stacks. But again in the same paper, it has been remarked that the proposition holds for non-Segal structures \cite[\S 20]{hirschowitz-1998} where in this case, the weak structure is assumed to be the one defined by Tamsamani. Its definition can be found in \cite{simpson-1997} and \cite{MR1673923}. However, we are being very informal and not discussing here the connection of the weak structure of our categories, pre-stacks and, stacks with the ones mentioned above.
\end{remark}

Finally, we define the trihomomorphism of tristacks by localizing the triequivalence (\ref{functor:triequivalence}).

\begin{equation}\label{functor:triequivalence_stacky}
\xymatrix@1{\mathfrak{T}^{[-2,0]}(\ES) \ar[r] & \PCh(\ES)},
\end{equation}
where $\PCh(\ES)$ is considered naturally as a tristack. We deduce then its stack analog

\begin{theorem}\label{theorem:triequivalence_stacky}
(\ref{functor:triequivalence_stacky}) is a triequivalence of tristacks.
\end{theorem}

\appendix
\section{Appendix}\label{Annex}
We give the proof of Proposition (\ref{proposition:ab_group_embedded_in_2_category}).

We assume that the set $E$ is well ordered and denote the order on $E$ by $\preceq$. In what follows, we define
\begin{enumerate}
\item a 2-functor $F: \mathbb{Z}(E) \ra \twoC$,
\item \label{item:additive}for any two words $w_1$ and $w_2$ in $\mathbf{Z}(E)$, a functorial 1-morphism $\lambda_{w_1,w_2}$
\[\lambda_{w_1,w_2}:\xymatrix@1{F(w_1) \otimes F(w_2) \ar[r] & F(w_1+w_2)},\]
\item \label{item:associative}for any three words $w_1$, $w_2$, and $w_3$ in $\mathbf{Z}(E)$, a 2-morphism $\psi_{w_1,w_2,w_3}$ (\ref{diagram:2-morphism_of_assoc}),
\item \label{item:braiding}for any two words $w_1$ and $w_2$ in $\mathbf{Z}(E)$, a 2-morphism $\phi_{w_1,w_2}$ (\ref{diagram:2-morphism_of_braiding}).
\end{enumerate}

\subsection{Definition of $F$}

We construct the 2-functor $F:\mathbb{Z}(E) \ra \twoC$ as follows:
\begin{itemize}
\item For any generator $a \in E$, $Fa:=F_0a$,
\item For any generator $a \in E$, $F(-a):={(Fa)}^*$, where ${(Fa)}^*$ is inverse of $Fa$ in $\twoC$,
\item $F(0)$ is the unit element in $\twoC$, where 0 denotes the unit element in $\mathbb{Z}(E)$.
\item For any word $w$ in $\mathbb{Z}(E)$, we
        \begin{itemize}
        \item simplify $w$ so that there are no cancelations and denote the simplified word by $w_s$,
        \item order the letters of $w_s$ from least to greatest and denote the simplified and ordered word by $w_{s,o}$.
        \end{itemize}
      $F(w)$ is defined by multiplying the letters of $w_{s,o}$ from left to right.
\end{itemize}
For instance let $w=2a+b-c-a-2b$. After cancelations and ordering the letters $w_{s,o}=a-b-c$ and
\[F(w)=F(w_{s,o})=((Fa \otimes {(Fb)}^*) \otimes Fc).\]

The order on the set $E$ is needed since without the order two words that differ by the position of letters would map to different objects in $\twoC$ although they are the same word in $\mathbb{Z}(E)$. For the reasons of compactness, we use juxtaposition for the group operation $\otimes$ on the 2-category $\twoC$.

\subsection{Monoidal Case}
The items (\ref{item:additive}-(\ref{item:braiding}) describes the additive structure of the 2-functor $F$. We first define them on the words that do not have letters with negative coefficients. That is, they are constructed first on the free abelian monoid $\mathbb{N}(E)$. In Appendix \ref{subsection:group_case}, we extend their definitions to the free abelian group $\mathbb{Z}(E)$. We leave the verification of their compatibility with the Picard structure to the author's thesis \cite{Tatar_thesis}.

\paragraph{Definition of $\lambda_{w_1,w_2}$:} Let $w_1=a_1+ \ldots +a_m$ and $w_2=b_1+ \ldots +b_n$ be two words in $\mathbb{N}(E)$. The word $w_1+w_2$ is defined by concatenation of $w_1$ and $w_2$ and then by an $(m,n)$-shuffle so that the letters of $w_1$ and $w_2$ are ordered from least to greatest. We denote $w_1+w_2$ by $c_1+ \ldots +c_{m+n}$. From the definition of $F$,
\begin{eqnarray}
\label{diagram:words_separated}F(w_1) \otimes F(w_2)&=&(\ldots((Fa_1Fa_2)Fa_3)\ldots Fa_m) \otimes (\ldots((Fb_1Fb_2)Fb_3)\ldots Fb_n)\\
\label{diagram:words_shuffeled}F(w_1+w_2)&=&(\ldots((Fc_1Fc_2)Fc_3)\ldots Fc_{m+n})
\end{eqnarray}

We define the functorial morphism $\lambda_{w_1+w_2}:F(w_1) \otimes F(w_2) \ra F(w_1+w_2)$ in two steps as follows:
\begin{description}
\item[Step 1: Correct Bracketing] \mbox{}

        In this step, we define the morphism
        \begin{equation}\label{diagram:correct_bracketing}
        \begin{split}
        (\ldots((Fa_1Fa_2)Fa_3)\ldots Fa_m) \otimes (\ldots((Fb_1Fb_2)Fb_3)\ldots Fb_n){}& \ra \\ ((((\ldots((Fa_1Fa_2)Fa_3)\ldots Fa_m)Fb_1)Fb_2)\ldots Fb_n),
        \end{split}
        \end{equation}
        which moves the pairs of parenthesis of $F(w_2)$ one by one to the left from the outer most to the inner most without changing the place of parenthesis of $F(w_1)$. (\ref{diagram:correct_bracketing}) is composition of $n-1$ many morphisms of the form
        \begin{equation}\label{equation:shift}
        (\ldots((F(w_1)(F(w_2')Fb_i))Fb_{i+1}) \ldots Fb_{n}) \ra (\ldots(((F(w_1)F(w_2'))Fb_{i})Fb_{i+1}) \ldots Fb_{n}),
        \end{equation}
        where $w_2'$ is a subword of $w_2$.
\item[Step 2: Ordering Letters]\mbox{}

        Once the morphism (\ref{diagram:correct_bracketing}) is applied, the letters of $w_1$ and $w_2$ are parenthesized from left. Next, we define the morphism
        \begin{equation}\label{diagram:ordering_letters}
        ((((\ldots((Fa_1Fa_2)Fa_3)\ldots Fa_m)Fb_1)Fb_2)\ldots Fb_n) \ra (\ldots((Fc_1Fc_2)Fc_3)\ldots Fc_{m+n}),
        \end{equation}
        that shuffles the letters of $w_1$ and $w_2$ to order them from least to greatest, that is $c_1 \preceq c_2 \preceq \ldots \preceq c_{m+n}$.

        The rule is, find the smallest letter of $w_2$ in $w_1+w_2$ such that it has a letter of $w_1$ on its left that is greater, change their places. Depending on the position of the letters, there are two cases. Either the letters are in the same parenthesis , then (\ref{diagram:ordering_letters}) simply permutes them
        \begin{equation}\label{equation:swap}
        (\ldots((Fc_1Fc_2)Fc_3) \ldots Fc_{m+n}) \ra (\ldots((Fc_2Fc_1)Fc_3) \ldots Fc_{m+n}),
        \end{equation}
        or they are in different pairs of parenthesis and (\ref{diagram:ordering_letters}) first groups them together by moving the appropriate pair of parenthesis to the right, then permutes the letters, and moves the pair of parenthesis moved to the right to the left, that is
        \begin{equation}\label{equation:shift_swap_shif}
        \begin{split}
        ((\ldots((((\ldots(Fc_1Fc_2)\ldots)Fc_{k-1})Fc_{k+1})Fc_{k})\ldots )Fc_{m+n}){}& \ra \\
        ((\ldots(((\ldots(Fc_1Fc_2)\ldots)Fc_{k-1})(Fc_{k+1}Fc_{k}))\ldots )Fc_{m+n}){}& \ra \\
        ((\ldots(((\ldots(Fc_1Fc_2)\ldots)Fc_{k-1})(Fc_{k}Fc_{k+1}))\ldots )Fc_{m+n}){}& \ra \\
        ((\ldots((((\ldots(Fc_1Fc_2)\ldots)Fc_{k-1})Fc_k)Fc_{k+1})\ldots )Fc_{m+n})
        \end{split}
        \end{equation}
        where $c_k$ is a letter of $w_2$ in $w_1+w_2$ with $1 < k < m+n$ and $c_{k-1}$ is a letter of $w_1$ such that $c_k \prec c_{k-1}$.

        We repeat the above process to every letter of $w_2$ in $w_1+w_2$. We define the morphism (\ref{diagram:ordering_letters}) as composition of the morphisms of the form (\ref{equation:swap}) or (\ref{equation:shift_swap_shif}).

        We can illustrate the map (\ref{diagram:ordering_letters}) by the lattice paths \cite[Chapter 7.3D]{MR2394437}. It is clear that there is a 1-1 correspondence between the lattice paths from $(0,0)$ to $(m,n)$ and the $(m,n)$-shuffles. (\ref{diagram:words_shuffeled}) can be seen as the lattice path corresponding to the $(m,n)$-shuffle of the words $w_1,w_2$ that defines $w_1+w_2$ and (\ref{diagram:words_separated}) as the lattice path corresponding to the concatenation of the words $w_1$ and $w_2$ (i.e. the empty $(m,n)$-shuffle). We denote these paths by $L_{w_1+w_2}$ and $L_{w_1,w_2}$, respectively. From this perspective, the map (\ref{diagram:ordering_letters}) can be thought as applying an $(m,n)$-shuffle to the concatenation of the words $w_1$ and $w_2$.

        \unitlength 10mm
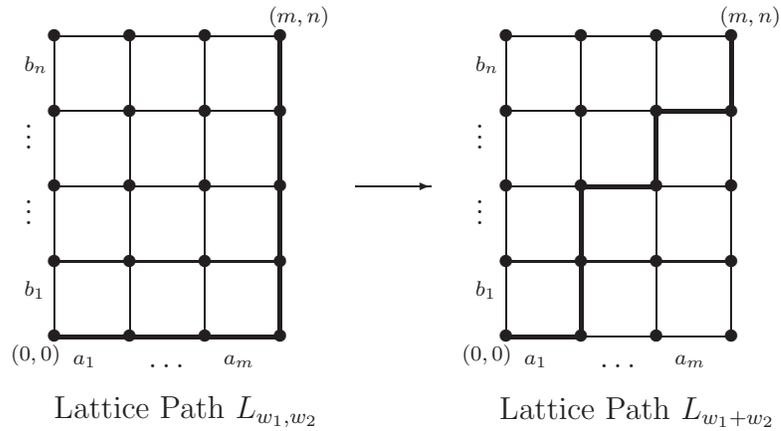
        \begin{picture}(7,5.5)
        \put(3.5,1){\line(1,0){3}}
        \put(3.5,2){\line(1,0){3}}
        \put(3.5,3){\line(1,0){3}}
        \put(3.5,4){\line(1,0){3}}
        \put(3.5,5){\line(1,0){3}}
        \put(3.5,1){\line(0,1){4}}
        \put(4.5,1){\line(0,1){4}}
        \put(5.5,1){\line(0,1){4}}
        \put(6.5,1){\line(0,1){4}}
        \put(3.25,0.75){\makebox(0,0)[cc]{\scriptsize$(0,0)$}}
        \put(6.75,5.25){\makebox(0,0)[cc]{\scriptsize$(m,n)$}}
        \put(3.5,1){\makebox(0,0)[cc]{$\bullet$}}
        \put(3.5,2){\makebox(0,0)[cc]{$\bullet$}}
        \put(3.5,3){\makebox(0,0)[cc]{$\bullet$}}
        \put(3.5,4){\makebox(0,0)[cc]{$\bullet$}}
        \put(3.5,5){\makebox(0,0)[cc]{$\bullet$}}
        \put(4.5,1){\makebox(0,0)[cc]{$\bullet$}}
        \put(4.5,2){\makebox(0,0)[cc]{$\bullet$}}
        \put(4.5,3){\makebox(0,0)[cc]{$\bullet$}}
        \put(4.5,4){\makebox(0,0)[cc]{$\bullet$}}
        \put(4.5,5){\makebox(0,0)[cc]{$\bullet$}}
        \put(5.5,1){\makebox(0,0)[cc]{$\bullet$}}
        \put(5.5,2){\makebox(0,0)[cc]{$\bullet$}}
        \put(5.5,3){\makebox(0,0)[cc]{$\bullet$}}
        \put(5.5,4){\makebox(0,0)[cc]{$\bullet$}}
        \put(5.5,5){\makebox(0,0)[cc]{$\bullet$}}
        \put(6.5,1){\makebox(0,0)[cc]{$\bullet$}}
        \put(6.5,2){\makebox(0,0)[cc]{$\bullet$}}
        \put(6.5,3){\makebox(0,0)[cc]{$\bullet$}}
        \put(6.5,4){\makebox(0,0)[cc]{$\bullet$}}
        \put(6.5,5){\makebox(0,0)[cc]{$\bullet$}}
        \put(3.75,0.6){\makebox(0,0)[lb]{\scriptsize$a_1$}}
        \put(4.75,0.6){\makebox(0,0)[lb]{\ldots}}
        \put(5.75,0.6){\makebox(0,0)[lb]{\scriptsize$a_m$}}
        \put(3.1,1.5){\makebox(0,0)[lb]{\scriptsize$b_1$}}
        \put(3.1,2.5){\makebox(0,0)[lb]{\vdots}}
        \put(3.1,3.5){\makebox(0,0)[lb]{\vdots}}
        \put(3.1,4.5){\makebox(0,0)[lb]{\scriptsize$b_n$}}
        \put(5.25,0){\makebox(0,0)[cc]{Lattice Path $L_{w_1,w_2}$}}
        \put(7.5,3){\vector(1,0){1}}
        \put(9.5,1){\line(1,0){3}}
        \put(9.5,2){\line(1,0){3}}
        \put(9.5,3){\line(1,0){3}}
        \put(9.5,4){\line(1,0){3}}
        \put(9.5,5){\line(1,0){3}}
        \put(9.5,1){\line(0,1){4}}
        \put(10.5,1){\line(0,1){4}}
        \put(11.5,1){\line(0,1){4}}
        \put(12.5,1){\line(0,1){4}}
        \put(9.25,0.75){\makebox(0,0)[cc]{\scriptsize$(0,0)$}}
        \put(12.75,5.25){\makebox(0,0)[cc]{\scriptsize$(m,n)$}}
        \put(9.5,1){\makebox(0,0)[cc]{$\bullet$}}
        \put(9.5,2){\makebox(0,0)[cc]{$\bullet$}}
        \put(9.5,3){\makebox(0,0)[cc]{$\bullet$}}
        \put(9.5,4){\makebox(0,0)[cc]{$\bullet$}}
        \put(9.5,5){\makebox(0,0)[cc]{$\bullet$}}
        \put(10.5,1){\makebox(0,0)[cc]{$\bullet$}}
        \put(10.5,2){\makebox(0,0)[cc]{$\bullet$}}
        \put(10.5,3){\makebox(0,0)[cc]{$\bullet$}}
        \put(10.5,4){\makebox(0,0)[cc]{$\bullet$}}
        \put(10.5,5){\makebox(0,0)[cc]{$\bullet$}}
        \put(11.5,1){\makebox(0,0)[cc]{$\bullet$}}
        \put(11.5,2){\makebox(0,0)[cc]{$\bullet$}}
        \put(11.5,3){\makebox(0,0)[cc]{$\bullet$}}
        \put(11.5,4){\makebox(0,0)[cc]{$\bullet$}}
        \put(11.5,5){\makebox(0,0)[cc]{$\bullet$}}
        \put(12.5,1){\makebox(0,0)[cc]{$\bullet$}}
        \put(12.5,2){\makebox(0,0)[cc]{$\bullet$}}
        \put(12.5,3){\makebox(0,0)[cc]{$\bullet$}}
        \put(12.5,4){\makebox(0,0)[cc]{$\bullet$}}
        \put(12.5,5){\makebox(0,0)[cc]{$\bullet$}}
        \put(9.75,0.6){\makebox(0,0)[lb]{\scriptsize$a_1$}}
        \put(10.75,0.6){\makebox(0,0)[lb]{\ldots}}
        \put(11.75,0.6){\makebox(0,0)[lb]{\scriptsize$a_m$}}
        \put(11.25,0){\makebox(0,0)[cc]{Lattice Path $L_{w_1+w_2}$}}
        \put(9.1,1.5){\makebox(0,0)[lb]{\scriptsize$b_1$}}
        \put(9.1,2.5){\makebox(0,0)[lb]{\vdots}}
        \put(9.1,3.5){\makebox(0,0)[lb]{\vdots}}
        \put(9.1,4.5){\makebox(0,0)[lb]{\scriptsize$b_n$}}
        \linethickness{.5mm}
        \put(3.5,1){\line(1,0){1}}
        \put(4.5,1){\line(1,0){1}}
        \put(5.5,1){\line(1,0){1}}
        \put(6.5,1){\line(0,1){1}}
        \put(6.5,2){\line(0,1){1}}
        \put(6.5,3){\line(0,1){1}}
        \put(6.5,4){\line(0,1){1}}
        \put(9.5,1){\line(1,0){1}}
        \put(10.5,1){\line(0,1){1}}
        \put(10.5,2){\line(0,1){1}}
        \put(10.5,3){\line(1,0){1}}
        \put(11.5,3){\line(0,1){1}}
        \put(11.5,4){\line(1,0){1}}
        \put(12.5,4){\line(0,1){1}}
        \end{picture}

        \vspace{0.5cm}

        The morphisms (\ref{equation:swap}) and (\ref{equation:shift_swap_shif}) describe the basic movement. They substitute the point $(i,j)$ on the lattice path with the point $(i-1,j+1)$ as shown in the picture below.

        \unitlength 10mm
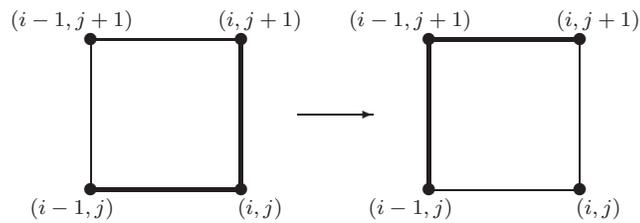
        \begin{picture}(7,2.5)
        \put(4.5,0){\line(1,0){2}}
        \put(4.5,0){\line(0,1){2}}
        \put(6.5,0){\line(0,1){2}}
        \put(4.5,2){\line(1,0){2}}
        \put(9,0){\line(1,0){2}}
        \put(9,0){\line(0,1){2}}
        \put(11,0){\line(0,1){2}}
        \put(9,2){\line(1,0){2}}
        \put(7.25,1){\vector(1,0){1}}
        \put(4.25,2.25){\makebox(0,0)[cc]{\scriptsize$(i-1,j+1)$}}
        \put(4.25,-0.25){\makebox(0,0)[cc]{\scriptsize$(i-1,j)$}}
        \put(6.75,2.25){\makebox(0,0)[cc]{\scriptsize$(i,j+1)$}}
        \put(6.75,-0.25){\makebox(0,0)[cc]{\scriptsize$(i,j)$}}
        \put(8.75,2.25){\makebox(0,0)[cc]{\scriptsize$(i-1,j+1)$}}
        \put(8.75,-0.25){\makebox(0,0)[cc]{\scriptsize$(i-1,j)$}}
        \put(11.25,2.25){\makebox(0,0)[cc]{\scriptsize$(i,j+1)$}}
        \put(11.25,-0.25){\makebox(0,0)[cc]{\scriptsize$(i,j)$}}
        \put(4.5,0){\makebox(0,0)[cc]{$\bullet$}}
        \put(6.5,0){\makebox(0,0)[cc]{$\bullet$}}
        \put(4.5,2){\makebox(0,0)[cc]{$\bullet$}}
        \put(6.5,2){\makebox(0,0)[cc]{$\bullet$}}
        \put(9,0){\makebox(0,0)[cc]{$\bullet$}}
        \put(11,0){\makebox(0,0)[cc]{$\bullet$}}
        \put(9,2){\makebox(0,0)[cc]{$\bullet$}}
        \put(11,2){\makebox(0,0)[cc]{$\bullet$}}
        \linethickness{0.5mm}
        \put(4.5,0){\line(1,0){2}}
        \put(6.5,0){\line(0,1){2}}
        \put(9,0){\line(0,1){2}}
        \put(9,2){\line(1,0){2}}
        \end{picture}

        \vspace{0.5cm}

        The overall movement is described by the morphism (\ref{diagram:ordering_letters}) where each step is a basic movement. We define the following special point on the lattice path in order to explain the mechanism of the movements. We call the point $(i,j)$ on the lattice path the \emph{corner point} if the points $(i-1,j)$ and $(i,j+1)$ are on the lattice path, as well. The morphism (\ref{diagram:ordering_letters}) picks at every step the corner point $(i,j)$ with the least $y$-coordinate that is not on the lattice path $L_{w_1+w_2}$ and substitutes it with $(i-1,j+1)$. We show in the picture below the transformation of the lattice path $L_{w_1,w_2}$ to the lattice path $L_{w_1+w_2}$.

        \unitlength 10mm
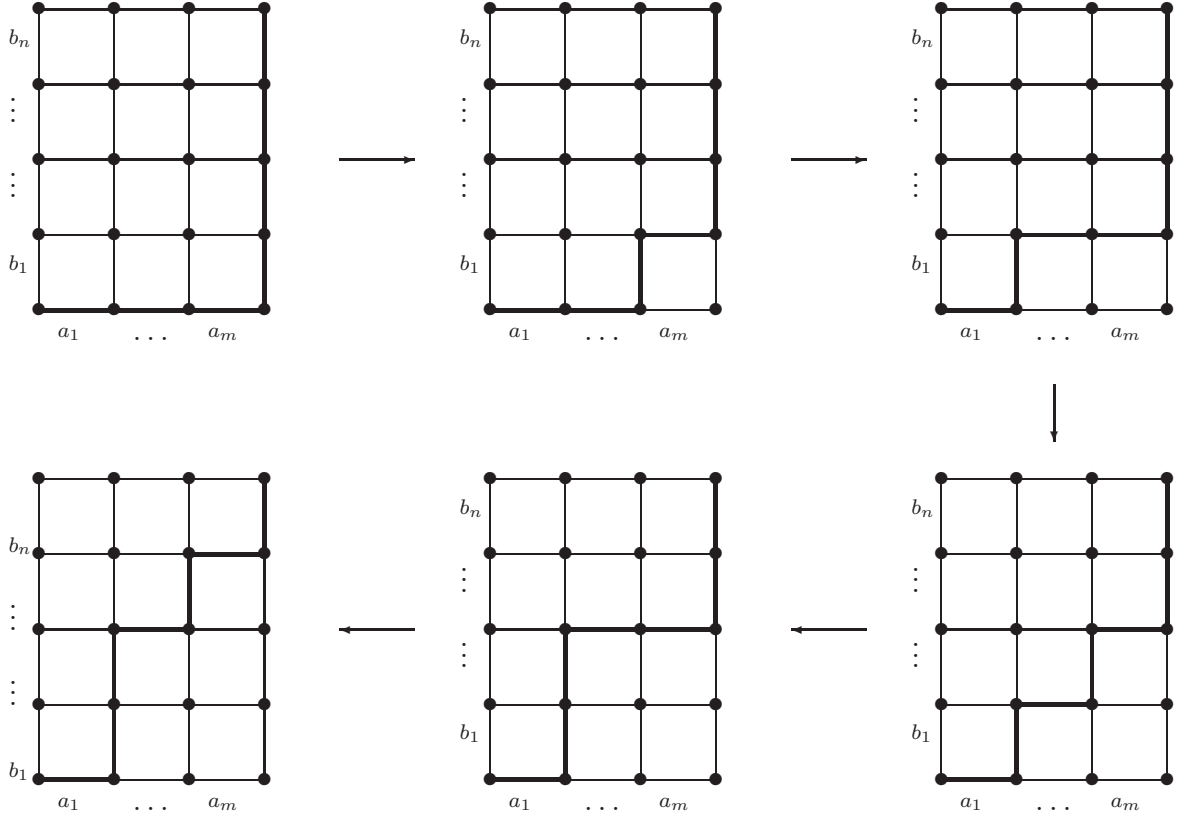
        \begin{picture}(7,12.5)
        \put(0.5,7.25){\line(1,0){3}}
        \put(0.5,8.25){\line(1,0){3}}
        \put(0.5,9.25){\line(1,0){3}}
        \put(0.5,10.25){\line(1,0){3}}
        \put(0.5,11.25){\line(1,0){3}}
        \put(0.5,7.25){\line(0,1){4}}
        \put(1.5,7.25){\line(0,1){4}}
        \put(2.5,7.25){\line(0,1){4}}
        \put(3.5,7.25){\line(0,1){4}}
        \put(0.5,7.25){\makebox(0,0)[cc]{$\bullet$}}
        \put(0.5,8.25){\makebox(0,0)[cc]{$\bullet$}}
        \put(0.5,9.25){\makebox(0,0)[cc]{$\bullet$}}
        \put(0.5,10.25){\makebox(0,0)[cc]{$\bullet$}}
        \put(0.5,11.25){\makebox(0,0)[cc]{$\bullet$}}
        \put(1.5,7.25){\makebox(0,0)[cc]{$\bullet$}}
        \put(1.5,8.25){\makebox(0,0)[cc]{$\bullet$}}
        \put(1.5,9.25){\makebox(0,0)[cc]{$\bullet$}}
        \put(1.5,10.25){\makebox(0,0)[cc]{$\bullet$}}
        \put(1.5,11.25){\makebox(0,0)[cc]{$\bullet$}}
        \put(2.5,7.25){\makebox(0,0)[cc]{$\bullet$}}
        \put(2.5,8.25){\makebox(0,0)[cc]{$\bullet$}}
        \put(2.5,9.25){\makebox(0,0)[cc]{$\bullet$}}
        \put(2.5,10.25){\makebox(0,0)[cc]{$\bullet$}}
        \put(2.5,11.25){\makebox(0,0)[cc]{$\bullet$}}
        \put(3.5,7.25){\makebox(0,0)[cc]{$\bullet$}}
        \put(3.5,8.25){\makebox(0,0)[cc]{$\bullet$}}
        \put(3.5,9.25){\makebox(0,0)[cc]{$\bullet$}}
        \put(3.5,10.25){\makebox(0,0)[cc]{$\bullet$}}
        \put(3.5,11.25){\makebox(0,0)[cc]{$\bullet$}}
        \put(0.75,6.85){\makebox(0,0)[lb]{\scriptsize$a_1$}}
        \put(1.75,6.85){\makebox(0,0)[lb]{\ldots}}
        \put(2.75,6.85){\makebox(0,0)[lb]{\scriptsize$a_m$}}
        \put(0.1,7.75){\makebox(0,0)[lb]{\scriptsize$b_1$}}
        \put(0.1,8.75){\makebox(0,0)[lb]{\vdots}}
        \put(0.1,9.75){\makebox(0,0)[lb]{\vdots}}
        \put(0.1,10.75){\makebox(0,0)[lb]{\scriptsize$b_n$}}
        \put(4.5,9.25){\vector(1,0){1}}
        \put(6.5,7.25){\line(1,0){3}}
        \put(6.5,8.25){\line(1,0){3}}
        \put(6.5,9.25){\line(1,0){3}}
        \put(6.5,10.25){\line(1,0){3}}
        \put(6.5,11.25){\line(1,0){3}}
        \put(6.5,7.25){\line(0,1){4}}
        \put(7.5,7.25){\line(0,1){4}}
        \put(8.5,7.25){\line(0,1){4}}
        \put(9.5,7.25){\line(0,1){4}}
        \put(6.5,7.25){\makebox(0,0)[cc]{$\bullet$}}
        \put(6.5,8.25){\makebox(0,0)[cc]{$\bullet$}}
        \put(6.5,9.25){\makebox(0,0)[cc]{$\bullet$}}
        \put(6.5,10.25){\makebox(0,0)[cc]{$\bullet$}}
        \put(6.5,11.25){\makebox(0,0)[cc]{$\bullet$}}
        \put(7.5,7.25){\makebox(0,0)[cc]{$\bullet$}}
        \put(7.5,8.25){\makebox(0,0)[cc]{$\bullet$}}
        \put(7.5,9.25){\makebox(0,0)[cc]{$\bullet$}}
        \put(7.5,10.25){\makebox(0,0)[cc]{$\bullet$}}
        \put(7.5,11.25){\makebox(0,0)[cc]{$\bullet$}}
        \put(8.5,7.25){\makebox(0,0)[cc]{$\bullet$}}
        \put(8.5,8.25){\makebox(0,0)[cc]{$\bullet$}}
        \put(8.5,9.25){\makebox(0,0)[cc]{$\bullet$}}
        \put(8.5,10.25){\makebox(0,0)[cc]{$\bullet$}}
        \put(8.5,11.25){\makebox(0,0)[cc]{$\bullet$}}
        \put(9.5,7.25){\makebox(0,0)[cc]{$\bullet$}}
        \put(9.5,8.25){\makebox(0,0)[cc]{$\bullet$}}
        \put(9.5,9.25){\makebox(0,0)[cc]{$\bullet$}}
        \put(9.5,10.25){\makebox(0,0)[cc]{$\bullet$}}
        \put(9.5,11.25){\makebox(0,0)[cc]{$\bullet$}}
        \put(6.75,6.85){\makebox(0,0)[lb]{\scriptsize$a_1$}}
        \put(7.75,6.85){\makebox(0,0)[lb]{\ldots}}
        \put(8.75,6.85){\makebox(0,0)[lb]{\scriptsize$a_m$}}
        \put(6.1,7.75){\makebox(0,0)[lb]{\scriptsize$b_1$}}
        \put(6.1,8.75){\makebox(0,0)[lb]{\vdots}}
        \put(6.1,9.75){\makebox(0,0)[lb]{\vdots}}
        \put(6.1,10.75){\makebox(0,0)[lb]{\scriptsize$b_n$}}
        \put(10.5,9.25){\vector(1,0){1}}
        \put(12.5,7.25){\line(1,0){3}}
        \put(12.5,8.25){\line(1,0){3}}
        \put(12.5,9.25){\line(1,0){3}}
        \put(12.5,10.25){\line(1,0){3}}
        \put(12.5,11.25){\line(1,0){3}}
        \put(12.5,7.25){\line(0,1){4}}
        \put(13.5,7.25){\line(0,1){4}}
        \put(14.5,7.25){\line(0,1){4}}
        \put(15.5,7.25){\line(0,1){4}}
        \put(12.5,7.25){\makebox(0,0)[cc]{$\bullet$}}
        \put(12.5,8.25){\makebox(0,0)[cc]{$\bullet$}}
        \put(12.5,9.25){\makebox(0,0)[cc]{$\bullet$}}
        \put(12.5,10.25){\makebox(0,0)[cc]{$\bullet$}}
        \put(12.5,11.25){\makebox(0,0)[cc]{$\bullet$}}
        \put(13.5,7.25){\makebox(0,0)[cc]{$\bullet$}}
        \put(13.5,8.25){\makebox(0,0)[cc]{$\bullet$}}
        \put(13.5,9.25){\makebox(0,0)[cc]{$\bullet$}}
        \put(13.5,10.25){\makebox(0,0)[cc]{$\bullet$}}
        \put(13.5,11.25){\makebox(0,0)[cc]{$\bullet$}}
        \put(14.5,7.25){\makebox(0,0)[cc]{$\bullet$}}
        \put(14.5,8.25){\makebox(0,0)[cc]{$\bullet$}}
        \put(14.5,9.25){\makebox(0,0)[cc]{$\bullet$}}
        \put(14.5,10.25){\makebox(0,0)[cc]{$\bullet$}}
        \put(14.5,11.25){\makebox(0,0)[cc]{$\bullet$}}
        \put(15.5,7.25){\makebox(0,0)[cc]{$\bullet$}}
        \put(15.5,8.25){\makebox(0,0)[cc]{$\bullet$}}
        \put(15.5,9.25){\makebox(0,0)[cc]{$\bullet$}}
        \put(15.5,10.25){\makebox(0,0)[cc]{$\bullet$}}
        \put(15.5,11.25){\makebox(0,0)[cc]{$\bullet$}}
        \put(12.75,6.85){\makebox(0,0)[lb]{\scriptsize$a_1$}}
        \put(13.75,6.85){\makebox(0,0)[lb]{\ldots}}
        \put(14.75,6.85){\makebox(0,0)[lb]{\scriptsize$a_m$}}
        \put(12.1,7.75){\makebox(0,0)[lb]{\scriptsize$b_1$}}
        \put(12.1,8.75){\makebox(0,0)[lb]{\vdots}}
        \put(12.1,9.75){\makebox(0,0)[lb]{\vdots}}
        \put(12.1,10.75){\makebox(0,0)[lb]{\scriptsize$b_n$}}
        \put(14,6.25){\vector(0,-1){0.75}}
        \linethickness{.5mm}
        \put(0.5,7.25){\line(1,0){1}}
        \put(1.5,7.25){\line(1,0){1}}
        \put(2.5,7.25){\line(1,0){1}}
        \put(3.5,7.25){\line(0,1){1}}
        \put(3.5,8.25){\line(0,1){1}}
        \put(3.5,9.25){\line(0,1){1}}
        \put(3.5,10.25){\line(0,1){1}}
        \put(6.5,7.25){\line(1,0){1}}
        \put(7.5,7.25){\line(1,0){1}}
        \put(8.5,7.25){\line(0,1){1}}
        \put(8.5,8.25){\line(1,0){1}}
        \put(9.5,8.25){\line(0,1){1}}
        \put(9.5,9.25){\line(0,1){1}}
        \put(9.5,10.25){\line(0,1){1}}
        \put(12.5,7.25){\line(1,0){1}}
        \put(13.5,7.25){\line(0,1){1}}
        \put(13.5,8.25){\line(1,0){1}}
        \put(14.5,8.25){\line(1,0){1}}
        \put(15.5,8.25){\line(0,1){1}}
        \put(15.5,9.25){\line(0,1){1}}
        \put(15.5,10.25){\line(0,1){1}}
        \thinlines
        \put(0.5,1){\line(1,0){3}}
        \put(0.5,2){\line(1,0){3}}
        \put(0.5,3){\line(1,0){3}}
        \put(0.5,4){\line(1,0){3}}
        \put(0.5,5){\line(1,0){3}}
        \put(0.5,1){\line(0,1){4}}
        \put(1.5,1){\line(0,1){4}}
        \put(2.5,1){\line(0,1){4}}
        \put(3.5,1){\line(0,1){4}}
        \put(0.5,1){\makebox(0,0)[cc]{$\bullet$}}
        \put(0.5,2){\makebox(0,0)[cc]{$\bullet$}}
        \put(0.5,3){\makebox(0,0)[cc]{$\bullet$}}
        \put(0.5,4){\makebox(0,0)[cc]{$\bullet$}}
        \put(0.5,5){\makebox(0,0)[cc]{$\bullet$}}
        \put(1.5,1){\makebox(0,0)[cc]{$\bullet$}}
        \put(1.5,2){\makebox(0,0)[cc]{$\bullet$}}
        \put(1.5,3){\makebox(0,0)[cc]{$\bullet$}}
        \put(1.5,4){\makebox(0,0)[cc]{$\bullet$}}
        \put(1.5,5){\makebox(0,0)[cc]{$\bullet$}}
        \put(2.5,1){\makebox(0,0)[cc]{$\bullet$}}
        \put(2.5,2){\makebox(0,0)[cc]{$\bullet$}}
        \put(2.5,3){\makebox(0,0)[cc]{$\bullet$}}
        \put(2.5,4){\makebox(0,0)[cc]{$\bullet$}}
        \put(2.5,5){\makebox(0,0)[cc]{$\bullet$}}
        \put(3.5,1){\makebox(0,0)[cc]{$\bullet$}}
        \put(3.5,2){\makebox(0,0)[cc]{$\bullet$}}
        \put(3.5,3){\makebox(0,0)[cc]{$\bullet$}}
        \put(3.5,4){\makebox(0,0)[cc]{$\bullet$}}
        \put(3.5,5){\makebox(0,0)[cc]{$\bullet$}}
        \put(0.75,0.6){\makebox(0,0)[lb]{\scriptsize$a_1$}}
        \put(1.75,0.6){\makebox(0,0)[lb]{\ldots}}
        \put(2.75,0.6){\makebox(0,0)[lb]{\scriptsize$a_m$}}
        \put(0.1,1){\makebox(0,0)[lb]{\scriptsize$b_1$}}
        \put(0.1,2){\makebox(0,0)[lb]{\vdots}}
        \put(0.1,3){\makebox(0,0)[lb]{\vdots}}
        \put(0.1,4){\makebox(0,0)[lb]{\scriptsize$b_n$}}
        \put(5.5,3){\vector(-1,0){1}}
        \put(6.5,1){\line(1,0){3}}
        \put(6.5,2){\line(1,0){3}}
        \put(6.5,3){\line(1,0){3}}
        \put(6.5,4){\line(1,0){3}}
        \put(6.5,5){\line(1,0){3}}
        \put(6.5,1){\line(0,1){4}}
        \put(7.5,1){\line(0,1){4}}
        \put(8.5,1){\line(0,1){4}}
        \put(9.5,1){\line(0,1){4}}
        \put(6.5,1){\makebox(0,0)[cc]{$\bullet$}}
        \put(6.5,2){\makebox(0,0)[cc]{$\bullet$}}
        \put(6.5,3){\makebox(0,0)[cc]{$\bullet$}}
        \put(6.5,4){\makebox(0,0)[cc]{$\bullet$}}
        \put(6.5,5){\makebox(0,0)[cc]{$\bullet$}}
        \put(7.5,1){\makebox(0,0)[cc]{$\bullet$}}
        \put(7.5,2){\makebox(0,0)[cc]{$\bullet$}}
        \put(7.5,3){\makebox(0,0)[cc]{$\bullet$}}
        \put(7.5,4){\makebox(0,0)[cc]{$\bullet$}}
        \put(7.5,5){\makebox(0,0)[cc]{$\bullet$}}
        \put(8.5,1){\makebox(0,0)[cc]{$\bullet$}}
        \put(8.5,2){\makebox(0,0)[cc]{$\bullet$}}
        \put(8.5,3){\makebox(0,0)[cc]{$\bullet$}}
        \put(8.5,4){\makebox(0,0)[cc]{$\bullet$}}
        \put(8.5,5){\makebox(0,0)[cc]{$\bullet$}}
        \put(9.5,1){\makebox(0,0)[cc]{$\bullet$}}
        \put(9.5,2){\makebox(0,0)[cc]{$\bullet$}}
        \put(9.5,3){\makebox(0,0)[cc]{$\bullet$}}
        \put(9.5,4){\makebox(0,0)[cc]{$\bullet$}}
        \put(9.5,5){\makebox(0,0)[cc]{$\bullet$}}
        \put(6.75,0.6){\makebox(0,0)[lb]{\scriptsize$a_1$}}
        \put(7.75,0.6){\makebox(0,0)[lb]{\ldots}}
        \put(8.75,0.6){\makebox(0,0)[lb]{\scriptsize$a_m$}}
        \put(6.1,1.5){\makebox(0,0)[lb]{\scriptsize$b_1$}}
        \put(6.1,2.5){\makebox(0,0)[lb]{\vdots}}
        \put(6.1,3.5){\makebox(0,0)[lb]{\vdots}}
        \put(6.1,4.5){\makebox(0,0)[lb]{\scriptsize$b_n$}}
        \put(11.5,3){\vector(-1,0){1}}
        \put(12.5,1){\line(1,0){3}}
        \put(12.5,2){\line(1,0){3}}
        \put(12.5,3){\line(1,0){3}}
        \put(12.5,4){\line(1,0){3}}
        \put(12.5,5){\line(1,0){3}}
        \put(12.5,1){\line(0,1){4}}
        \put(13.5,1){\line(0,1){4}}
        \put(14.5,1){\line(0,1){4}}
        \put(15.5,1){\line(0,1){4}}
        \put(12.5,1){\makebox(0,0)[cc]{$\bullet$}}
        \put(12.5,2){\makebox(0,0)[cc]{$\bullet$}}
        \put(12.5,3){\makebox(0,0)[cc]{$\bullet$}}
        \put(12.5,4){\makebox(0,0)[cc]{$\bullet$}}
        \put(12.5,5){\makebox(0,0)[cc]{$\bullet$}}
        \put(13.5,1){\makebox(0,0)[cc]{$\bullet$}}
        \put(13.5,2){\makebox(0,0)[cc]{$\bullet$}}
        \put(13.5,3){\makebox(0,0)[cc]{$\bullet$}}
        \put(13.5,4){\makebox(0,0)[cc]{$\bullet$}}
        \put(13.5,5){\makebox(0,0)[cc]{$\bullet$}}
        \put(14.5,1){\makebox(0,0)[cc]{$\bullet$}}
        \put(14.5,2){\makebox(0,0)[cc]{$\bullet$}}
        \put(14.5,3){\makebox(0,0)[cc]{$\bullet$}}
        \put(14.5,4){\makebox(0,0)[cc]{$\bullet$}}
        \put(14.5,5){\makebox(0,0)[cc]{$\bullet$}}
        \put(15.5,1){\makebox(0,0)[cc]{$\bullet$}}
        \put(15.5,2){\makebox(0,0)[cc]{$\bullet$}}
        \put(15.5,3){\makebox(0,0)[cc]{$\bullet$}}
        \put(15.5,4){\makebox(0,0)[cc]{$\bullet$}}
        \put(15.5,5){\makebox(0,0)[cc]{$\bullet$}}
        \put(12.75,0.6){\makebox(0,0)[lb]{\scriptsize$a_1$}}
        \put(13.75,0.6){\makebox(0,0)[lb]{\ldots}}
        \put(14.75,0.6){\makebox(0,0)[lb]{\scriptsize$a_m$}}
        \put(12.1,1.5){\makebox(0,0)[lb]{\scriptsize$b_1$}}
        \put(12.1,2.5){\makebox(0,0)[lb]{\vdots}}
        \put(12.1,3.5){\makebox(0,0)[lb]{\vdots}}
        \put(12.1,4.5){\makebox(0,0)[lb]{\scriptsize$b_n$}}
        \linethickness{.5mm}
        \put(0.5,1){\line(1,0){1}}
        \put(1.5,1){\line(0,1){1}}
        \put(1.5,2){\line(0,1){1}}
        \put(1.5,3){\line(1,0){1}}
        \put(2.5,3){\line(0,1){1}}
        \put(2.5,4){\line(1,0){1}}
        \put(3.5,4){\line(0,1){1}}
        \put(6.5,1){\line(1,0){1}}
        \put(7.5,1){\line(0,1){1}}
        \put(7.5,2){\line(0,1){1}}
        \put(7.5,3){\line(1,0){1}}
        \put(8.5,3){\line(1,0){1}}
        \put(9.5,3){\line(0,1){1}}
        \put(9.5,4){\line(0,1){1}}
        \put(12.5,1){\line(1,0){1}}
        \put(13.5,1){\line(0,1){1}}
        \put(13.5,2){\line(1,0){1}}
        \put(14.5,2){\line(0,1){1}}
        \put(14.5,3){\line(1,0){1}}
        \put(15.5,3){\line(0,1){1}}
        \put(15.5,4){\line(0,1){1}}
        \end{picture}
\end{description}

The morphism (\ref{diagram:correct_bracketing}) obtained in the first step followed by the morphism (\ref{diagram:ordering_letters}) constructed in the second step defines $\lambda_{w_1,w_2}$.

We remark that if all the letters of $w_1$ are less than all the letters of $w_2$, then $w_1+w_2$ is obtained by concatenating the words $w_1$ and $w_2$ without the shuffle. That is $L_{w_1+w_2}$ coincides with $L_{w_1,w_2}$. In this case $\lambda_{w_1,w_2}$ is of the form (\ref{diagram:correct_bracketing}).

We also observe that the morphism $\lambda_{w_1,w_2}$ is a path in the 1-skeleton of permuto-associahedron $K\Pi_{m+n-1}$ where $m$ and $n$ are lengths of the words $w_1$ and $w_2$, respectively. $K\Pi_{m+n-1}$ is a  polytope whose vertices are all possible orderings and groupings of strings of length $m+n$ and whose edges are all possible adjacent permutations and all possible parenthesis movements. For more details about permuto-associahedron, we refer to \cite{MR1207505} and \cite{MR1311028}.

\paragraph{Definition of $\psi_{w_1,w_2,w_3}$:}
For any three words $w_1,w_2,w_3$ in $\mathbb{N}(E)$, we define the 2-morphism $\psi_{w_1,w_2,w_3}$

\begin{equation}\label{diagram:2-morphism_of_assoc}
\xymatrix{((F(w_1)F(w_2)) F(w_3)) \ar[d]_{\mathsf{a}} \ar[r]^(0.55){\lambda_{w_1,w_2}} & F(w_1+w_2) F(w_3) \ar[r]^(0.65){\lambda_{w_1+w_2,w_3}} \drtwocell<\omit>{\hspace{1.2cm}\psi_{w_1,w_2,w_3}}& **[r]F(w_1+w_2+w_3) \ar@{=}[d]\\
            (F(w_1)(F(w_2)F(w_3))) \ar[r]_(0.55){\lambda_{w_2,w_3}} & F(w_1) F(w_2+w_3) \ar[r]_(0.65){\lambda_{w_1,w_2+w_3}} & **[r] F(w_1+w_2+w_3)}
\end{equation}
between the 1-morphisms $\lambda_{w_1,w_2+w_3} \circ \lambda_{w_2,w_3} \circ \mathsf{a}$ and $\lambda_{w_1+w_2,w_3} \circ \lambda_{w_1,w_2}$ from $((F(w_1)F(w_2))F(w_3))$ to $F(w_1+w_2+w_3)$\footnote[1]{We commit an abuse of notation in diagram (\ref{diagram:2-morphism_of_assoc}). By $\lambda_{w_1,w_2}$ and $\lambda_{w_2,w_3}$ we mean $\lambda_{w_1,w_2} \otimes \id_{w_3}$ and $\id_{w_1} \otimes \lambda_{w_2,w_3}$, respectively.}. These 1-morphisms are paths in the 1-skeleton of $K\Pi_{m+n+p-1}$ where $n$,$m$, and $p$ are the lengths of the words $w_1$, $w_2$, and $w_3$, respectively. This follows from the fact that every map in the diagram (\ref{diagram:2-morphism_of_assoc}) is in the 1-skeleton of $K\Pi_{m+n+p-1}$.

In order to better understand these paths, we interpret them in terms of 3-dimensional lattice paths. Assume that the letters of the words $w_1$, $w_2$, and $w_3$ represent respectively the unit intervals on the $x$-, $y$-, and $z$-axis. $F(w_1+w_2+w_3)$ can be represented by the 3-dimensional lattice path corresponding to the $(m,n,p)$-shuffle of the words $w_1,w_2,w_3$ that defines $w_1+w_2+w_3$ and  $((F(w_1)(F(w_2))F(w_3))$ by the 3-dimensional lattice path corresponding to the empty shuffle of the words $w_1,w_2,w_3$. Therefore, the paths $\lambda_{w_1,w_2+w_3} \circ \lambda_{w_2,w_3} \circ \mathsf{a}$ and $\lambda_{w_1+w_2,w_3} \circ \lambda_{w_1,w_2}$ can be thought as two different ways of shuffling $w_1, w_2, w_3$ to obtain $w_1+w_2+w_3$. The path $\lambda_{w_1,w_2+w_3} \circ \lambda_{w_2,w_3} \circ \mathsf{a}$ first does the $(n,p)$-shuffle then the $(m,n)$-shuffle. On the other hand the path $\lambda_{w_1+w_2,w_3} \circ \lambda_{w_1,w_2}$ does the $(m,n)$-shuffle first, then the $(n,p)$-shuffle. In this sense the 2-morphism $\psi_{w_1,w_2,w_3}$ can be seen as the connection between the two different ways of doing the $(m,n,p)$-shuffle.

To define the 2-morphism $\psi_{w_1,w_2,w_3}$, we need the following lemmas .

\begin{lemma}\label{lemma:path_equality}
Let $w_1$ and $w_2$ be two elements of $\mathbb{N}(E)$. $\lambda_{w_2,w_3}=\mathsf{c}$ and $\lambda_{w_1,w_2}=\id$ if and only if $\lambda_{w_1,w_2+w_3} \circ \lambda_{w_2,w_3} \circ \mathsf{a}=\lambda_{w_1+w_2,w_3} \circ \lambda_{w_1,w_2}$
\end{lemma}

\begin{proof}
We first remark that $\lambda_{w_2,w_3}=\mathsf{c}$ and $\lambda_{w_1,w_2}=\id$ is equivalent to assuming $w_2$ and $w_3$ are letters such that $w_2$ is greater than $w_3$ and $w_2$ is greater than or equal to all letters of $w_1$. These facts imply that the map $\lambda_{w_1+w_2,w_3}$ first permutes $F(w_2)$ and $F(w_3)$ then shuffles $F(w_1)$ and $F(w_3)$ without changing the position of $F(w_2)$. Thus $\lambda_{w_1,w_2+w_3} \circ \lambda_{w_2,w_3} \circ \mathsf{a}=\lambda_{w_1+w_2,w_3} \circ \lambda_{w_1,w_2}$.

In the other direction, we observe that the morphism $\mathsf{a}$ can be only part of the morphism $\lambda_{w_1,w_2+w_3}$ which means $\lambda_{w_1,w_2}=\id$. This requires $w_2$ to be a letter greater than or equal to all letters of $w_1$ and $\lambda_{w_1,w_2+w_3} \circ \lambda_{w_2,w_3} \circ \mathsf{a}=\lambda_{w_1+w_2,w_3}$. We also observe that a parenthesis movement caused by $\lambda_{w_2,w_3}$ effects only the places of the parenthesis around the letters of $w_2$ and $w_3$ and such a movement cannot be caused by $\lambda_{w_1+w_2,w_3}$. This means $\lambda_{w_2,w_3}$ does not cause any parenthesis movements. Hence, we deduce that $w_3$ is also a letter. If $w_2 \preceq w_3$ then $\lambda_{w_2,w_3}$ and $\lambda_{w_1+w_2,w_3}$ become identity morphisms and we obtain $\lambda_{w_1,w_2+w_3} \circ \mathsf{a}=\id$ which is not possible. Therefore $\lambda_{w_2,w_3}$ should consist of a single permutation.
\end{proof}

\begin{lemma}\label{lemma:path_inclusion}
   Let $w_1$, $w_2$, and $w_3$ be three elements of $\mathbb{N}(E)$. Then the followings are equivalent.
    \begin{enumerate}
    \item\label{lemma:path_inclusion_1} The path $\lambda_{w_1+w_2,w_3} \circ \lambda_{w_1,w_2}$ is strictly included in $\lambda_{w_1,w_2+w_3} \circ \lambda_{w_2,w_3} \circ \mathsf{a}$. That is $V_{(w_1,w_2|w_3)}$ the vertex set of the path $\lambda_{w_1+w_2,w_3} \circ \lambda_{w_1,w_2}$ is strictly included in $V_{(w_1|w_2,w_3)}$ the vertex set of the path $\lambda_{w_1,w_2+w_3} \circ \lambda_{w_2,w_3} \circ \mathsf{a}$.
    \item\label{lemma:path_inclusion_2} $\lambda_{w_1,w_2+w_3} \circ \lambda_{w_2,w_3}= \lambda_{w_1+w_2,w_3} \circ \lambda_{w_1,w_2} \circ \mathsf{a}^{-1}$.
    \item\label{lemma:path_inclusion_3} $\lambda_{w_2,w_3}=\id$.
    \end{enumerate}
\end{lemma}
\begin{proof}
    It is clear that (\ref{lemma:path_inclusion_2}) implies (\ref{lemma:path_inclusion_1}).

    (\ref{lemma:path_inclusion_3}) $\Ra$ (\ref{lemma:path_inclusion_2}): $\lambda_{w_2,w_3}=\id$ is equivalent to assuming that both $w_2$ and $w_3$ are letters and $w_2 \prec w_3$. This requires $F(w_1) F(w_2+w_3)$ to be of the form $F(w_1)(F(w_2)F(w_3))$. Since all the morphisms $\lambda$'s start with moving parenthesis to the left, $\lambda_{w_1,w_2+w_3}$ starts exactly with $\mathsf{a}^{-1}$. Therefore $\lambda_{w_1,w_2+w_3} \circ \lambda_{w_2,w_3}= \lambda_{w_1+w_2,w_3} \circ \lambda_{w_1,w_2} \circ \mathsf{a}^{-1}$.

    (\ref{lemma:path_inclusion_1}) $\Ra$ (\ref{lemma:path_inclusion_3}): In all the vertices that $\lambda_{w_2,w_3}$ pass through, $F(w_1)$ is grouped separately from $F(w_2)$ and $F(w_3)$. Therefore any parenthesis movement or permutation that is part of $\lambda_{w_2,w_3}$ does not change the parenthesis around $F(w_1)$. However, on the path $\lambda_{w_1+w_2,w_3} \circ \lambda_{w_1,w_2}$ the same movements that describe $\lambda_{w_2,w_3}$ are part of the morphism $\lambda_{w_1+w_2,w_3}$. Since this path passes through the vertices that group $F(w_1)$ and $F(w_2)$, the parenthesis movements and permutations change the parenthesis around $F(w_1)$. This contradicts to the fact that $\lambda_{w_1+w_2,w_3} \circ \lambda_{w_1,w_2}$ is included in $\lambda_{w_1,w_2+w_3} \circ \lambda_{w_2,w_3} \circ \mathsf{a}$.
    \end{proof}

We remark that Lemma \ref{lemma:path_inclusion} can be also expressed as $\lambda_{w_1+w_2,w_3} \circ \lambda_{w_1,w_2}$ is strictly included in $\lambda_{w_1,w_2+w_3} \circ \lambda_{w_2,w_3} \circ \mathsf{a}$ if and only if $V_{(w_1|w_2,w_3)}= V_{(w_1,w_2|w_3)} \cup \{(F(w_1) (F(w_2)F(w_3)))\}$.

We can return to the definition of the 2-morphism $\psi_{w_1,w_2,w_3}$. By Lemmas \ref{lemma:path_equality} and \ref{lemma:path_inclusion}, the paths $\lambda_{w_1,w_2+w_3} \circ \lambda_{w_2,w_3} \circ \mathsf{a}$ and $\lambda_{w_1+w_2,w_3} \circ \lambda_{w_1,w_2}$ are going to satisfy one of the following three cases.
\begin{enumerate}
\item The paths may be the same. In this case, the 2-morphism $\psi_{w_1,w_2,w_3}$ is identity.
\item The path $\lambda_{w_1+w_2,w_3} \circ \lambda_{w_1,w_2}$ is strictly included in $\lambda_{w_1,w_2+w_3} \circ \lambda_{w_2,w_3} \circ \mathsf{a}$. In this case, by Lemma \ref{lemma:path_inclusion}, the 2-morphism $\psi_{w_1,w_2,w_3}$ is $\mathsf{a}\mathsf{a}^{-1} \Ra \id$.
\item The paths may enclose a 2-cell. This 2-cell is a tiling of pentagonal and rectangular 2-cells. The pentagonal 2-cells are either MacLane Pentagons or their derivatives obtained by inverting the direction of an edge. The rectangular 2-cells are of the form

\begin{equation}\label{diagram:natural_squares}
\begin{tabular}{ccc}
$\xymatrix{\bullet \ar[r]^{\mathsf{a_1}} \ar[d]_{\mathsf{a_2}} & \bullet \ar[d]^{\mathsf{a_2}}\\
          \bullet \ar[r]_{\mathsf{a_1}} & \bullet}$&
$\xymatrix{\bullet \ar[r]^{\mathsf{a_1}} \ar[d]_{\mathsf{c_1}} & \bullet \ar[d]^{\mathsf{c_1}}\\
          \bullet \ar[r]_{\mathsf{a_1}} & \bullet}$&
$\xymatrix{\bullet \ar[r]^{\mathsf{c_1}} \ar[d]_{\mathsf{c_2}} & \bullet \ar[d]^{\mathsf{c_2}}\\
          \bullet \ar[r]_{\mathsf{c_1}} & \bullet}$
\end{tabular}
\end{equation}
where $\mathsf{a_1}$, $\mathsf{a_2}$ are either leftward or rightward parenthesis movements and $\mathsf{c_1}$, $\mathsf{c_2}$ permute adjacent objects. Rectangular 2-cells can be also derived from (\ref{diagram:natural_squares}) by inverting the direction of an edge. These 2-cells commute up to structural 2-morphisms defined by the Picard structure of the 2-category $\twoC$. Theorem 3.3 in \cite{MR1040947} implies that these 2-morphisms compose in a unique way. We let $\psi_{w_1,w_2,w_3}$ be this composition.
\end{enumerate}

\paragraph{Definition of $\phi_{w_1,w_2}$:}
The last piece of the additive structure of $F$ is the 2-morphism $\phi_{w_1,w_2}$

\begin{equation}\label{diagram:2-morphism_of_braiding}
\xymatrix{F(w_1)F(w_2) \ar[d]_{\mathsf{c}} \ar[r]^(0.55){\lambda_{w_1,w_2}} \drtwocell<\omit>{\qquad\phi_{w_1,w_2}} & F(w_1+w_2)  \ar@{=}[d]\\
            F(w_2)F(w_1) \ar[r]_(0.55){\lambda_{w_2,w_1}} & F(w_2+w_1)}
\end{equation}
between the 1-morphisms $\lambda_{w_2,w_1} \circ \mathsf{c}$ and $\lambda_{w_1,w_2}$ from $F(w_1)F(w_2)$ to $F(w_1+w_2)$ where $w_1$ and $w_2$ are any two words in $\mathbb{N}(E)$. We notice that the path $\lambda_{w_2,w_1} \circ \mathsf{c}$ is not necessarily in the 1-skeleton of $K\Pi_{m+n-1}$. The reason is that the braiding $\mathsf{c}$ is not an adjacent permutation unless $w_1$ and $w_2$ are letters.

In the case where the words $w_1$ and $w_2$ are letters, $\phi_{w_1,w_2}$ is defined by the table

\begin{center}
\begin{tabular}{|l|c|c|}
\hline
\rule[-3mm]{0mm}{8mm} $\mathbf{w_1}$ & $\mathbf{w_2}$ & $\phi_{w_1,w_2}$\\ \hline
\rule[-3mm]{0mm}{8mm} a & a & $\id$ \\ \hline
\rule[-3mm]{0mm}{8mm} a & b & $\id \Ra \mathsf{c}^2$\\ \hline
\rule[-3mm]{0mm}{8mm} b & a & $\id$ \\ \hline
\end{tabular}
\end{center}
where $\id \Ra \mathsf{c}^2$ is given by the Picard structure of the 2-category $\twoC$.

Now, we assume that $w_1$ and $w_2$ are two words such that their sum of lengths is $m+n \geq 3$. The 2-morphism $\phi_{w_1,w_2}$ is defined in the following way. We first transform the path $\lambda_{w_2,w_1} \circ \mathsf{c}$ to a path in the 1-skeleton of $K\Pi_{m+n-1}$. Second we apply the process that defines $\psi_{w_1,w_2,w_3}$ to the new path and the path $\lambda_{w_1,w_2}$. $\phi_{w_1,w_2}$ is then defined as the appropriate composition of the 2-morphisms obtained at the first and the second step. Therefore to define $\phi_{w_1,w_2}$, it suffices to describe how we transform the path $\lambda_{w_2,w_1} \circ \mathsf{c}$ into a path in the 1-skeleton of $K\Pi_{m+n-1}$.

The main idea is to substitute the edge $\mathsf{c}$ that is not in the 1-skeleton by a sequence of five other edges. This sequence is an alternating collection of three leftward or rightward parenthesis movements and two braidings. The parenthesis movements are certainly in the 1-skeleton; however the braidings may not be. If they are not, then we substitute each of those braidings by a sequence of five other edges as above. We keep substituting until all the braidings become permutations of adjoint objects, therefore part of the 1-skeleton. We know that the substitution process is going to terminate because after each substitution braidings permute parenthesized objects with shorter length.

We describe this process on the sample $w_1=b+e$ and $w_2=a+c+d$. The braiding $\mathsf{c}$ permutes $F(w_1)$ and $F(w_2)$. First, we substitute $\mathsf{c}$ by the braidings $\mathsf{c}_{(a,c,d|e)}$ and $\mathsf{c}_{(a,c,d|b)}$. $\mathsf{c}_{(a,c,d|e)}$ permutes the parenthesized object $((FaFc)Fd)$ with $Fe$ and $\mathsf{c}_{(a,c,d|b)}$ permutes $((FaFc)Fd)$ with $Fb$. They are going to be substituted by $\mathsf{c}_{(d|e)}$ and $\mathsf{c}_{(a,c|e)}$ and by $\mathsf{c}_{(a,c|b)}$ and $\mathsf{c}_{(d|b)}$, respectively. Since $\mathsf{c}_{(d|e)}$ permutes $Fd$ and $Fe$ and $\mathsf{c}_{(d|b)}$ permutes $Fd$ and $Fb$, they are edges in the 1-skeleton and therefore cannot be substituted. In the diagram below, we illustrate the complete process of substituting $\mathsf{c}$ by adjacent permutations $\mathsf{c}_{(a|b)}$, $\mathsf{c}_{(c|b)}$, $\mathsf{c}_{(d|b)}$, $\mathsf{c}_{(a|e)}$, $\mathsf{c}_{(c|e)}$, and $\mathsf{c}_{(d|e)}$ using lattice paths.

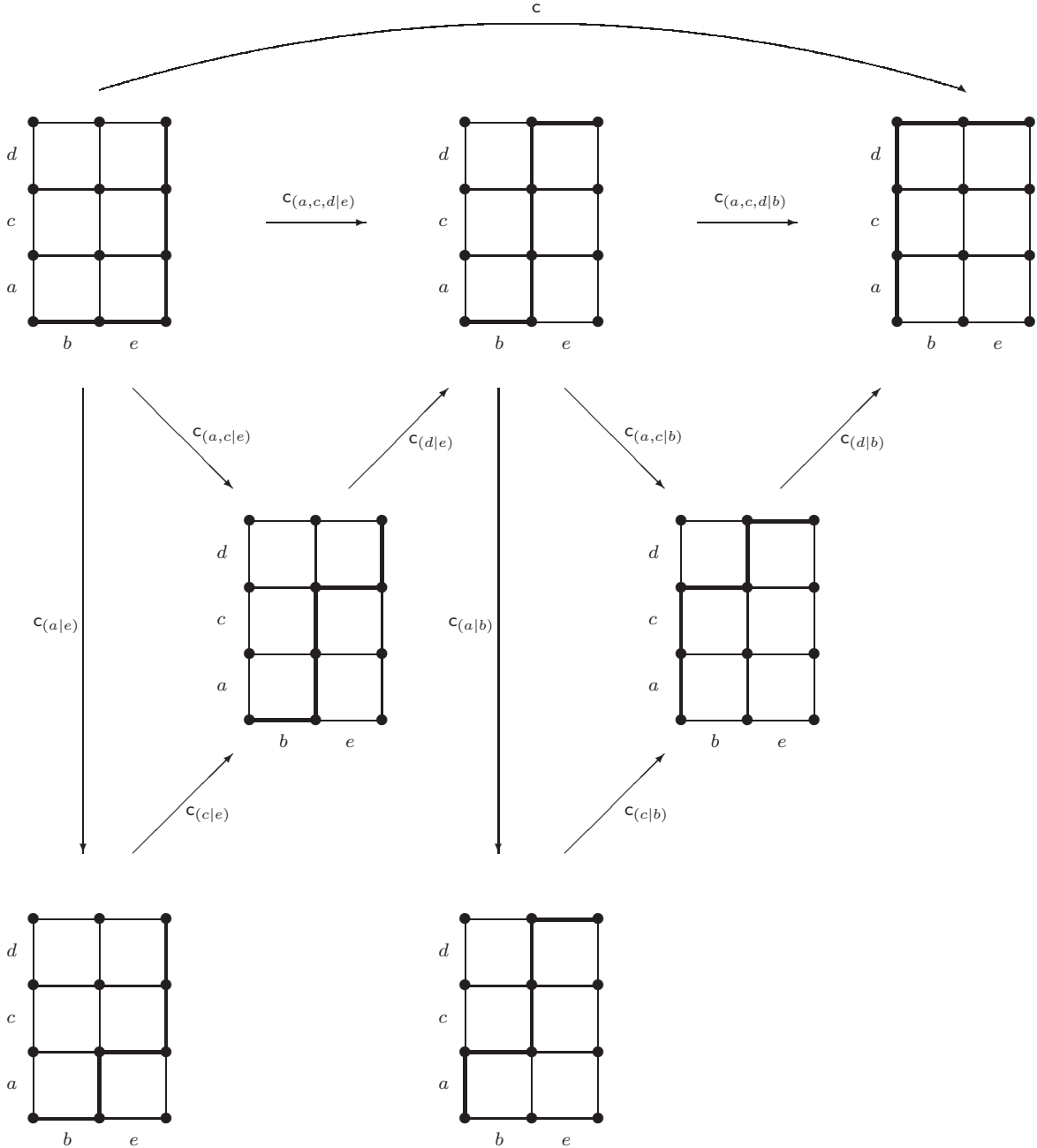
\unitlength 10mm
\begin{picture}(7,18)
\put(0.5,12.5){\line(1,0){2}}
\put(0.5,13.5){\line(1,0){2}}
\put(0.5,14.5){\line(1,0){2}}
\put(0.5,15.5){\line(1,0){2}}
\put(0.5,12.5){\line(0,1){3}}
\put(1.5,12.5){\line(0,1){3}}
\put(2.5,12.5){\line(0,1){3}}
\put(0.5,12.5){\makebox(0,0)[cc]{$\bullet$}}
\put(1.5,12.5){\makebox(0,0)[cc]{$\bullet$}}
\put(2.5,12.5){\makebox(0,0)[cc]{$\bullet$}}
\put(0.5,13.5){\makebox(0,0)[cc]{$\bullet$}}
\put(1.5,13.5){\makebox(0,0)[cc]{$\bullet$}}
\put(2.5,13.5){\makebox(0,0)[cc]{$\bullet$}}
\put(0.5,14.5){\makebox(0,0)[cc]{$\bullet$}}
\put(1.5,14.5){\makebox(0,0)[cc]{$\bullet$}}
\put(2.5,14.5){\makebox(0,0)[cc]{$\bullet$}}
\put(0.5,15.5){\makebox(0,0)[cc]{$\bullet$}}
\put(1.5,15.5){\makebox(0,0)[cc]{$\bullet$}}
\put(2.5,15.5){\makebox(0,0)[cc]{$\bullet$}}
\qbezier(1.5,16)(8,18)(14.5,16)
\put(14.5,16){\vector(4,-3){0.05}}
\put(8,17.3){\makebox(0,0)[tl]{\scriptsize$\mathsf{c}$}}
\put(0.95,12.1){\makebox(0,0)[lb]{\scriptsize$b$}}
\put(1.95,12.1){\makebox(0,0)[lb]{\scriptsize$e$}}
\put(0.1,12.95){\makebox(0,0)[lb]{\scriptsize$a$}}
\put(0.1,13.95){\makebox(0,0)[lb]{\scriptsize$c$}}
\put(0.1,14.95){\makebox(0,0)[lb]{\scriptsize$d$}}
\put(4,14){\vector(1,0){1.5}}
\put(4.25,14.2){\makebox(0,0)[lb]{\scriptsize$\mathsf{c}_{(a,c,d|e)}$}}
\put(1.25,11.5){\vector(0,-1){7}}
\put(0.5,8.05){\makebox(0,0)[tl]{\scriptsize$\mathsf{c}_{(a|e)}$}}
\put(2,11.5){\vector(1,-1){1.5}}
\put(2.9,10.9){\makebox(0,0)[tl]{\scriptsize$\mathsf{c}_{(a,c|e)}$}}
\put(7,12.5){\line(1,0){2}}
\put(7,13.5){\line(1,0){2}}
\put(7,14.5){\line(1,0){2}}
\put(7,15.5){\line(1,0){2}}
\put(7,12.5){\line(0,1){3}}
\put(8,12.5){\line(0,1){3}}
\put(9,12.5){\line(0,1){3}}
\put(7,12.5){\makebox(0,0)[cc]{$\bullet$}}
\put(8,12.5){\makebox(0,0)[cc]{$\bullet$}}
\put(9,12.5){\makebox(0,0)[cc]{$\bullet$}}
\put(7,13.5){\makebox(0,0)[cc]{$\bullet$}}
\put(8,13.5){\makebox(0,0)[cc]{$\bullet$}}
\put(9,13.5){\makebox(0,0)[cc]{$\bullet$}}
\put(7,14.5){\makebox(0,0)[cc]{$\bullet$}}
\put(8,14.5){\makebox(0,0)[cc]{$\bullet$}}
\put(9,14.5){\makebox(0,0)[cc]{$\bullet$}}
\put(7,15.5){\makebox(0,0)[cc]{$\bullet$}}
\put(8,15.5){\makebox(0,0)[cc]{$\bullet$}}
\put(9,15.5){\makebox(0,0)[cc]{$\bullet$}}
\put(7.45,12.1){\makebox(0,0)[lb]{\scriptsize$b$}}
\put(8.45,12.1){\makebox(0,0)[lb]{\scriptsize$e$}}
\put(6.6,12.95){\makebox(0,0)[lb]{\scriptsize$a$}}
\put(6.6,13.95){\makebox(0,0)[lb]{\scriptsize$c$}}
\put(6.6,14.95){\makebox(0,0)[lb]{\scriptsize$d$}}
\put(10.5,14){\vector(1,0){1.5}}
\put(10.75,14.2){\makebox(0,0)[lb]{\scriptsize$\mathsf{c}_{(a,c,d|b)}$}}
\put(7.5,11.5){\vector(0,-1){7}}
\put(6.75,8.05){\makebox(0,0)[tl]{\scriptsize$\mathsf{c}_{(a|b)}$}}
\put(8.5,11.5){\vector(1,-1){1.5}}
\put(9.4,10.9){\makebox(0,0)[tl]{\scriptsize$\mathsf{c}_{(a,c|b)}$}}
\put(5.25,10){\vector(1,1){1.5}}
\put(6.15,10.8){\makebox(0,0)[tl]{\scriptsize$\mathsf{c}_{(d|e)}$}}
\put(11.75,10){\vector(1,1){1.5}}
\put(12.65,10.8){\makebox(0,0)[tl]{\scriptsize$\mathsf{c}_{(d|b)}$}}
\put(13.5,12.5){\line(1,0){2}}
\put(13.5,13.5){\line(1,0){2}}
\put(13.5,14.5){\line(1,0){2}}
\put(13.5,15.5){\line(1,0){2}}
\put(13.5,12.5){\line(0,1){3}}
\put(14.5,12.5){\line(0,1){3}}
\put(15.5,12.5){\line(0,1){3}}
\put(13.5,12.5){\makebox(0,0)[cc]{$\bullet$}}
\put(14.5,12.5){\makebox(0,0)[cc]{$\bullet$}}
\put(15.5,12.5){\makebox(0,0)[cc]{$\bullet$}}
\put(13.5,13.5){\makebox(0,0)[cc]{$\bullet$}}
\put(14.5,13.5){\makebox(0,0)[cc]{$\bullet$}}
\put(15.5,13.5){\makebox(0,0)[cc]{$\bullet$}}
\put(13.5,14.5){\makebox(0,0)[cc]{$\bullet$}}
\put(14.5,14.5){\makebox(0,0)[cc]{$\bullet$}}
\put(15.5,14.5){\makebox(0,0)[cc]{$\bullet$}}
\put(13.5,15.5){\makebox(0,0)[cc]{$\bullet$}}
\put(14.5,15.5){\makebox(0,0)[cc]{$\bullet$}}
\put(15.5,15.5){\makebox(0,0)[cc]{$\bullet$}}
\put(13.95,12.1){\makebox(0,0)[lb]{\scriptsize$b$}}
\put(14.95,12.1){\makebox(0,0)[lb]{\scriptsize$e$}}
\put(13.1,12.95){\makebox(0,0)[lb]{\scriptsize$a$}}
\put(13.1,13.95){\makebox(0,0)[lb]{\scriptsize$c$}}
\put(13.1,14.95){\makebox(0,0)[lb]{\scriptsize$d$}}
\linethickness{.5mm}
\put(0.5,12.5){\line(1,0){1}}
\put(1.5,12.5){\line(1,0){1}}
\put(2.5,12.5){\line(0,1){1}}
\put(2.5,13.5){\line(0,1){1}}
\put(2.5,14.5){\line(0,1){1}}
\put(7,12.5){\line(1,0){1}}
\put(8,12.5){\line(0,1){1}}
\put(8,13.5){\line(0,1){1}}
\put(8,14.5){\line(0,1){1}}
\put(8,15.5){\line(1,0){1}}
\put(13.5,12.5){\line(0,1){1}}
\put(13.5,13.5){\line(0,1){1}}
\put(13.5,14.5){\line(0,1){1}}
\put(13.5,15.5){\line(1,0){1}}
\put(14.5,15.5){\line(1,0){1}}
\thinlines
\put(3.75,6.5){\line(1,0){2}}
\put(3.75,7.5){\line(1,0){2}}
\put(3.75,8.5){\line(1,0){2}}
\put(3.75,9.5){\line(1,0){2}}
\put(3.75,6.5){\line(0,1){3}}
\put(4.75,6.5){\line(0,1){3}}
\put(5.75,6.5){\line(0,1){3}}
\put(3.75,6.5){\makebox(0,0)[cc]{$\bullet$}}
\put(4.75,6.5){\makebox(0,0)[cc]{$\bullet$}}
\put(5.75,6.5){\makebox(0,0)[cc]{$\bullet$}}
\put(3.75,7.5){\makebox(0,0)[cc]{$\bullet$}}
\put(4.75,7.5){\makebox(0,0)[cc]{$\bullet$}}
\put(5.75,7.5){\makebox(0,0)[cc]{$\bullet$}}
\put(3.75,8.5){\makebox(0,0)[cc]{$\bullet$}}
\put(4.75,8.5){\makebox(0,0)[cc]{$\bullet$}}
\put(5.75,8.5){\makebox(0,0)[cc]{$\bullet$}}
\put(3.75,9.5){\makebox(0,0)[cc]{$\bullet$}}
\put(4.75,9.5){\makebox(0,0)[cc]{$\bullet$}}
\put(5.75,9.5){\makebox(0,0)[cc]{$\bullet$}}
\put(4.2,6.1){\makebox(0,0)[lb]{\scriptsize$b$}}
\put(5.2,6.1){\makebox(0,0)[lb]{\scriptsize$e$}}
\put(3.26,6.95){\makebox(0,0)[lb]{\scriptsize$a$}}
\put(3.26,7.95){\makebox(0,0)[lb]{\scriptsize$c$}}
\put(3.26,8.95){\makebox(0,0)[lb]{\scriptsize$d$}}
\put(10.25,6.5){\line(1,0){2}}
\put(10.25,7.5){\line(1,0){2}}
\put(10.25,8.5){\line(1,0){2}}
\put(10.25,9.5){\line(1,0){2}}
\put(10.25,6.5){\line(0,1){3}}
\put(11.25,6.5){\line(0,1){3}}
\put(12.25,6.5){\line(0,1){3}}
\put(10.25,6.5){\makebox(0,0)[cc]{$\bullet$}}
\put(11.25,6.5){\makebox(0,0)[cc]{$\bullet$}}
\put(12.25,6.5){\makebox(0,0)[cc]{$\bullet$}}
\put(10.25,7.5){\makebox(0,0)[cc]{$\bullet$}}
\put(11.25,7.5){\makebox(0,0)[cc]{$\bullet$}}
\put(12.25,7.5){\makebox(0,0)[cc]{$\bullet$}}
\put(10.25,8.5){\makebox(0,0)[cc]{$\bullet$}}
\put(11.25,8.5){\makebox(0,0)[cc]{$\bullet$}}
\put(12.25,8.5){\makebox(0,0)[cc]{$\bullet$}}
\put(10.25,9.5){\makebox(0,0)[cc]{$\bullet$}}
\put(11.25,9.5){\makebox(0,0)[cc]{$\bullet$}}
\put(12.25,9.5){\makebox(0,0)[cc]{$\bullet$}}
\put(10.7,6.1){\makebox(0,0)[lb]{\scriptsize$b$}}
\put(11.7,6.1){\makebox(0,0)[lb]{\scriptsize$e$}}
\put(9.76,6.95){\makebox(0,0)[lb]{\scriptsize$a$}}
\put(9.76,7.95){\makebox(0,0)[lb]{\scriptsize$c$}}
\put(9.76,8.95){\makebox(0,0)[lb]{\scriptsize$d$}}
\linethickness{.5mm}
\put(3.75,6.5){\line(1,0){1}}
\put(4.75,6.5){\line(0,1){1}}
\put(4.75,7.5){\line(0,1){1}}
\put(4.75,8.5){\line(1,0){1}}
\put(5.75,8.5){\line(0,1){1}}
\put(10.25,6.5){\line(0,1){1}}
\put(10.25,7.5){\line(0,1){1}}
\put(10.25,8.5){\line(1,0){1}}
\put(11.25,8.5){\line(0,1){1}}
\put(11.25,9.5){\line(1,0){1}}
\thinlines
\put(0.5,0.5){\line(1,0){2}}
\put(0.5,1.5){\line(1,0){2}}
\put(0.5,2.5){\line(1,0){2}}
\put(0.5,3.5){\line(1,0){2}}
\put(0.5,0.5){\line(0,1){3}}
\put(1.5,0.5){\line(0,1){3}}
\put(2.5,0.5){\line(0,1){3}}
\put(0.5,0.5){\makebox(0,0)[cc]{$\bullet$}}
\put(1.5,0.5){\makebox(0,0)[cc]{$\bullet$}}
\put(2.5,0.5){\makebox(0,0)[cc]{$\bullet$}}
\put(0.5,1.5){\makebox(0,0)[cc]{$\bullet$}}
\put(1.5,1.5){\makebox(0,0)[cc]{$\bullet$}}
\put(2.5,1.5){\makebox(0,0)[cc]{$\bullet$}}
\put(0.5,2.5){\makebox(0,0)[cc]{$\bullet$}}
\put(1.5,2.5){\makebox(0,0)[cc]{$\bullet$}}
\put(2.5,2.5){\makebox(0,0)[cc]{$\bullet$}}
\put(0.5,3.5){\makebox(0,0)[cc]{$\bullet$}}
\put(1.5,3.5){\makebox(0,0)[cc]{$\bullet$}}
\put(2.5,3.5){\makebox(0,0)[cc]{$\bullet$}}
\put(7.45,0.1){\makebox(0,0)[lb]{\scriptsize$b$}}
\put(8.45,0.1){\makebox(0,0)[lb]{\scriptsize$e$}}
\put(6.6,0.95){\makebox(0,0)[lb]{\scriptsize$a$}}
\put(6.6,1.95){\makebox(0,0)[lb]{\scriptsize$c$}}
\put(6.6,2.95){\makebox(0,0)[lb]{\scriptsize$d$}}
\put(2,4.5){\vector(1,1){1.5}}
\put(2.8,5.2){\makebox(0,0)[tl]{\scriptsize$\mathsf{c}_{(c|e)}$}}
\put(7,0.5){\line(1,0){2}}
\put(7,1.5){\line(1,0){2}}
\put(7,2.5){\line(1,0){2}}
\put(7,3.5){\line(1,0){2}}
\put(7,0.5){\line(0,1){3}}
\put(8,0.5){\line(0,1){3}}
\put(9,0.5){\line(0,1){3}}
\put(7,0.5){\makebox(0,0)[cc]{$\bullet$}}
\put(8,0.5){\makebox(0,0)[cc]{$\bullet$}}
\put(9,0.5){\makebox(0,0)[cc]{$\bullet$}}
\put(7,1.5){\makebox(0,0)[cc]{$\bullet$}}
\put(8,1.5){\makebox(0,0)[cc]{$\bullet$}}
\put(9,1.5){\makebox(0,0)[cc]{$\bullet$}}
\put(7,2.5){\makebox(0,0)[cc]{$\bullet$}}
\put(8,2.5){\makebox(0,0)[cc]{$\bullet$}}
\put(9,2.5){\makebox(0,0)[cc]{$\bullet$}}
\put(7,3.5){\makebox(0,0)[cc]{$\bullet$}}
\put(8,3.5){\makebox(0,0)[cc]{$\bullet$}}
\put(9,3.5){\makebox(0,0)[cc]{$\bullet$}}
\put(0.95,0.1){\makebox(0,0)[lb]{\scriptsize$b$}}
\put(1.95,0.1){\makebox(0,0)[lb]{\scriptsize$e$}}
\put(0.1,0.95){\makebox(0,0)[lb]{\scriptsize$a$}}
\put(0.1,1.95){\makebox(0,0)[lb]{\scriptsize$c$}}
\put(0.1,2.95){\makebox(0,0)[lb]{\scriptsize$d$}}
\put(8.5,4.5){\vector(1,1){1.5}}
\put(9.4,5.2){\makebox(0,0)[tl]{\scriptsize$\mathsf{c}_{(c|b)}$}}
\linethickness{.5mm}
\put(0.5,0.5){\line(1,0){1}}
\put(1.5,0.5){\line(0,1){1}}
\put(1.5,1.5){\line(1,0){1}}
\put(2.5,1.5){\line(0,1){1}}
\put(2.5,2.5){\line(0,1){1}}
\put(7,0.5){\line(0,1){1}}
\put(7,1.5){\line(1,0){1}}
\put(8,1.5){\line(0,1){1}}
\put(8,2.5){\line(0,1){1}}
\put(8,3.5){\line(1,0){1}}
\end{picture}

This process defines a 2-morphism as follows. Substituting a braiding by an alternating sequence of three leftward or rightward parenthesis movements and two braidings means substituting an edge in a hexagonal 2-cell by the other five edges. Such hexagonal 2-cells commute up to a 2-morphism given by the Picard structure of the 2-category $\twoC$. The appropriate composition of these 2-morphisms defines the 2-morphism of the first step.

\subsection{Extending the Additive Structure to Free Abelian Group}\label{subsection:group_case}
Here we extend the definition of the 2-functor $F$ so that it transforms the trivial Picard structure of the free abelian group $\mathbb{Z}(E)$ generated by the set $E$ to the Picard structure of the 2-category $\twoC$.
\paragraph{Extending $\lambda_{w_1,w_2}$:} The extension of $\lambda_{w_1,w_2}$, denoted by $\widetilde{\lambda}_{w_1,w_2}$, to the words in $\mathbb{Z}(E)$ should take into consideration the cancelations that might occur in $w_1+w_2$. If $w_2$ does not have a letter that appears with an opposite sign in $w_1$ then there aren't any cancelations in $w_1+w_2$ and $\widetilde{\lambda}_{w_1,w_2}=\lambda_{w_1,w_2}$. Otherwise, $\widetilde{\lambda}_{w_1,w_2}$ orders the letters of $w_1$ and $w_2$ from least to greatest, left parenthesizes, and does the cancelations starting with the image of the least letter.  That is $\widetilde{\lambda}_{w_1,w_2}$ is equal to post composition of $\lambda_{w_1,w_2}$ with the morphisms of the form
\begin{equation}\label{diagram:tail_of_lambda}
\xymatrix{&(\ldots(((F(w)Fc_i){(Fc_i)}^*)Fc_{i+1})\ldots Fc_{n+m}) \ar[r] & (\ldots((F(w)(Fc_i{(Fc_i)}^*))Fc_{i+1})\ldots Fc_{n+m}) \ar`[d]`[dll]`[ddll]_{\textrm{inv}_{Fc_i}}`[ddllr][ddl]\\&&\\
          &(\ldots((F(w)I)Fc_{i+1})\ldots Fc_{n+m}) \ar[r]_{\mathrm{r}_{F(w)}} & (\ldots(F(w)Fc_{i+1})\ldots Fc_{n+m})}
\end{equation}
for every cancelation. In (\ref{diagram:tail_of_lambda}) $w$ is a subword of $w_1+w_2$, $I$ is a unit element in the Picard 2-category and $\textrm{inv}_{Fc_i}$ and $\mathrm{r}_{F(w)}$ are structural morphisms due to the Picard structure of the 2-category. By the Picard structure, we can also assume for simplicity that when $\widetilde{\lambda}_{w_1,w_2}$ orders letters from least to greatest the inverse of an object is always adjacent to the object and it is on its left. We note that using $\lambda_{w_1,w_2}$ for the morphism that orders the letters of $w_1$ and $w_2$ from least to greatest and left parenthesizes them is an abuse of notation. Here $\lambda_{w_1,w_2}$ does not map to the object $F(w_1+w_2)$ but to an object that we denote $F(w_{1,2})$. $F(w_{1,2})$ is product of the images of all letters in $w_1$ and $w_2$ parenthesized from the left, ordered from least to greatest, and if there exists inverse of an object is placed on its left. For instance, if $w_1=b+c$ and $w_2=a-b$, then
\[\lambda_{w_1,w_2}:\xymatrix@1{(Fb Fc)(Fa(Fb)^*) \ar[r] & (((FaFb)Fb)^*)Fc),}\]
where $F(w_{1,2})=((FaFb)(Fb)^*)Fc)$. Thus $\widetilde{\lambda}_{w_1,w_2}$ can be expressed as composition of
\begin{equation}\label{diagram:lambda_extended}
\xymatrix@1{F(w_1)F(w_2) \ar[rr]^(0.55){\lambda_{w_1,w_2}} && F(w_{1,2}) \ar[rr]^(0.45){\tau_{w_1,w_2}} && F(w_1+w_2),}
\end{equation}
where $\tau_{w_1,w_2}$ is composition of morphisms of the form (\ref{diagram:tail_of_lambda}) for every cancelation. We remark that $\lambda_{w_1,w_2}$ as in the monoidal case defines a path in the 1-skeleton of the permuto-associahedron $K\Pi_{m+n-1}$. However if there are cancelations, $\widetilde{\lambda}_{w_1,w_2}$ is not a path in the 1-skeleton of $K\Pi_{m+n-1}$.

\paragraph{Extending $\psi_{w_1,w_2,w_3}$:} The extension of $\psi_{w_1,w_2,w_3}$, denoted by $\widetilde{\psi}_{w_1,w_2,w_3}$, to the words $w_1$, $w_2$, $w_3$ in $\mathbb{Z}(E)$ is a 2-morphism

\begin{equation}\label{diagram:2-morphism_of_assoc_extended}
\xymatrix{((F(w_1)F(w_2)) F(w_3)) \ar[d]_{\mathsf{a}} \ar[r]^(0.55){\widetilde{\lambda}_{w_1,w_2}} & F(w_1+w_2) F(w_3) \ar[r]^(0.65){\widetilde{\lambda}_{w_1+w_2,w_3}} \drtwocell<\omit>{\hspace{1.2cm}\widetilde{\psi}_{w_1,w_2,w_3}}& **[r]F(w_1+w_2+w_3) \ar@{=}[d]\\
            (F(w_1)(F(w_2)F(w_3))) \ar[r]_(0.55){\widetilde{\lambda}_{w_2,w_3}} & F(w_1) F(w_2+w_3) \ar[r]_(0.65){\widetilde{\lambda}_{w_1,w_2+w_3}} & **[r] F(w_1+w_2+w_3)}
\end{equation}
between the 1-morphisms $\widetilde{\lambda}_{w_1,w_2+w_3} \circ \widetilde{\lambda}_{w_2,w_3} \circ \mathsf{a}$ and $\widetilde{\lambda}_{w_1+w_2,w_3} \circ \widetilde{\lambda}_{w_1,w_2}$. As noticed, these paths may not be in the 1-skeleton of $K\Pi_{m+n+p-1}$. However, there exists a vertex $V_0$ of the permuto-associahedron $K\Pi_{m+n+p-1}$ that both paths $\widetilde{\lambda}_{w_1+w_2,w_3} \circ \widetilde{\lambda}_{w_1,w_2}$ and $\widetilde{\lambda}_{w_1,w_2+w_3} \circ \widetilde{\lambda}_{w_2,w_3}$ pass through. Therefore the diagram (\ref{diagram:2-morphism_of_assoc_extended}) can be rewritten as:

\begin{equation}\label{diagram:2-morphism_of_assoc_extended_1}
\xymatrix{((F(w_1)F(w_2)) F(w_3)) \ar[d]_{\mathsf{a}} \ar[r] \drtwocell<\omit>{\hspace{1.2cm}\psi'_{w_1,w_2,w_3}}& V_0 \ar[r] \ar@{=}[d] & F(w_1+w_2) F(w_3) \ar[r]^(0.65){\widetilde{\lambda}_{w_1+w_2,w_3}} \drtwocell<\omit>{\mspace{80mu}\rho_{w_1,w_2,w_3}}& **[r]F(w_1+w_2+w_3) \ar@{=}[d]\\
            (F(w_1)(F(w_2)F(w_3))) \ar[r]& V_0 \ar[r] & F(w_1) F(w_2+w_3) \ar[r]_(0.65){\widetilde{\lambda}_{w_1,w_2+w_3}} & **[r] F(w_1+w_2+w_3)}
\end{equation}
where both horizontal morphisms to $V_0$ are paths on $K\Pi_{m+n+p-1}$. So we compute $\psi'_{w_1,w_2,w_3}$ in the same way as $\psi$ of the monoidal case. After the vertex $V_0$, the morphisms on the diagram (\ref{diagram:2-morphism_of_assoc_extended_1}) are not any more in the 1-skeleton of $K\Pi_{m+n+p-1}$ because of the cancelations. The region between the two paths from $V_0$ to $F(w_1+w_2+w_3)$ can be filled with the structural 2-morphisms of the Picard structure in particular by the ones involving the inverse and unit objects. The 2-morphism $\rho_{w_1,w_2,w_3}$ is then the unique pasting of those structural 2-morphisms. Hence, we define $\widetilde{\psi}_{w_1,w_2,w_3}$ as pasting of $\psi'_{w_1,w_2,w_3}$ and $\rho_{w_1,w_2,w_3}$.

\paragraph{Extending $\phi_{w_1,w_2}$:} The extension of $\phi_{w_1,w_2}$, denoted by $\widetilde{\phi}_{w_1,w_2}$ is a 2-morphism

\begin{equation}\label{diagram:2-morphism_of_braiding_extended}
\xymatrix{F(w_1)F(w_2) \ar[d]_{\mathsf{c}} \ar[r]^(0.55){\widetilde{\lambda}_{w_1,w_2}} \drtwocell<\omit>{\qquad\widetilde{\phi}_{w_1,w_2}} & F(w_1+w_2)  \ar@{=}[d]\\
            F(w_2)F(w_1) \ar[r]_(0.55){\widetilde{\lambda}_{w_2,w_1}} & F(w_2+w_1)}
\end{equation}
between the 1-morphisms $\widetilde{\lambda}_{w_2,w_1} \circ \mathsf{c}$ and $\widetilde{\lambda}_{w_1,w_2}$ from $F(w_1)F(w_2)$ to $F(w_1+w_2)$ where $w_1$ and $w_2$ are any two words in $\mathbb{Z}(E)$. We rewrite the diagram (\ref{diagram:2-morphism_of_braiding_extended}) by expressing $\widetilde{\lambda}_{w_1,w_2}$ and $\widetilde{\lambda}_{w_2,w_1}$ as compositions using (\ref{diagram:lambda_extended}).

\begin{equation}\label{diagram:2-morphism_of_braiding_extended_1}
\xymatrix{F(w_1)F(w_2) \ar[d]_{\mathsf{c}} \ar[r]^(0.6){{\lambda}_{w_1,w_2}} \drtwocell<\omit>{\qquad \phi'_{w_1,w_2}} & F(w_{1,2})  \ar@{=}[d] \ar[r]^{\tau_{w_1,w_2}} & F(w_1+w_2) \ar@{=}[d]\\
            F(w_2)F(w_1) \ar[r]_(0.6){{\lambda}_{w_2,w_1}} & F(w_{2,1}) \ar[r]_{\tau_{w_2,w_1}} \ar[r] & F(w_2+w_1)}
\end{equation}
The square on the left commutes up to the 2-morphism $\phi'_{w_1,w_2}$ defined in the same way as $\phi$ of the monoidal case. The square on the right commutes since $F(w_{1,2})=F(w_{2,1})$ and therefore $\tau_{w_1,w_2}=\tau_{w_2,w_1}$. Hence, $\widetilde{\phi}_{w_1,w_2}$ is the whiskering $\phi'_{w_1,w_2} * \tau_{w_1,w_2}$.
\bibliographystyle{plain}

\begin{thebibliography}{10}

\bibitem{MR2387582}
Ettore Aldrovandi.
\newblock 2-gerbes bound by complexes of {$gr$}-stacks, and cohomology.
\newblock {\em J. Pure Appl. Algebra}, 212(5):994--1038, 2008.

\bibitem{Aldrovandi2009687-III}
Ettore Aldrovandi and Behrang Noohi.
\newblock Butterflies {III}: Higher butterflies and higher gr-stacks.
\newblock {\em In preparation}.

\bibitem{Aldrovandi2009687}
Ettore Aldrovandi and Behrang Noohi.
\newblock Butterflies {I}: Morphisms of 2-group stacks.
\newblock {\em Advances in Mathematics}, 221(3):687 -- 773, 2009.

\bibitem{MR1402727}
John~C. Baez and Martin Neuchl.
\newblock Higher-dimensional algebra. {I}. {B}raided monoidal {$2$}-categories.
\newblock {\em Adv. Math.}, 121(2):196--244, 1996.

\bibitem{MR0220789}
Jean B{\'e}nabou.
\newblock Introduction to bicategories.
\newblock In {\em Reports of the {M}idwest {C}ategory {S}eminar}, pages 1--77.
  Springer, Berlin, 1967.

\bibitem{MR1086889}
Lawrence Breen.
\newblock Bitorseurs et cohomologie non ab\'elienne.
\newblock In {\em The {G}rothendieck {F}estschrift, {V}ol.\ {I}}, volume~86 of
  {\em Progr. Math.}, pages 401--476. Birkh\"auser Boston, Boston, MA, 1990.

\bibitem{MR1191733}
Lawrence Breen.
\newblock Th\'eorie de {S}chreier sup\'erieure.
\newblock {\em Ann. Sci. \'Ecole Norm. Sup. (4)}, 25(5):465--514, 1992.

\bibitem{MR1301844}
Lawrence Breen.
\newblock On the classification of {$2$}-gerbes and {$2$}-stacks.
\newblock {\em Ast\'erisque}, (225):160, 1994.

\bibitem{breen-2006}
Lawrence Breen.
\newblock Notes on 1- and 2-gerbes.
\newblock In {\em Towards Higher Categories, J.C. Baez and J.P. May (eds.)},
  volume 152 of {\em The IMA Volumes in Mathematics and its Applications},
  pages 193--235. Springer, New York, 2010.

\bibitem{MR1458415}
Brian Day and Ross Street.
\newblock Monoidal bicategories and {H}opf algebroids.
\newblock {\em Adv. Math.}, 129(1):99--157, 1997.

\bibitem{Deligne}
Pierre Deligne.
\newblock La formule de dualit\'e globale, 1973.
\newblock SGA 4 III, Expos\'e XVIII.

\bibitem{MR546620}
Pierre Deligne.
\newblock Vari\'et\'es de {S}himura: interpr\'etation modulaire, et techniques
  de construction de mod\`eles canoniques.
\newblock In {\em Automorphic forms, representations and {$L$}-functions
  ({P}roc. {S}ympos. {P}ure {M}ath., {O}regon {S}tate {U}niv., {C}orvallis,
  {O}re., 1977), {P}art 2}, Proc. Sympos. Pure Math., XXXIII, pages 247--289.
  Amer. Math. Soc., Providence, R.I., 1979.

\bibitem{MR2394437}
I.~M. Gelfand, M.~M. Kapranov, and A.~V. Zelevinsky.
\newblock {\em Discriminants, resultants and multidimensional determinants}.
\newblock Modern Birkh\"auser Classics. Birkh\"auser Boston Inc., Boston, MA,
  2008.
\newblock Reprint of the 1994 edition.

\bibitem{MR1261589}
R.~Gordon, A.~J. Power, and Ross Street.
\newblock Coherence for tricategories.
\newblock {\em Mem. Amer. Math. Soc.}, 117(558):vi+81, 1995.

\bibitem{Gurski_thesis}
Nick Gurski.
\newblock An algebraic theory of tricategories.
\newblock {\em PhD Thesis}, 2006.

\bibitem{MR0364245}
Monique Hakim.
\newblock {\em Topos annel\'es et sch\'emas relatifs}.
\newblock Springer-Verlag, Berlin, 1972.
\newblock Ergebnisse der Mathematik und ihrer Grenzgebiete, Band 64.

\bibitem{hirschowitz-1998}
Andre Hirschowitz and Carlos Simpson.
\newblock Descente pour les n-champs (descent for n-stacks), 1998.

\bibitem{MR1278735}
M.~M. Kapranov and V.~A. Voevodsky.
\newblock {$2$}-categories and {Z}amolodchikov tetrahedra equations.
\newblock In {\em Algebraic groups and their generalizations: quantum and
  infinite-dimensional methods ({U}niversity {P}ark, {PA}, 1991)}, volume~56 of
  {\em Proc. Sympos. Pure Math.}, pages 177--259. Amer. Math. Soc., Providence,
  RI, 1994.

\bibitem{MR1207505}
Mikhail~M. Kapranov.
\newblock The permutoassociahedron, {M}ac {L}ane's coherence theorem and
  asymptotic zones for the {KZ} equation.
\newblock {\em J. Pure Appl. Algebra}, 85(2):119--142, 1993.

\bibitem{MR2276246}
Stephen Lack.
\newblock Bicat is not triequivalent to {G}ray.
\newblock {\em Theory Appl. Categ.}, 18:No. 1, 1--3 (electronic), 2007.

\bibitem{noohi-2005}
Behrang Noohi.
\newblock On weak maps between 2-groups, 2005.

\bibitem{MR2280287}
Behrang Noohi.
\newblock Notes on 2-groupoids, 2-groups and crossed modules.
\newblock {\em Homology, Homotopy Appl.}, 9(1):75--106 (electronic), 2007.

\bibitem{MR1040947}
A.~J. Power.
\newblock A {$2$}-categorical pasting theorem.
\newblock {\em J. Algebra}, 129(2):439--445, 1990.

\bibitem{simpson-1997}
Carlos Simpson.
\newblock A closed model structure for $n$-categories, internal $hom$,
  $n$-stacks and generalized seifert-van kampen, 1997.

\bibitem{MR1673923}
Zouhair Tamsamani.
\newblock Sur des notions de {$n$}-cat\'egorie et {$n$}-groupo\"\i de non
  strictes via des ensembles multi-simpliciaux.
\newblock {\em $K$-Theory}, 16(1):51--99, 1999.

\bibitem{Tatar_thesis}
Ahmet~E. Tatar.
\newblock On the picard 2-stacks.
\newblock {\em PhD Thesis, In preparation}.

\bibitem{MR1311028}
G{\"u}nter~M. Ziegler.
\newblock {\em Lectures on polytopes}, volume 152 of {\em Graduate Texts in
  Mathematics}.
\newblock Springer-Verlag, New York, 1995.

\end{thebibliography}

\end{document}